\documentclass[preprint,sort&compress,number]{elsarticle}

\usepackage{amsmath,amsthm,amssymb,latexsym,amscd}
\usepackage{subfigure} 
\usepackage{booktabs}
\usepackage[list-final-separator={, and }]{siunitx}
\usepackage[shadow,textsize=small,textwidth=2.5cm]{todonotes}
\usepackage{stmaryrd}

\usepackage{geometry}
\usepackage{mathtools}
\usepackage{algorithm}
\usepackage{algpseudocode}
\algtext*{EndIf}
\algtext*{EndWhile}
\algtext*{EndFor}

\usepackage{color}
\usepackage{nomencl}
\usepackage{bm}
\usepackage{lineno}
\usepackage{cancel}

\usepackage{amsmath,amsthm,amssymb,latexsym,amscd}
\usepackage{lineno}
\usepackage[defaultcolor=red]{changes}

\usepackage{todonotes}

\newcommand{\ATone}{\texttt{AT$_1$}}
\newcommand{\ATtwo}{\texttt{AT$_2$}}

\newcommand{\mrm}{\mathrm}
\newcommand{\Gc}{G_\mrm{c}}

\newcommand{\eps}{\varepsilon}
\newcommand{\Mac}[1]{\langle #1 \rangle}

\newcommand{\mathd}{\mathrm{d}}
\newcommand{\mbf}[1]{{\mathbf{#1}}}
\newcommand{\jump}[1]{[\![#1]\!]}
\newcommand{\mbfs}{\boldsymbol}
\newcommand{\Tr}[1]{\text{Tr} \left( #1\right) }

\DeclareMathOperator*{\argmin}{argmin}

\usepackage{hyperref}
\hypersetup{
	colorlinks=true,
	linkcolor=blue,
	filecolor=magenta,      
	urlcolor=cyan,
}

\begin{document}
\begin{frontmatter}
\title{On poroelastic strain energy degradation in the variational phase--field models for hydraulic fracture}
\author[Hohai,UFZ,Hohai2]{Tao You}
\author[MUL]{Keita Yoshioka\corref{cor1}}
\address[Hohai]{Key Laboratory of Ministry of Education for Geomechanics and Embankment Engineering, Hohai University, Nanjing 210098, China}
\address[UFZ]{Department of Environmental Informatics, Helmholtz Centre for Environmental Research -- UFZ, Leipzig 04318, Germany}
\address[Hohai2]{College of Civil and Transportation Engineering, Hohai University, Nanjing, China}
\address[MUL]{Department of Petroleum Engineering, Montanuniversit\"at Leoben, Leoben 8700, Austria}
\cortext[cor1]{Corresponding author. }

\begin{abstract}
Though a number of formulations have been proposed for phase--field models for hydraulic fracture, the definition of the degraded poroelastic strain energy varies from one model to another. 
This study explores previously proposed forms of the poroelastic strain energy with diffused fracture and assesses their ability to recover the explicit fracture opening aperture. 
We then propose a new form of degraded poroelastic strain energy derived from micromechanical analyses. 
Unlike the previously proposed models, our poroelastic strain energy degradation depends not only on the phase--field variable (damage) but also on the type of strain energy decomposition. 
Comparisons against closed form solutions suggest that our proposed model can recover crack opening displacement more accurately irrespective of Biot's coefficient or the pore--pressure distribution.
We then verify our model against the plane strain hydraulic fracture propagation, known as the KGD fracture, in the toughness dominated regime.
Finally, we demonstrate the model’s ability to handle complex hydraulic fracture interactions with a pre--existing natural fracture.

\end{abstract}

\end{frontmatter}



\section{Introduction}

The phase--field model for fracture has become a standard method to simulate fracturing.
The model was originally proposed by Bourdin et al.~(2000)~\cite{Bourdin2000} as a numerical regularization of the variational approach to fracture~\cite{francfort1998}, and can be regarded as a special case for the gradient damage model~\cite{Pham2011,lorentz2011gradient,de2016gradient}.
Its applications range from ductile~\cite{Ambati2015, Kuhn2016, Alessi2017,Yin2020}, dynamic~\cite{Bourdin2011, Borden2012,Li2016,Nguyen2018}, fatigue~\cite{Seiler2018,Carrara2020,khalil2022generalised}, dessication~\cite{Maurini2013,Cajuhi2017,Heider2020,luo2023phase}, and to environment assisted fracturing~\cite{Martinez2018phase,Schuler2020chemo,Cui2021phase,wu2021phase,wu2023crack,quinteros2023electromechanical}.
For hydraulic fracture, which is the main focus of this paper, Bourdin et al.~(2012)~\cite{Bourdin2012} adapted the phase--field model by including the effects of fluid pressure within fracture and verified the fracture propagation behaviors in an impermeable elastic medium with injection of inviscid fluid. 
Wheeler et al.~(2014)~\cite{Wheeler2014} extended the phase--field approach to porous media but with an invariant pore--pressure field, and it was further extended for variant pore--pressure by Mikeli\'c et al.(2015a, 2015b)~\cite{Mikelic2015_Multi, Mikelic2015_NonLin}.
Since these pioneering works, many phase--field models for hydraulic fracturing have been proposed.
Amongst these models, the differences amount mainly to how to treat these two points: 1) distinct hydraulic responses in reservoir (porous media) and fracture (open channel), and 2) degradation of the poroelastic strain energy.
We briefly review these two points in the following.

For two distinct hydraulic responses, various ways to delineate the fracture from the reservoir domain have been proposed. 
Mikeli\'c et al. (2015a, 2015b)~\cite{Mikelic2015_Multi, Mikelic2015_NonLin} applied a Darcy flow in reservoir and a Poiseuille flow in fracture and relied on the phase--field value to transition between them.
Another way to delineate the fracture was proposed by Wilson and Landis (2016)~\cite{Wilson2016} who scaled the fluid viscosity to recover Poiseuille flow in fracture.
Santill\`an et al. (2017, 2018)~\cite{santillan2017phase,santillan2018phase} also delineated the fracture domain from the reservoir domain based on the phase--field value, but in their approach, the fracture domain is extracted and built separately with one--dimensional line elements to form a hybrid system with two--dimensional reservoir elements.
This one--dimensional fracture domain then has to be updated dynamically as the phase--field profile evolves.
Another hybrid method was proposed by~\cite{costa2022multi}. 
They applied the phase--field model for fracture propagation and once fractures fully develop, then the fractured domain is treated with the discrete contact model.
Instead of delineating the fracture domain, Miehe et al. (2015)~\cite{miehe2015minimization} enhanced the permeability as a function of the deformation. 
Similarly, Yoshioka and Bourdin (2016)~\cite{Yoshioka2016} enhanced the permeability with a simple threshold based on the phase--field variable.
On the other hand, Lee et al. (2017)~\cite{Lee2017_ls} set a level--set function to average the two flow systems, which has been followed by many studies~\cite{zhou2018phase, li2019numerical,zhuang2020hydraulic,wheeler2020ipacs, liu2020investigation, li2021phase, xu2022phase,zhuang2022three}.
Alternatively, Chukwudozie et al. (2019)~\cite{Chukwudozie2019} averaged the flow in fracture domain using phase--field calculus.
In~\cite{Heider2017, Ehlers2017}, the theory of porous media was used instead to transition between the reservoir flow (Darcy) and the fracture flow (Navier--Stokes).
Despite all the different ways to treat the two distinct flows (porous medium and fracture), the models appear to generate reasonable fluid pressure responses at least in the toughness dominated hydraulic fracture where the pressure drop within the fracture is negligible~\cite{detournay2003near}.

For degradation of the \emph{poroelastic} strain energy, we have yet seen a convergence while the degradation of the \emph{elastic} strain energy is more or less standardized.
For an elastic medium, the strain energy density without fracture or damage is: 
\begin{equation*}
\psi(\mbfs{u}) =  \frac12 \mathbb{C}_\mrm{m}: \mbfs{\eps} (\mbfs{u}) : \mbfs{\eps} (\mbfs{u})
,
\end{equation*}
where $\mathbb{C}_\mrm{m}$ denotes the elastic tensor, and $\mbfs{\eps}$ is the linearlized strain defined as $\mbfs{\eps}=1/2\left(\nabla \mbfs{u} + \nabla^\mrm{T} \mbfs{u}\right)$ with $\nabla (\cdot )$ being the gradient of $(\cdot )$ and $(\cdot )^\mrm{T}$ being the transpose of $(\cdot )$.  
With the phase--field variable $v \in [0, 1]$, which varies from intact $v= 1$ to fully broken $v = 0$, the general degraded strain energy takes the form:
\begin{equation*}
\psi(\mbfs{u},v) =  \frac12 \mathbb{C}_{\mrm{eff}} (v): \mbfs{\eps} (\mbfs{u}) : \mbfs{\eps} (\mbfs{u})
,
\end{equation*}
where $\mathbb{C}_{\mrm{eff}} (v) = v^2 \mathbb{C}_{\mrm{m}}$ in the original model~\cite{Bourdin2000}.
However, other forms of $\mathbb{C}_\mrm{eff} (v)$ have been proposed depending on the degradation function~\cite{Wu2017, wu2018length, Sargado2018}, the type of strain energy decomposition~\cite{Amor2009, miehe2015minimization, Freddi2010, deLorenzis2021nucleation}, and the material constitutive relationship~\cite{Nguyen2020,DIJK2020,ziaei2023orthogonal}.
For poroelastic media, on the other hand, the so called effective stress, which is the part of the stress borne by the solid, is given as $\mbfs{\sigma}_\mrm{eff} = \mathbb{C}_\mrm{m} :\mbfs{\eps}$ and the total stress, which is balanced with the body force and external loadings, is by $\mbfs{\sigma} = \mbfs{\sigma}_\mrm{eff} - {\alpha_\mrm{m}} p \mbfs{\delta}$ where $\alpha_\mrm{m}$ is Biot's coefficient and $p$ is the fluid pressure.
While the effective stress (solid part) is degraded with $\mathbb{C}_{\mrm{eff}} (v)$, degradation for the fluid pressure has not arrived at a consensus within the community.
Some models do not degrade the fluid pressure at all~\cite{heider2020phase,zhou2019phase} and some do but in various ways~\cite{Mikelic2015_CompGeo, Yoshioka2016, santillan2017phase, Lee2017_ls}\footnote{For desiccation crack in partially saturated porous media, we refer to comprehensive stability analyses by~\cite{luo2023phase}.}.


Our objective in this study is to derive a poroelastic strain energy that is consistent with the explicit fracture representation where the fluid pressure is localized within the fracture.
In other words, the poroelastic deformation from the phase--field model should converge to the one from the explicit fracture as the regularization length scale parameter approaches 0. 
We first revisit the previously proposed forms of poroelastic strain energy and point out the differences in terms of degradation of the fluid pressure related terms.
We then derive a degraded poroelastic strain energy from micromechanical analyses.
As it turns out, Biot's parameters depend not only on the phase--field variable (damage) but also on the type of strain energy decomposition.
Our proposed poroelastic strain energy is verified to reproduce the theoretical fracture aperture openings.
On the other hand, the previously proposed models require the fluid pressure in fracture to be diffused beyond fully fractured domain or certain value of Biot's coefficient.

This paper is structured as follows. The next section recalls three existing formulations of poroelastic strain energy in the phase--field model for hydraulic fracture. We then propose a model with phase--field dependent Biot's parameters and porosity via micromechanical analysis. Section \ref{sec:num_impl} deals with the numerical implementation of the new proposed model in which the fixed stress split strategy is derived for solving the hydro--mechanical phase--field coupling  problem. In section \ref{sec:verification}, the proposed model is verified against the analytical solutions. The accuracy of the four involved formulations in estimating the fracture aperture is highlighted. Section \ref{sec:hydraulic_interface} studies interactions between hydraulic fracture and natural fracture. The paper ends with a conclusion and possible future works.

Throughout the paper, the following notations are used.
Tensorial product of any second--order tensors $\boldsymbol{A}$ and $\boldsymbol{B}$ and fourth--order tensor $\mathbb{C}$ are $\boldsymbol{A}:\boldsymbol{B}=A_{ij}B_{ij}$ and $\mathbb{C}:\boldsymbol{B}=\mathbb{C}_{ijkl}\boldsymbol{B}_{kl}$. 
The symbol $\left\Vert \boldsymbol{A}\right\Vert=\sqrt{\boldsymbol{A}:\boldsymbol{A}}$ is used to calculate the norm of any second--order tensor $\boldsymbol{A}$. 
The identity tensors for the second and fourth order are denoted by $\mbfs{\delta}$ and $\mathbb{I}$, respectively.
The fourth order projector $\mathbb{J}$ and $\mathbb{K}$ are expressed in the component form as $J_{ijkl}=\delta _{ij}\delta _{kl}/3$ and $K_{ijkl}=(\delta _{ik}\delta_{jl}+\delta _{jl}\delta _{jk})/2-\delta _{ij}\delta _{kl}/3$, respectively. 
The trace operator $\Tr{\cdot}$ is defined as $\Tr{\mbfs{A}} = \mbfs{\delta}: \mbfs{A}$.
$\nabla (\cdot )$ is the gradient of $(\cdot )$, $(\cdot )^\mrm{T}$ is the transpose of $(\cdot )$, and $\nabla^\mrm{s} (\cdot ) $ is the symmetric part of the gradient of $(\cdot )$ defined as $\nabla^\mrm{s} (\cdot ) := \frac12 (\nabla (\cdot ) + \nabla^\mrm{T} (\cdot ))$. 
The overhead dot symbol $\dot{A}$ indicates the time differentiation of variable $A$. 

\section{Model formulation} \label{model}

In this section, we start with the now established phase--field model by introducing a total energy of an elastic medium.
The total energy is composed of the strain energy and the surface energy, and the strain energy needs to be adjusted when applied to poroelastic medium.
To explore the required adjustment, we present a three--scale (macro--, meso--, and micro--) poroelastic medium. 
To clarify, we use ``meso--" to indicate an object or quantity at the mesoscale, and ``micro--" at the microscale.
Two--level homogenization will be presented first to derive the poroelastic behavior, considering meso--cracks and micro--pores. 
Following that, the phase--field dependent Biot's parameters and porosity will be introduced to account for macroscale fracture. 
In this study, we deal with the hydraulic fracture problem in fluid--saturated porous media. 

\subsection{Phase--field model for fracture} \label{microPF}
\begin{figure}[htp!]
	\centering
	\includegraphics[scale=0.5]{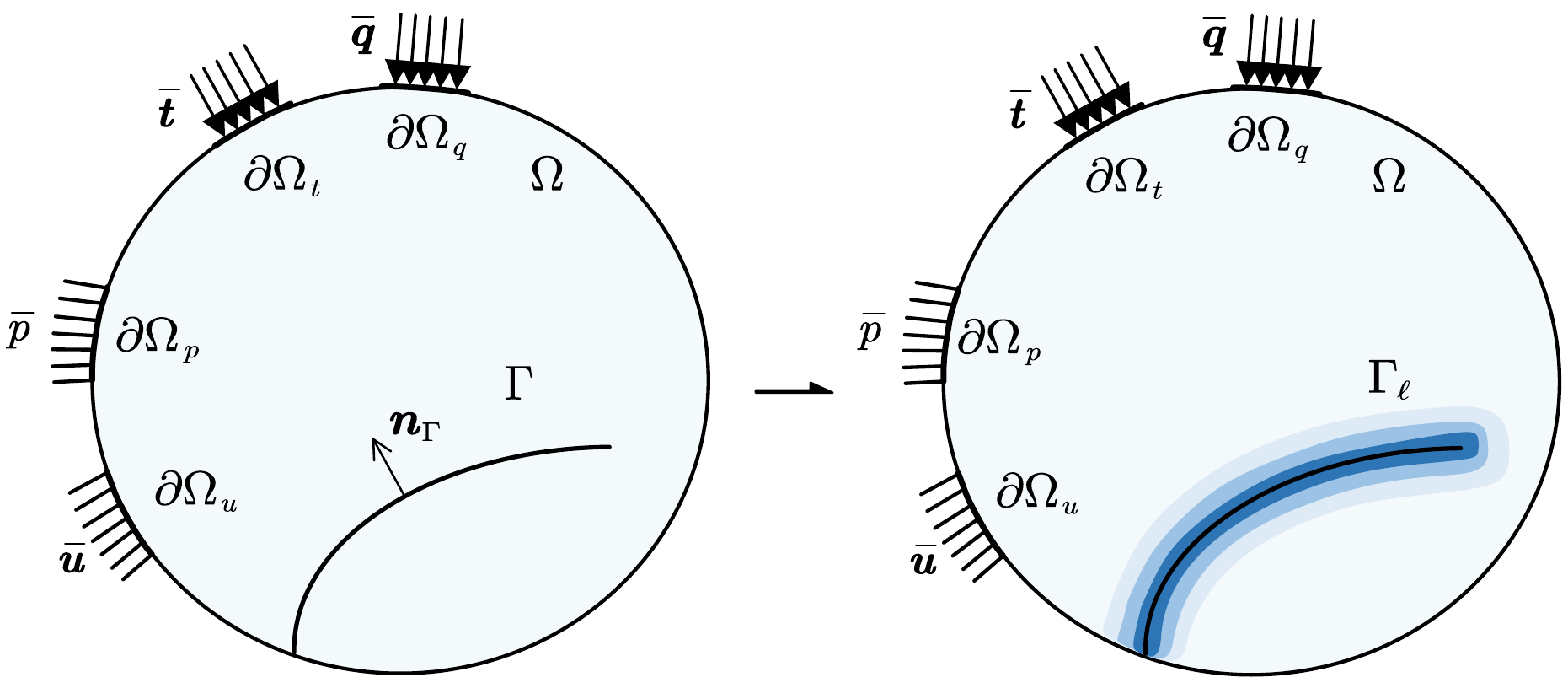}
	\caption{Schematic illustration of the cracked porous media where the macroscopic fracture is regularized by the phase--field method.}
	\label{fig:object}
	
\end{figure}

Consider an elastic body $\Omega \subset \mathbb{R}^{n_\mrm{dim}}$ $\left(n_{\mrm{dim}} = 2,3\right)$ embedded with a lower--dimensional macroscopic fracture $\Gamma \subset \mathbb{R}^{n_{\mrm{dim}-1}}$ which is regularized by the well--known phase--field method into $\Gamma_{\ell} \subset \mathbb{R}^{n_{\mrm{dim}}}$ (Fig.~\ref{fig:object} a).  
The body is subjected to a surface force $\mbfs{\bar{t}}$ on boundary $\partial \Omega_t$ and a flux $\mbfs{\bar{q}}$ on $\partial \Omega_q$.  
In addition, $\bar{p}$ and $\bar{\mbfs{u}}$ are the prescribed pressure and displacement on boundary $\partial \Omega_p$ and $\partial \Omega_u$, respectively. 

The phase--field model is based on the variational approach to fracture proposed by Francfort and Marigo~\citep{francfort1998} and they defined the total energy of the system (Fig.~\ref{fig:object}) as a sum of the strain energy and the surface energy:
\begin{equation}	
	\label{eq:total_energy}
	\mathcal{E} (\mbfs{u},\Gamma):= \int_{\Omega\setminus\Gamma} \psi(\mbfs{u})  \, \mathd V  
	+ \int_\Gamma \Gc \, \mathd S
 ,
\end{equation}
where $\psi(\mbfs{u}) $ is the strain energy density and $\Gc$ is the fracture surface energy density.
Bourdin et al.~\cite{Bourdin2000} proposed to regularize this energy functional following~\cite{mumford1989optimal} by introducing a phase--field variable $v \in [0, 1]$\footnote{In this study, we follow the original definition of the phase--field variable where $v=0$ represents fracture and $v=1$ intact material~\cite{Bourdin2000}. Thus, $(1-v)$ corresponds to the degree of damage in the context of damage mechanics.} and a regularization length parameter $\ell$ as
\begin{equation}	
	\label{eq:total_energy_reg}
	\mathcal{E}_\ell(\mbfs{u},v) := \int_{\Omega}  \psi(\mbfs{u},v)  \, \mathd V  
	+  \int_\Omega \frac{\Gc}{4c_n} \left[\frac{(1-v)^n}{\ell}+ \ell \nabla v\cdot \nabla v \right]\,\mathd V
	,
\end{equation}
where $\psi(\mbfs{u}, v) $ is given as $ \psi(\mbfs{u},v)= v^2\psi(\mbfs{u})$ in~\cite{Bourdin2000} and $c_n$ is the normalizing parameter given by $c_n := \int_0^1(1-\omega)^{n/2}\mathd \omega$~\cite{mesgarnejad2015validation,Tanne2018}.
The regularized model is referred as \ATone{} for $n=1$ and as \ATtwo{} for $n=2$~\cite{Pham2011,Bourdin2014}.

\subsection{Extension to poroelastic media} \label{sec:pf models}
When the phase--field model is applied for a poroelastic medium, we need to account for poroelastic effects in the total energy functional $\mathcal{E}_\ell$. 
We consider a poroelastic medium ($\Omega$) that consists of a porous matrix ($\Omega\setminus\Gamma$) and a crack set ($\Gamma$) (Fig. \ref{fig:object}).
Then, we wish to find the regularized effective total energy functional for the poroelastic medium, which can be written as 
\begin{equation}
\label{eq:F_l}
\mathcal{F}_\ell (\mbfs{u},v,p)	:= W(\mbfs{u},v,p) +  \int_\Omega \frac{\Gc}{4c_n} \left[\frac{(1-v)^n}{\ell}+ \ell \nabla v\cdot \nabla v \right]\,\mathd V
,
\end{equation}
where $W(\mbfs{u},v,p)$ is the poroelastic strain energy.
Several variations of $W(\mbfs{u},v,p)$ have been proposed previously and we explore them in the following.

We start with the linear momentum balance and mechanical boundaries of the poroelastic medium in a discrete setting depicted in Fig.~\ref{fig:object} (left):
\begin{align}
	\label{eq:momentum}
	\begin{cases}
		\nabla \cdot \boldsymbol{\sigma} + \mbfs{b} = \mathbf{0} & \,\mrm{in}\,\Omega\setminus\Gamma 
		\\
		\boldsymbol{\sigma} \cdot \mbfs{n}_\Gamma  = -p \mbfs{n}_\Gamma & \,\mrm{on}\, \Gamma
		\\
		\boldsymbol{\sigma} \cdot \mbfs{n}  = \mbfs{\bar{t}} & \,\mrm{on}\, \partial{\Omega_t}
		\\
		\mbfs{u}  = \mbfs{\bar{u}} & \,\mrm{on}\, \partial{\Omega_u}
	\end{cases}
\end{align}
where $\mbfs{\sigma}$ denotes the Cauchy (total) stress, $\mbfs{b}$ the body force and $\mbfs{n}_\Gamma$ the normal vector to the fracture.
Let the elasticity stiffness tensor for the poroelastic medium be $\mathbb{C}_\mrm{m} = 3k_\mrm{m} \mathbb{J}+2\mu_\mrm{m} \mathbb{K}$ where $k_\mrm{m}$ and $\mu_\mrm{m}$ are the bulk modulus and the shear modulus for the porous matrix respectively. 
Following Biot's poroelasticity~\cite{Biot1941}, the total stress $\mbfs{\sigma}$ in the porous matrix ($\Omega\setminus\Gamma$) is given as
\begin{equation}
	\label{eq:macroscopic stress}
	\mbfs{\sigma} = \mathbb{C}_\mrm{m} :\mbfs{\eps} - {\alpha_\mrm{m}} p \mbfs{\delta}
 .
\end{equation}

Let us neglect the body force $\mbfs{b}$ and the traction $\mbfs{\bar{t}}$ in the comparisons that follow in this section for the sake of simplicity.
Multiplying \eqref{eq:momentum}-1 by a test function $\mbfs{w}_u \in H^1(\Omega\setminus\Gamma)$ and using Green's formula, \eqref{eq:momentum}--2 and \eqref{eq:momentum}--3, we have
\begin{align}
	\label{eq:weak_momentum}
	&\int _{\Omega\setminus\Gamma} \left[ \mathbb{C}_\mrm{m} : \mbfs{\varepsilon}(\mbfs{u}) -  {\alpha_\mrm{m}} p \mbfs{\delta} \right] \cdot \mbfs{\eps}(\mbfs{w}_u)  \mathd V  
	-  \int _{\Gamma} p \jump{  \mbfs{w}_u \cdot \mbfs{n}_\Gamma } \, \mathd S= \mathbf{0} 
	,
\end{align}
where $\jump{\cdot}$ denotes a jump quantity over $\Gamma$.
In~\cite{Yoshioka2016, Chukwudozie2019}, they considered the poroelastic strain energy of:
\begin{align}
	\label{eq: W1_discrete}
	W_1(\mbfs{u};p) :=  
	\int _{\Omega\setminus\Gamma}  \frac{1}{2}\mathbb{C}_\mrm{m} : \left( \mbfs{\varepsilon}(\mbfs{u}) -  \frac{\alpha_\mrm{m} p}{ 3 k_\mrm{m}}  \mbfs{\delta} \right) : \left( \mbfs{\varepsilon}(\mbfs{u}) -  \frac{\alpha_\mrm{m} p}{ 3 k_\mrm{m}} \mbfs{\delta} \right)  \mathd V 
	-  \int _{\Gamma} p \jump{  \mbfs{u} \cdot \mbfs{n}_\Gamma } \, \mathd S
	,
\end{align}
so that \eqref{eq:weak_momentum} is the first variation of \eqref{eq: W1_discrete} with respect to $\mbfs{u}$ assuming that $W_1$ is a function of $\mbfs{u}$.
Using a $\Gamma$-convergence property and the co-area formula~\cite{Bourdin2012, Chukwudozie2019}, the last term of \eqref{eq: W1_discrete} can be approximated as
\begin{equation*}
	\int _{\Gamma} p \jump{  \mbfs{u} \cdot \mbfs{n}_\Gamma } \mathd S
	\approx
	\int_{\Omega} p \mbfs{u} \cdot \nabla v \, \mathd V
	.
\end{equation*}
Then the regularized poroelastic strain energy can be
\begin{align}
	\label{eq: W1_reg1}
	W_1(\mbfs{u}, v;p) \approx  
	\int _{\Omega}  \frac{1}{2} \mathbb{C}_\mrm{m} : \left(v \mbfs{\varepsilon}(\mbfs{u}) - \frac{\alpha_\mrm{m} p}{ 3 k_\mrm{m}} \mbfs{\delta} \right) : \left( v\mbfs{\varepsilon}(\mbfs{u}) - \frac{\alpha_\mrm{m} p}{ 3 k_\mrm{m}} \mbfs{\delta} \right)  \mathd V 
	-\int_{\Omega} p \mbfs{u} \cdot \nabla v \, \mathd V
\end{align}
Because $p$ was considered constant with respect to $\mbfs{u}$ in~\cite{Yoshioka2016, Chukwudozie2019} (i.e., ${\alpha_\mrm{m} p}/{ 3 k_\mrm{m}} = \text{const.}$), the poroelastic strain energy can be equivalently written:
\begin{align*}
	\label{eq: W1_reg2}
	W_1(\mbfs{u}, v;p) =  
	\int _{\Omega}  \frac{v^2}{2} \mathbb{C}_\mrm{m} :  \mbfs{\varepsilon}(\mbfs{u}) :  \mbfs{\varepsilon}(\mbfs{u})\, \mathd V  
	- \int _{\Omega}  v \alpha_\mrm{m} p  \mbfs{\varepsilon}(\mbfs{u}) : \mbfs{\delta}\, \mathd V  
	-\int_{\Omega} p \mbfs{u} \cdot \nabla v \, \mathd V
	.
\end{align*}
Integrating the last term by parts and denoting $\psi(\mbfs{u},v) =  \frac{v^2}{2} \mathbb{C}_\mrm{m}: \mbfs{\eps} (\mbfs{u}) : \mbfs{\eps} (\mbfs{u})$, we have
\begin{equation}
	\label{eq:W1_reg_final}
   W_1(\mbfs{u}, v;p) = 
	 \int_{\Omega } \psi (\mbfs{u},v) \, \mathd V
	+ \int_{\Omega } v \left( 1- \alpha_\mrm{m}  \right) p \nabla \cdot \mbfs{u} \, \mathd V 
	+ \int_{\Omega} v \nabla p \cdot \mbfs{u} \, \mathd V - \int_{\partial \Omega}v  p \mbfs{u}\cdot \mbfs{n} \, \mathd S
 - \cancel{\int_{\Gamma} v p \mbfs{u} \cdot \mbfs{n}_\Gamma \, \mathd S}
	.
\end{equation}	
where the last term is canceled because $v=0$ in the fracture domain.

On the other hand in~\cite{Mikelic2015_CompGeo}, though derived differently, the poroelastic strain energy is defined as
\begin{equation}	
	\label{eq:W2_reg_final}
	W_2(\mbf{u}, v;p) := 
	 \int_{\Omega } \psi (\mbfs{u},v) \, \mathd V
	+ \int_{\Omega } v^2 \left( 1- \alpha_\mrm{m}  \right) p \nabla \cdot \mbfs{u} \, \mathd V 
	+ \int_{\Omega} v^2 \nabla p \cdot \mbfs{u} \, \mathd V -  \int_{\partial \Omega}  p \mbfs{u}\cdot \mbfs{n} \, \mathd S
	.
\end{equation}
 Note that, in \cite{Mikelic2015_CompGeo}, the authors assumed fractures do not reach the domain boundary $\partial \Omega$ so that $v = 1$ on the boundary $\partial \Omega$, which means that the last term of $W_2$ corresponds to that of $W_1$ ($\int_{\partial \Omega}  p \mbfs{u}\cdot \mbfs{n} \, \mathd S = \int_{\partial \Omega} v p \mbfs{u}\cdot \mbfs{n} \, \mathd S$).
Thus, one can observe that the difference between $W_1$ and $W_2$ amounts to the degradation function for the pressure related terms ($v$ or $v^2$)\footnote{In~\cite{Mikelic2015_Multi,  Mikelic2015_NonLin}, $v^2$ was chosen over $v$ because $v$ may go below 0 during the minimization, which causes numerical instabilities. In~\cite{Bourdin2012, Yoshioka2016, Chukwudozie2019}, the solution space for $v$ is bounded in $[0,1]$ by applying a variational inequality solver for $v$, which permits $v$ to be a degradation function for the fluid pressure.}. 

Furthermore, we remark another popular form for the poroelastic strain energy, which does not degrade the fluid pressure related terms. 
Assuming the continuity of $\mbfs{\sigma}$ over $\Omega$ (i.e. without imposing~\eqref{eq:momentum}--2), one obtains~\cite{zhou2019phase, Heider2020}:
\begin{equation}	
	\label{eq:W0_reg_final}
	W_0(\mbf{u}, v;p) := 
	 \int_{\Omega } \psi (\mbfs{u},v) \, \mathd V
	- \int_{\Omega }  \alpha_\mrm{m}   p \nabla \cdot \mbfs{u} \, \mathd V 
	.
\end{equation}

Now, we have introduced three common forms of the poroelastic strain energy ($W_0$, $W_1$, and $W_2$).
In these notations of the poroelastic strain energies, we use different subscripts (0, 1, and 2) to indicate their exponent of $v$ that multiples the pressure terms ($v^0=1$, $v^1=v$, and $v^2$).

We propose yet another form here.
Following~\citep{biot1962mechanics, Mikelic2015_CompGeo}, we define $W$ as
\begin{align}
	\label{eq:poroelasticityVariational}
	W \left(\mbfs{u}, v,p\right)
	:= 	 & \int_{\Omega } \psi (\mbfs{u},v) \, \mathd V  +\int_\Omega \frac{M}{2} \left[ \alpha \nabla \cdot \mbfs{u} - \zeta\right]^2 \, \mathd V \\
	=& 	 \int_{\Omega } \psi (\mbfs{u},v) \, \mathd V + \int_\Omega \frac{p^2}{2M} \, \mathd V
	,
\end{align}
where $\alpha$ and $M$ is Biot's coefficient and modulus for poroelastic medium ($\Omega$) and $\zeta$ is the variation of fluid content given as
\begin{equation}
    \label{eq:zeta}
    \zeta = \alpha \nabla \cdot \mbfs{u} + \frac{p}{M}
    .
\end{equation}
As $p$ depends on $\mbfs{\eps}$ (i.e. ${\partial p}/{\partial \mbfs{\eps}} = {\alpha}/{M}\mbfs{\delta}$)\footnote{In \cite{Mikelic2015_CompGeo}, $W$ is defined similarly to \eqref{eq:poroelasticityVariational}, but $p$ is regarded constant with respect to $\mbfs{\eps}$.}, the first variation of $W$ with respect to $\mbfs{\eps}$ gives
\begin{equation}
\frac{\partial W}{\partial \mbfs{\eps}} = 
    \int _{\Omega} \left[ \mathbb{C}_\mrm{eff} (v) : \mbfs{\eps}(\mbfs{u}) -  {\alpha(v)} p \mbfs{\delta} \right] \mathd V 
    ,
\end{equation}
where $\mathbb{C}_\mrm{eff}$ is the effective stiffness tensor for the poroelastic medium ($\Omega$) including the porous matrix ($\Omega\setminus\Gamma$) and the crack set ($\Gamma$).
Note that $\mathbb{C}_\mrm{eff}$ is not the same as the effective stiffness tensor for the porous matrix ($\Omega\setminus\Gamma$) i.e., $\mathbb{C}_\mrm{eff} \neq \mathbb{C}_\mrm{m}$.
Here, we indicate that both $\mathbb{C}_\mrm{eff}$ and $\alpha$ depend on the damage ($v$) so that the total stress is continuous over $\Omega$, and the condition $\mbfs{\sigma} \cdot \mbfs{n}_{\Gamma} = -p \mbfs{n}_{\Gamma}$ in Eq. \eqref{eq:momentum}--2 is ensured.
Biot's parameters proposed in \cite{yi2020consistent,ulloa2022variational,zhang2023phase} also depend on the damage ($v$), but we derive their dependency on the damage ($v$) from micromechanics analyses in the following section.




\subsection{Micromechanical analysis}

\begin{figure}[htp!]
	\centering
	\includegraphics[scale=0.5]{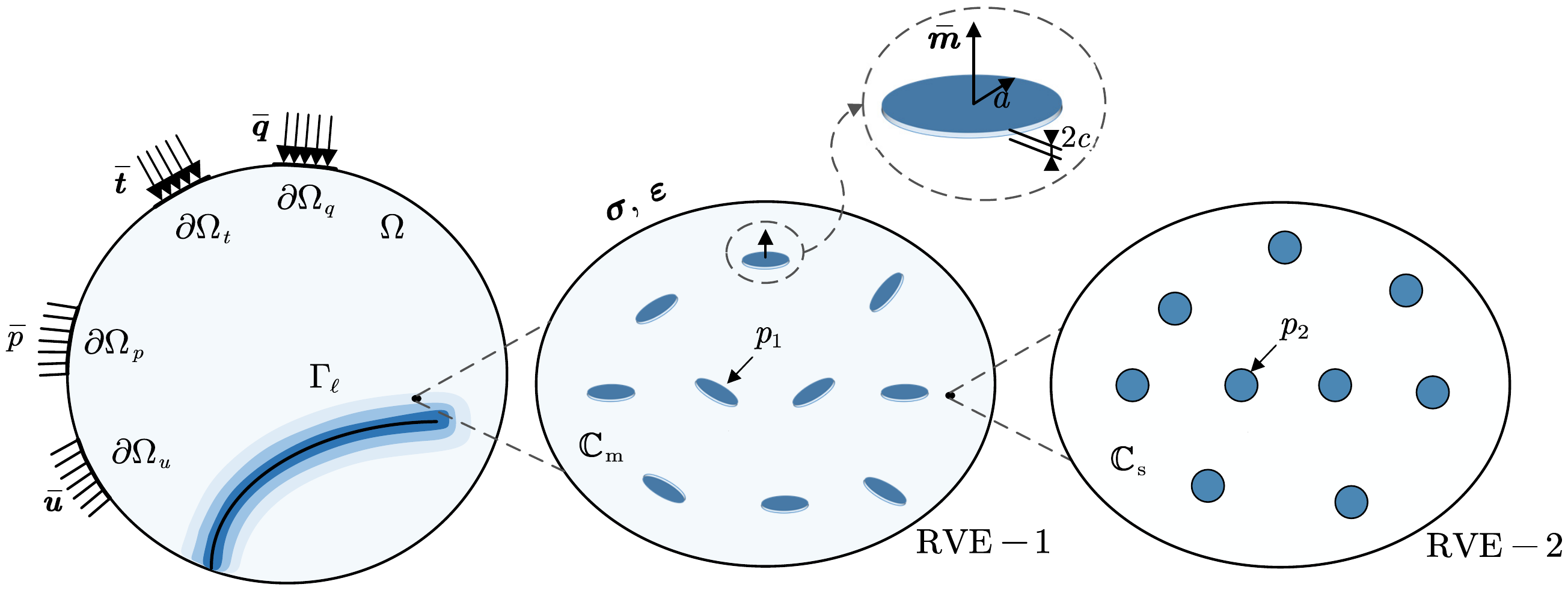}
	\caption{Schematic illustration of the scale separation of the cracked porous media where the microcrack and micropores are filled with fluid.}
	\label{fig:scale-seperation}
\end{figure}

We consider the material at the mesoscale, i.e., RVE--1, which is composed of porous matrix and distributed meso--cracks (Fig.~\ref{fig:scale-seperation} (middle)). 
At the lower microscale, i.e., RVE--2, the porous material is made up of solid matrix and micro--pore (Fig.~\ref{fig:scale-seperation} (right)). 
The porous media is fully fluid--saturated and drained, which means that the fluids in meso--cracks and micro--pores are interconnected, and the fluid pressures in meso--cracks (denoted as $p_1$) and micro--pores (denoted as $p_2$) are hydraulically equilibrated, i.e., $p_1 = p_2 = p$. 
These assumptions of interconnectivity are important to simplify the following derivation of homogenized poroelastic behaviors. 

\subsubsection{Homogenized poroelastic behavior}\label{homogenization}

Given the two--scale RVEs composed of randomly distributed penny--shaped meso--cracks and micro--pores respectively  (Fig.~\ref{fig:scale-seperation}), the porous matrix is assumed to be linearly elastic with the elasticity tensor $\mathbb{C}_\mrm{m} = 3k_\mrm{m} \mathbb{J}+2\mu_\mrm{m} \mathbb{K}$. 
The stiffness tensor for the solid matrix in RVE--2 is $\mathbb{C}_\mrm{s} = 3k_\mrm{s} \mathbb{J}+2\mu_\mrm{s} \mathbb{K}$ where $k_\mrm{s}$ and $\mu_\mrm{s}$ denote the bulk modulus and shear modulus of mineral grains. 
Under the isotropic assumption, the volume fraction $\varphi$ of penny-shaped meso-cracks in RVE--1 can be expressed as \citep{budiansky1976elastic}
\begin{equation}
	\varphi =\frac{4}{3}\pi\mathcal{N}a^{2} X =\frac{4}{3} \pi \omega X
	\label{eq: volume fraction}
	,
\end{equation}%
where $X=c/a$ is the aspect ratio of meso-crack, $\mathcal{N}$ is the density of meso--cracks (i.e., the number of meso-cracks per unit
volume), and $\omega =\mathcal{N}a^{3}$ is the microcrack-density parameter. 
Additionally, the pore volume in RVE--2 is denoted as $\phi_\mrm{m}$.

The classical Biot theory~\cite{Biot1941} can be derived from the response of a representative element volume (i.e., RVE) subjected to a mechanical loading defined by uniform strain boundary conditions ($\nabla ^\mrm{s}\boldsymbol{u}$) and the fluid pressure ($p$) as has been investigated by Dormieux et al. \citep{dormieux2005poroelasticity,dormieux2006microporomechanics}. 
For this two-scale double--porosity problem, Biot's coefficient tensor for RVE--1 can be written as
\begin{equation}
	\mbfs{\alpha} = \mbfs{\alpha}_1 + \mbfs{\alpha}_2
	\label{eq:alpha}
	,
\end{equation}
where $\mbfs{\alpha}_1$ and $\mbfs{\alpha}_2$ are Biot's coefficient tensors associated with the pressure in the meso-cracks and micro-pores, respectively.
Assuming $p_1=p_2=p$, we have
\begin{align}
	\label{eq:alpha1}
	\mbfs{\alpha}_1 & = \mbfs{\delta}:\left( \mathbb{I} - (1-\varphi)\mathbb{\bar{A}}_\mrm{m}\right)
	,
	\\
	\label{eq:alpha2}
	\mbfs{\alpha}_2 & = (1-\varphi) \mbfs{\alpha}_\mrm{m}:\mathbb{\bar{A}}_\mrm{m}
	,
\end{align}
where $\mathbb{\bar{A}}_\mrm{m}$ is the average strain concentration tensor of the porous matrix and $\mbfs{\alpha}_\mrm{m}$ is isotropic Biot's coefficient tensor for the porous matrix given by
\begin{equation}
	\mbfs{\alpha}_\mrm{m} = \alpha_\mrm{m} \mbfs{\delta}
	\label{eq:alpha_m}
	.
\end{equation}
Following~\citep{coussy2004poromechanics}, we have
\begin{equation}
	\label{alpha_m}
	\alpha_\mrm{m} = 1 - \frac{k_\mrm{m}}{k_\mrm{s}}
	.
\end{equation}

Assuming $p_1=p_2=p$ again and  incompressible fluid, Biot's modulus of RVE--1 is given as~\citep{dormieux2006microporomechanics}
\begin{equation}
	\frac{1}{N} = \frac{1}{N_{11}} + \frac{1}{N_{12}} + \frac{1}{N_{21}} + \frac{1}{N_{22}}
	\label{eq:M}
	,
\end{equation}
with
\begin{align}
	\begin{cases}
		\frac{1}{N_{11}} =  \mbfs{\alpha}_\mrm{m}:\mathbb{S}_\mrm{m}:\left((1-\varphi)\mbfs{\alpha}_\mrm{m} - \mbfs{\alpha}_2\right) + \frac{1 - \varphi}{N_\mrm{m}} \\
		\frac{1}{N_{12}} = \mbfs{\alpha}_\mrm{m}:\mathbb{S}_\mrm{m}:\left(\varphi\mbfs{\delta} - \mbfs{\alpha}_1\right) \\
		\frac{1}{N_{21}} = \mbfs{\delta}:\mathbb{S}_\mrm{m}:\left(\mbfs{\alpha}_2 - (1-\varphi) \mbfs{\alpha}_\mrm{m}\right) \\
		\frac{1}{N_{22}} = \mbfs{\delta}:\mathbb{S}_\mrm{m}:\left(\mbfs{\alpha}_1 - \varphi \mbfs{\delta}\right)
	\end{cases}
	\label{eq:M_matrix}
\end{align}
where $N_\mrm{m}$ is Biot's modulus of the porous matrix, and $\mathbb{S}_\mrm{m}$ is the compliance tensor defined as $\mathbb{S}_\mrm{m} = \left(\mathbb{C}_\mrm{m}\right)^{-1}$.
Using Eq. \eqref{eq:alpha_m}, we have
\begin{equation}
	\frac{1}{N_\mrm{m}} = \left( \mbfs{\alpha}_\mrm{m}- \phi_\mrm{m} \mbfs{\delta} \right) :\mathbb{S}_\mrm{s}:\mbfs{\delta} =  \frac{\alpha_\mrm{m} - \phi_\mrm{m}}{k_\mrm{s}}
	\label{eq:Mm}
	.
\end{equation}

We now turn our attention to the macroscopic deformation in RVE--1. 
The macroscopic strain $\mbfs{\varepsilon}$ and $p$ are the boundary conditions of RVE--1 as illustrated in Fig. \ref{fig:scale-seperation}. 
The homogenized or effective stiffness tensor $\mathbb{C}_\mrm{eff}$ of RVE--1 for open meso-cracks takes the form~\citep{zaoui02,zhu2008micromechanical}
\begin{equation}
	\mathbb{C}_\mrm{eff}=(1-\varphi)\mathbb{C}_\mrm{m}:\mathbb{\bar{A}}_\mrm{m}
	\label{eq:chom_general}
	,
\end{equation}
which can also be expressed as $\mathbb{C}_\mrm{eff} = 3k_\mrm{eff} \mathbb{J}+2\mu_\mrm{eff} \mathbb{K}$ where $k_\mrm{eff}$ and $\mu_\mrm{eff}$ denote the effective bulk modulus and shear modulus of RVE--1 respectively.

Finally, combining Eqs. \eqref{eq:alpha}, \eqref{eq:alpha1}, \eqref{eq:alpha2}, \eqref{eq:alpha_m}, and \eqref{eq:chom_general}, we can link Biot's coefficient with the effective stiffness tensor as

\begin{equation}
\begin{aligned}
\mbfs{\alpha} 
& = \mbfs{\delta} - (1 - \alpha_\mrm{m})(1-\varphi) \mathbb{\bar{A}}_\mrm{m} : \mbfs{\delta}
\\
& =\mbfs{\delta} - \frac{k_\mrm{m}}{k_\mrm{s}} \mathbb{C}_\mrm{eff}:\mathbb{S}_\mrm{m} : \mbfs{\delta}
\\
& = {\left(1-\frac{k_\mrm{eff}}{k_\mrm{s}}\right)}\mbfs{\delta}   
\label{eq:alpha_micro}
,
 \end{aligned}
\end{equation}
Therefore, we have the effective Biot's coefficient of RVE--1 as
\begin{equation}
\alpha =   \left(1-\frac{k_\mrm{eff}}{k_\mrm{s}}\right)
.
\end{equation}
Furthermore, combining Eqs.~\eqref{eq:M}, \eqref{eq:M_matrix}, \eqref{eq:Mm}, and \eqref{eq:chom_general}, we obtain
\begin{equation}
\begin{aligned}
\frac{1}{N}	 
&  = (1-\alpha_\mrm{m})\left( \mbfs{\alpha} - \varphi \mbfs{\delta} - (1-\varphi)  \mbfs{\alpha}_\mrm{m}  \right): \mathbb{S}_\mrm{m}:\mbfs{\delta}
    + (1-\varphi) \left( \mbfs{\alpha}_\mrm{m}- \phi_\mrm{m} \mbfs{\delta} \right) :\mathbb{S}_\mrm{s}:\mbfs{\delta}
    \\
& =  (1-\alpha_\mrm{m})\left( \mbfs{\alpha} - \varphi \mbfs{\delta} - (1-\varphi) \phi_\mrm{m} \mbfs{\delta} \right): \mathbb{S}_\mrm{m}:\mbfs{\delta}
\\
& =  \frac{k_\mrm{m}}{k_\mrm{s}} \left( \alpha - \varphi  - (1-\varphi) \phi_\mrm{m} \right)\mbfs{\delta}: \mathbb{S}_\mrm{m}:\mbfs{\delta}
\\
& = \frac{\alpha - \phi}{k_\mrm{s}} 
   \label{eq:M_inv}
   .
\end{aligned}
\end{equation}
where $\phi$ is the effective porosity of the RVE--1 and can be expressed as 
\begin{equation}
	\phi = \varphi + (1-\varphi) \phi_\mrm{m}
	\label{eq:effective_porosity}
	.
\end{equation}

Subsequently, by combining Eqs. \eqref{eq:alpha_m}, \eqref{eq:alpha_micro}, and \eqref{eq:M_inv}, Biot's modulus with considering slight fluid compressibility can be written as
\begin{equation}
	\label{eq:M_micro}
	\frac{1}{M} = \frac{1}{N} + \phi c_\mrm{f}  =  \frac{\alpha - \phi}{k_\mrm{s}} + \phi c_\mrm{f}   
	,
\end{equation}
where $c_\mrm{f}$ is the fluid compressibility.

At this point, we have derived from the two-scale poroelastic behavior that Biot's parameters ($\mbfs{\alpha}$, $M$), porosity ($\phi$) and effective stiffness tensor ($\mathbb{C}_\mrm{eff}$) are a function of mesoscopic damage parameter $\omega$ which physically stands for the density of microcracks per unit volume.



\subsubsection{Phase--field-dependent Biot's parameters and porosity}
Note that we have not specified a particular homogenization scheme, e.g., dilute \cite{eshelby1957determination}, Mori-Tanaka \cite{mori1973average} or Ponte Castaneda and Willis (PCW) \cite{castaneda1995effect} in our derivations, which means the expression of $\mathbb{C}_\mrm{eff}$ in Eq. \eqref{eq:chom_general} is a general form with $\mathbb{A}_\mrm{m}$ unspecified. 
In the present work, $\mathbb{C}_\mrm{eff}$ will be derived directly from the phase--field model.

Recall, in the original phase--field model~\cite{Bourdin2000}, the strain energy is degraded regardless of the deformation direction, i.e., $ \psi(\mbfs{u},v)= v^2\psi(\mbfs{u})$.
However, a more general form reads~\cite{Amor2009,Freddi2010,deLorenzis2021nucleation},
\begin{equation}
	\psi  (\mbfs{u},v)=g\left( v\right) \psi_+(\mbfs{u}) + \psi_-(\mbfs{u})
	\label{eq:free energy pf}
	,
\end{equation}%
where $g(v)$ is the degradation function representing the material stiffness deterioration due to the regularized macroscopic crack.
In this study, we define it as
\begin{equation}
	g(v) = v^2 
	\label{mv}
\end{equation}
In addition, $\psi_+$ and $\psi_-$ in Eq. \eqref{eq:free energy pf} are the positive and negative parts of strain energy and only the positive part contributes to the phase--field evolution. 
For the isotropic strain energy model~\cite{Bourdin2000}, we have $\psi_+(\mbfs{u}) = \psi (\mbfs{u})$ and $\psi_- (\mbfs{u}) = 0$, which reads
\begin{align}
	\psi^\mrm{iso} ((\mbfs{u}),v)&=g\left( v\right) \frac{1}{2}\boldsymbol{\varepsilon }:\mathbb{C}_\mrm{m}:%
	\boldsymbol{\varepsilon } 
\end{align}
For the volumetric-deviatoric (V--D) decomposition model~\cite{Amor2009}, we have
\begin{align}
	\psi^\mrm{vd} ((\mbfs{u}),v)&=g(v) \left[\frac{k_\mrm{m}}{2} \Mac{\Tr{\bm{\eps}}}_+^2 + \mu_\mrm{m} \boldsymbol{\varepsilon }: \mathbb{K}:\boldsymbol{\varepsilon}\right] + \frac{k_\mrm{m}}{2} \Mac{\Tr{\bm{\eps}}}_-^2
	,
\end{align}
where $\Mac{\cdot}$ denotes the Macaulay brackets defined as $\Mac{\cdot}_\pm = (|\cdot|\pm \cdot)/2$. 
The stiffness tensors can then be derived as
\begin{equation}
	\mathbb{C}_\mrm{eff}(v) = 
	\begin{cases}
		3 \, \underset{k_\mrm{eff}(v)}{\underbrace{g(v) k_\mrm{m}}} \, \mathbb{J} + 2 g(v) \mu_\mrm{m} \mathbb{K} & \text{for isotropic model}\\
		3 \, \underset{k_\mrm{eff}(v)}{\underbrace{ \left[ g(v) H (\Tr{\bm{\eps}} ) + H (\Tr{-\bm{\eps}} ) \right] k_\mrm{m}}} \,\mathbb{J} + 2 g(v) \mu_\mrm{m} \mathbb{K} & \text{for V--D decomposition model}
	\end{cases}
	\label{eq:relation of stiffness}
\end{equation}
where the underlined quantities correspond to the effective bulk modulus $k_\mrm{eff}$, and $H(\cdot)$ is the Heaviside step function defined as
\begin{align}
	H(x) := \begin{dcases}
		0, & \, x < 0, \\
		1, & \, x \ge 0.
	\end{dcases}
\end{align}

Note that Eq. \eqref{eq:chom_general} is only valid for open cracks and therefore a tension-compression transition criterion should be considered. Apparently, for the isotropic model in Eq. \eqref{eq:relation of stiffness}--1, $k_\mathrm{eff}(v)=3g(v)k_\mrm{m}$ holds both in tension and compression which fails to capture the unilateral effect of fracture in geomaterials. 
On the contrary, the V--D decomposition model uses the sign of volumetric strain as the tension-compression transition criterion. Without loss of generality, we take the V--D decomposition model in the subsequent derivations.  
As a result, Biot's coefficient tensor \eqref{eq:alpha_micro} and Biot's modulus \eqref{eq:M_micro} can be expressed by using Eq. \eqref{eq:relation of stiffness} and Eq. \eqref{alpha_m} as
\begin{align}
	\mbfs{\alpha}(v) &= \left[ 1- \left[g(v) H (\Tr{\bm{\eps}} ) + H (\Tr{-\bm{\eps}} )\right]\frac{k_\mrm{m}}{k_\mrm{s}}\right] \mbfs{\delta} = \alpha(v) \mbfs{\delta} 
	\label{eq:alpha_macro}
	\\
	\frac{1}{M(v) }& = \frac{\alpha(v)-\phi}{k_\mrm{s}} + \phi c_f
	\label{eq:M_macro},
\end{align}
where now $\alpha(v)$ is associated with the phase--field variable and can be expressed as 
\begin{equation}
	\begin{aligned}
	\alpha(v)
	&= 1- \left[g(v) H (\Tr{\bm{\eps}} ) + H (\Tr{-\bm{\eps}} )\right] \frac{k_\mrm{m}}{k_\mrm{s}} 
	\\
	&= 1- \left[g(v) H (\Tr{\bm{\eps}} ) + H (\Tr{-\bm{\eps}} )\right] (1-\alpha_\mrm{m})
 .
	\end{aligned}
	\label{eq:alpha_v}
\end{equation}

For a fractured material ($v=0$ and $g(v)=0$), the fracture is open when $\Tr{\bm{\eps}} \geq 0$, and Biot's coefficient is enhanced, i.e. $\alpha(v)=1-g(v)+g(v)\alpha_\mrm{m} = 1$. 
On the other hand, when $\Tr{\bm{\eps}} < 0$, the fracture is closed and $\alpha(v)= 1-(1-\alpha_\mrm{m})=\alpha_\mrm{m}$ which corresponds to Biot's coefficient in the bulk material.  Unlike the pre-existing models that consider the degradation of Biot's parameters \cite{yi2020consistent,zhang2023phase},
the value of $\alpha$ expressed in Eq. \eqref{eq:alpha_v} depends both on the value of the phase--field variable and on the specified strain energy decomposition scheme. Therefore, the unilateral effect of fracture can be taken into account if one chooses the V--D decomposition for example.

Additionally, the effective porosity $\phi$ in Eq. \eqref{eq:effective_porosity} is related to the volume fraction of the penny-shaped meso-cracks, i.e., $\varphi$.
And $\varphi $ is approximately 0 because the aspect ratio $X$ of penny-shaped meso-crack is much smaller than 1, and then $\phi \approx \phi_\mrm{m}$ implying that the porosity of RVE--1 remains constant as assumed by~\cite{xie2012micromechanical,ulloa2022variational}. 
Such an assumption is valid for RVEs far from the fracture because the value of $\omega$ is small. 
For nearly fractured RVEs, however, it is invalid because $\omega$ tends to infinity and therefore $\varphi$ tends to 1. 

Since the fracture is diffusely represented in the phase--field model, it is reasonable to introduce a transition function for the porosity with its value ranging from $\phi_\mrm{m}$ to 1.
One of the possibilities is to use the following relationship:
\begin{equation}
	\phi(v) =  1- \left[g(v) H (\Tr{\bm{\eps}} ) + H (\Tr{-\bm{\eps}} )\right] (1-\phi_\mrm{m})
	\label{eq:phi_v}
	,
\end{equation}
which takes a similar expression to $\alpha(v)$ in Eq. \eqref{eq:alpha_v}. 
Therefore, Eq. \eqref{eq:M_macro} reads
\begin{align}
	\frac{1}{M(v) }& = \frac{\alpha(v)-\phi(v)}{k_\mrm{s}} + \phi(v) c_f
	\label{eq:M_macro_f}
	.
\end{align}

As a result, Biot's parameters ($\mbfs{\alpha}$, $M$), porosity ($\phi$) and effective stiffness tensor ($\mathbb{C}_\mrm{eff}$) are now associated with the phase--field variable ($v$).

\subsection{Fluid flow model in cracked porous media}
The mass balance for fractured porous medium and the boundary conditions are given as
\begin{align}
	\label{eq:flow equation}
	\begin{cases}
		\dfrac{\partial  \zeta}{\partial t}+\nabla\cdot( -\dfrac{ \mbfs{K}}{\mu}\nabla p) = Q	& \text{in }  \Omega
		\\
		p = \bar p & \text{on } \partial \Omega_p
		\\
		\mbfs{v}\cdot \mathbf{n} = \mbfs{\bar{q}} & \text{on }  \partial \Omega_q
	\end{cases}
\end{align}
where $Q$ is the source or sink term, $\bar{p}$ is the prescribed pressure, $\mbfs{\bar{q}}$ is the prescribed normal flux (Fig. \ref{fig:object}). $\mbfs{K}$ and $\mu$ are the permeability tensor and fluid viscosity, respectively.

In this work, we consider Darcy flow both in fracture and porous matrix rather than using a diffraction system~\cite{Lee2017_ls} or phase--field calculus~\cite{Chukwudozie2016,Chukwudozie2019} to average or transition between Dacry flow in porous media and Reynolds/Poiseuille flow in fracture~\cite{Heider2017, Ehlers2017}.
Instead, the permeability is enhanced by the evolution of meso-cracks\footnote{Further investigations on this point can be found for example in~\cite{jiang2010experimental,chen2014micromechanical,chen2014experimental} where damage-induced permeability variation was theoretically analyzed from the micromechanical aspect and identified through laboratory experiments. 
Basically, dozens of parameters need to be calibrated experimentally.} as in~\cite{miehe2015minimization, Yoshioka2016}.

We employ an empirical formulation of anisotropic permeability which implicitly takes into account the Poiseuille-type flow \citep{miehe2015minimization,miehe2016phase} in the macroscopic crack, i.e.,
\begin{equation}
	\mbfs{K}= K_\mrm{m} \mbfs{\delta} + (1-v)^{\xi} \frac{w^2}{12}\left(\mbfs{\delta} - \mbfs{n}_{\Gamma} \otimes \mbfs{n}_{\Gamma}\right)
 ,
	\label{eq:Kf}
\end{equation}
where $K_\mrm{m}$ is the isotropic permeability of the porous media, $\xi \ge 1$ is a weighting exponent, and $w$ denotes the fracture aperture.

The weighting exponent has no impact in fully fractured elements where $v=0$ because $(1-v)^\xi = 1$. 
On the contrary, it lessens the permeability enhancement in the transition zone where $v<1$ because $(1-v)^\xi \rightarrow 0$ as $\xi \rightarrow \infty$.
An alternative way to give such a high contrast to the permeability in fully fractured elements ($v=0$) would be to set $w=0$ everywhere except where $v=0$.
However, we find that this would put too much contrast in the permeability and some regularization in the permeability field such as this weighting is necessary for numerical stability.

Because the flow velocity is proportional to its square as in Eq.~\eqref{eq:flow equation}, estimation of $w$ is important and should not be treated lightly~\cite{Yoshioka2020}.
In this study, we use an approach proposed by \cite{miehe2015minimization}, i.e.,
\begin{equation}
	w = h_e \Tr{\bm{\eps}}
 ,
 \label{eq:width}
\end{equation}
where $h_e$ is the characteristic length which is taken as the element size in the present work. 

Furthermore, the crack normal vector $\mbf{n}_{\Gamma}$ needs to be resolved for Eq.~\eqref{eq:Kf}.
Bourdin et al.~\cite{Bourdin2012} approaximated $\mbf{n}_{\Gamma}$ from the gradient of the phase--field (e.g. $\mbf{n}_{\Gamma} \approx \nabla v/|\nabla v|$), which has been followed by many~\cite{Mikelic2015_NonLin,miehe2016phase, santillan2017phase, guo2019modelling,Kienle2021, costa2022multi}.
However, as pointed out by~\cite{Chukwudozie2019}, the gradient of the phase--field deviates from the normal direction near the crack tip.
Another issue is that $\nabla v$ is not defined at fully fractured elements where $v=0$, which requires an additional operation such as $L$-2 projection~\cite{Kienle2021}.
In this study, we estimate the normal vector from the principal strains at integration points.
We first decompose a strain tensor as
\begin{equation}
    \label{eq:eigen_strain}
    \mbfs{\eps} = \sum_{i=1}^3 \eps_i \mbfs{e}_i
    ,
\end{equation}
where $\eps_i$ are the principal strains with $\eps_1>\eps_2>\eps_3$ and $\mbf{e}_i$ are the associated eigenvectors.
Then we assign 
\begin{equation}
\mbfs{n}_{\Gamma} = \mbfs{e}_1
.
\end{equation}
One advantage of this approach is that it is easy to parallelize the computation because it is performed locally at each integration point while the line integral~\cite{Bourdin2012} or level--set~\cite{Lee2017_ls} based approach requires global phase--field profile.
This approximation is sufficiently accurate for hydraulic fracture application as demonstrated in Section~\ref{sec:verification}.
However, for other applications such as frictional sliding~\cite{fei2020phase}, one may require more involved approaches based on image processing techniques~\cite{ziaei2016identifying, bryant2021phase,xu2023reconstruct}.

\section{Numerical Implementation}
\label{sec:num_impl}
The regularized crack evolution problem can be regarded as the minimization of functional $\mathcal{F}_{\ell}$ (expressed in Eq. \eqref{eq:F_l}) with respect to the variables $\mbfs{u}$ and $v$ with the damage irreversible condition, i.e.,
\begin{equation}
	\label{eq:glob_min}
	(\mbfs{u}, v) =
	\mathop {\argmin} \limits_{\mbfs{u}, v} \left\lbrace  \mathcal{F}_{\ell}(\mbfs{u},v; p) : \mathbf{u} \in \mathcal{U}, v \in \mathcal{V}(t_i) \right\rbrace,
\end{equation}
where $\mathcal{U}$ is the kinematically admissible displacement set:
\begin{equation}
	\label{eq:_adm}
	\mathcal{U} = \left\{\mbfs{u} \in H^1 (\Omega) : \mbfs{u} = \bar{\mbfs{u}} \quad \mrm{on} \quad \partial \Omega_u \right\}
	,
\end{equation}
and the admissible set of $v$ is
\begin{equation}
	\label{eq:v_adm}
	\mathcal{V}(t_i) = \left\{v \in H^1 (\Omega) : 0 \le v(x, t_i) \le v(x, t_{i-1}) \le 1 \ \forall x \ \mrm{s.t.} \  v(x, t_{i-1}) \le v_\mrm{irr} \right\}
	,
\end{equation}
where $v_\mrm{irr}$ is a variable that controls the irreversibility.
Rather than imposing a strict irreversibility such as $v(x, t_{i}) < v(x, t_{i-1})$, we apply the irreversibility only where $v(x, t_{i-1}) \le v_\mrm{irr} $ as proposed by~\cite{bourdin2007numerical,burke2013adaptive}.
Note that if we choose a large $v_\mrm{irr}$ (i.e. $v_\mrm{irr} \rightarrow 1$), it is essentially the same as the strict irreversibility.
On the other hand, if $v_\mrm{irr}$ is small (i.e. $v_\mrm{irr} \rightarrow 0$), then the damage is reversible until complete failure.

Eq.~\eqref{eq:glob_min} is solved via an alternate minimization scheme~\cite{Bourdin2000}, taking advantage of the convexity of $\mathcal{F}_\ell$ with respect to $\mbfs{u}$ and $v$, which minimizes $\mathcal{F}_\ell$ first with respect to $\mbfs{u}$ while fixing $v$ and then minimizes $\mathcal{F}_\ell$ with respect to $v$ while fixing $\mbfs{u}$ with the irreversible condition.
Accordingly, we can obtain:
\begin{align}
	\label{eq: phase-field}
	&\nabla \cdot \left[ \mathbb{C}_\mrm{eff}(v) : \mbfs{\varepsilon}(\mbfs{u}) -  {\alpha}(v) p \mbfs{\delta} \right]  + \mbfs{b} = \mathbf{0} \\
	\label{eq: displacement}
	&2v\psi^+_{e} + \frac{p^2}{2} \frac{\partial {1/M(v)}}{\partial v}  + \frac{G_c}{4c_n} \left( -\frac{n (1-v)^{n-1}}{\ell_s} + 2\ell_s \Delta v\right) = 0
\end{align}
with boundary conditions $\mbfs{\sigma}\cdot n=\bar{\mbfs{t}}$ on $\partial \Omega_t$ and $\nabla d \cdot \mbfs{n} =0 $ on $\partial \Omega$. 

The derivative $\frac{\partial {1/M(v)}}{\partial v}$ in the second term of Eq. \eqref{eq: displacement} can be written according to Eq. \eqref{eq:M_macro_f} as
\begin{align}
    \frac{\partial {1/M(v)}}{\partial v} &= \frac{\partial}{\partial v}\left[\frac{g(v)H (\Tr{\bm{\eps}} ) (\alpha_\mrm{m}-\phi_\mrm{m})}{k_\mrm{s}} - g(v)H (\Tr{\bm{\eps}} )(1-\phi_\mrm{m})c_f +c_f \right] \nonumber \\
    &=2vH (\Tr{\bm{\eps}} ) \left[\frac{\alpha_\mrm{m}-\phi_\mrm{m}}{k_\mrm{s}} - (1-\phi_\mrm{m}) c_f\right]
    .
\end{align}

\subsection{Weak form and solution strategy of the three-field problem}

The weak forms of Eqs. \eqref{eq: phase-field} and \eqref{eq: displacement} can be obtained by multiplying with the weight functions $\mbfs{w}_u$ and $w_{v}$ and integrating over the problem domain. Accordingly, we have
\begin{align}
	\label{eq: phase-field weak form}
	&\int _{\Omega} \nabla \mbfs{w}_u \cdot \left[ \mathbb{C}_\mrm{eff}(v) : \mbfs{\varepsilon}(\mbfs{u}) -  {\alpha}(v) p \mbfs{\delta} \right] \mathd V  - \int _{\Omega} \mbfs{b} \cdot \mbfs{w}_u \mathd V -  \int _{\partial\Omega_t} \mbfs{\bar{t}} \cdot \mbfs{w}_u \mathd S= \mathbf{0} \\
	\label{eq: displacement weak form}
	&\int _{\Omega} w_{v}\left[2v\psi^+(\boldsymbol{\varepsilon},v) +\frac{p^2}{2} \frac{\partial {1/M(v)}}{\partial v} \right] \mathd V - \int _{\Omega}w_{v} \frac{G_c}{4c_n} \frac{n (1-v)^{n-1}}{\ell_s} \mathd V -\int _{\Omega} \frac{G_c}{2c_n} \ell_s \nabla w_{v} \cdot \nabla v \mathd V = 0 
\end{align}

Additionally, the weak form of the  flow equation has been presented in Eq. \eqref{eq:flow equation}. 
For the weak form of the flow equation Eq.~\eqref{eq:flow equation}, multiplying the weight function $w_p$ yields
\begin{equation}
	\int_{\Omega}  w_p \frac{\partial \zeta(\mbfs{\varepsilon}, v)}{\partial t}\, \mathd V 
	+\int_{\Omega} \frac{ \mbfs{K}}{\mu}\nabla p \cdot \nabla w_p\,\mathd V 
	=  \int_{\Omega}\, w_p Q \mathd \Omega  	
	-\int_{\partial \Omega_q} w_p \bar{\mbfs{q}}\, \mathd S
	\label{eq:combined flow2}
\end{equation}

We employ a staggered scheme in the present work to solve this $(\mbfs{u}-p)-v$ three-field problem. Globally, the $(\mbfs{u}-p)$ problem will be solved first in a staggered manner following the fixed-stress split iterative coupling strategy \cite{Kim2011}, and then the phase-field $v$ problem is solved and iterated with $(\mbfs{u}-p)$ problem until three fields converge.  
Additionally, to bound the solution space for $v$ (i.e., $v(x,t_i) \in [0,1] $ if $v(x,t_{i-1}) > v_\mrm{irr}$, and $v(x,t_i) \in [0,v_\mrm{irr}] $ otherwise), we apply the variational inequality solver of PETSc \cite{balay2019petsc}.
This irreversibly has been enforced alternatively by an augmented Lagrangian approach~\cite{Wheeler2014}, a history variable~\cite{Miehe2010, Ambati2014}, or a penalty based method~\cite{Gerasimov2019}.

The proposed model is implemented in an open source code, OpenGeoSys~\cite{ogs:6.4.3}. 
Further information and examples are freely accessible at \url{https://www.opengeosys.org/}.

\subsection{ Fixed stress split for solving {u}-p  problem}
The backward Euler scheme tackles the time derivatives in Eq. \eqref{eq:combined flow2}. The time derivative of the first term in Eq. \eqref{eq:combined flow2} yields
\begin{align}
	\label{eq:time derivative}
	\int_{\Omega}  w_p \frac{\partial \zeta(\mbfs{\varepsilon}, v)}{\partial t}\,\mathd V & = \int_{\Omega}  w_p \frac{\partial}{\partial t} \left[ \alpha(v)\mbfs{\delta}:\mbfs{\varepsilon}+\frac{p}{M(v)}\right] \psi \mathd V \nonumber \\
	&=  \int_{\Omega}  w_p \alpha(v) \frac{\partial(\mbfs{\delta}:\mbfs{\varepsilon})}{\partial t} \mathd V  +  \int_{\Omega}  w_p \frac{1}{M(\nu)}\frac{\partial p}{\partial t}  \mathd V
 .
\end{align}
As mentioned, the phase-field $v$ is frozen during the iterative solution of $\mbfs{u}-p$ problem, and therefore the variation of $v$ is not involved in Eq. \eqref{eq:time derivative}.

At step $i$, the back Euler scheme gives rise to
\begin{equation}
	\label{eq:back Euler}
	\begin{cases}
		\dfrac{\partial(\mbfs{\delta}:\mbfs{\varepsilon})}{\partial t} = \dfrac{\mbfs{\delta}:(\mbfs{\varepsilon}^{i} - \mbfs{\varepsilon}^{i-1})}{\Delta t} \\
		\dfrac{\partial p}{\partial t} = \dfrac{p^{i} - p^{i-1}}{\Delta t} 
	\end{cases}
\end{equation}
where $\Delta t$ denotes the time increment.

Substituting Eq. \eqref{eq:back Euler} into Eq. \eqref{eq:time derivative} gives
\begin{equation}
	\label{eq:time derivative2}
	\int_{\Omega}  w_p \frac{\partial \zeta(\mbfs{\varepsilon}, v)}{\partial t}\,\mathd V
	=  \int_{\Omega}  w_p \alpha(v) \frac{\mbfs{\delta}:(\mbfs{\varepsilon}^{i} - \mbfs{\varepsilon}^{i-1})}{\Delta t} \mathd V 
	+  \int_{\Omega}  w_p \frac{1}{M(v)} \frac{p^{i} - p^{i-1}}{\Delta t} \mathd V
 .
\end{equation}
Substituting Eq. \eqref{eq:time derivative2} into Eq. \eqref{eq:combined flow2} then yields
\begin{multline}
	\int_{\Omega}  w_p \frac{1}{M(v)}  \frac{p^{i} - p^{i-1}}{\Delta t}  \psi \mathd V 
	+\int_{\Omega}  \frac{ \mbfs{K}}{\mu}\nabla p\cdot \nabla w_p\,\mathd V  \\
	=  \int_{\Omega}\, w_p Q \mathd \Omega  - \int_{\Omega}  w_p \alpha(v) \frac{\mbfs{\delta}:(\mbfs{\varepsilon}^{i} - \mbfs{\varepsilon}^{i-1})}{\Delta t} \mathd V 
	-\int_{\partial \Omega_q} w_p \bar{\mbfs{q}} \, \mathd S
 .
	\label{eq:combined flow_time discretization}
\end{multline}

To ensure the stability of the flow field, Eq. \eqref{eq:combined flow_time discretization} will be solved following the fixed stress split scheme proposed in~\cite{Kim2011}. To this end, we define $j$ as the coupling iteration step for $\mbfs{u}-p$ solution and ensures the stress in iteration~$j$ is fixed:
\begin{equation}
	\label{eq:mean_stress}
	k_\mrm{eff} (v) \mbfs{\delta}:\mbfs{\varepsilon}^{i,j}  - \alpha(v) p^{i,j} = k_\mrm{eff} (v) \mbfs{\delta}:\mbfs{\varepsilon}^{i,j-1}  - \alpha(v) p^{i,j-1}
 ,
\end{equation}
which gives rise to the following relation as
\begin{equation}
	\label{eq:epsilon_v}
	\mbfs{\delta}:\mbfs{\varepsilon}^{i,j}  =  \mbfs{\delta}:\mbfs{\varepsilon}^{i,j-1}  - \frac{\alpha(\nu)}{k_\mrm{eff} (v)}  \left(p^{i,j-1}-p^{i,j} \right)
 ,
\end{equation}
where $k_\mrm{eff}$ is the effective bulk modulus defined in Eq. \eqref{eq:relation of stiffness}.

Plugging Eq. \eqref{eq:epsilon_v} into Eq. \eqref{eq:combined flow_time discretization} gives
\begin{multline}
	\int_{\Omega}  w_p \frac{ 1}{M(v)} \frac{p^{i,j} - p^{i-1}}{\Delta t}  \mathd V + \int_{\Omega} w_p \frac{\alpha^2(v)}{k_\mrm{eff}(v)} \frac{p^{i,j} - p^{i,j-1}}{\Delta t} \mathd V  
	+\int_{\Omega} \frac{ \mbfs{K}}{\mu}\nabla p^{i,j} \cdot \nabla w_p\,\mathd V 
	\\ =  \int_{\Omega}\, w_p Q \mathd \Omega  	
	-\int_{\partial \Omega_q} w_p \bar{\mbfs{q}}\, \mathd S  - \int_{\Omega}  w_p \alpha(v) \frac{\mbfs{\delta}:(\mbfs{\varepsilon}^{i,j-1} - \mbfs{\varepsilon}^{i-1})}{\Delta t} \mathd V  
 .
	\label{eq:fix stress split}
\end{multline}
\subsection{Spatial discretization}

An isoparametric element is used to discretize the problem in spatial, and bilinear shape functions are employed for spatial interpolation. Therefore, the approximations of $\mbfs{u}$, $v$ and $p$ can be expressed as
\begin{equation}
	\begin{array}{cccc}
		\boldsymbol{u}\approx \mbf{N}_{u}\boldsymbol{\hat{u}}, & \boldsymbol{w%
		}_{u}\approx \mbf{N}_{u}\boldsymbol{\hat{w}}_{u}, &
		\nabla \boldsymbol{u}\approx \mbf{B}_{u}\boldsymbol{\hat{u}}, &
		\nabla \boldsymbol{w}_{u}\approx \mbf{B}_{u}\boldsymbol{\hat{w}}_{u}%
		\\
		v\approx \mbf{N}_{v}\hat{v}, & w_{v}\approx \mbf{N}_{v}\hat{w}%
		_{v},  &
		\nabla v\approx \mbf{B}_{v}\hat{v}, & \nabla w_{v}\approx \mbf{%
			B}_{v}\hat{w}_{v}%
		\\
		p\approx \mbf{N}_{p}\hat{p}, & w_{p}\approx \mbf{N}_{p}\hat{w}%
		_{p},  &
		\nabla p\approx \mbf{B}_{p}\hat{p}, & \nabla w_{p}\approx \mbf{%
			B}_{p}\hat{w}_{p}%
	\end{array}
	\label{eq:Galerkin}
\end{equation}%
Substituting the Galerkin approximation \eqref{eq:Galerkin} into Eqs. \eqref{eq: phase-field weak form}, \eqref{eq: displacement weak form} and \eqref{eq:fix stress split} gives the discrete system as follows
\begin{equation}
	\label{eq: phase-field discrete form}
	\int _{\Omega} \mbf{B}_u^{\mrm{T}} \cdot \left[ \mathbb{C}_\mrm{eff} : \mbf{B}_u\mbfs{\hat{u}} -  {\alpha}(v) p \mbfs{\delta}\right] \mathd V  - \int _{\Omega} \mbf{N}_u^{\mrm{T}}  \mbfs{b} \mathd V -  \int _{\partial\Omega_t} \mbf{N}_u^{\mrm{T}}  \mbfs{\bar{t}} \mathd S= \mathbf{0}
\end{equation}
\begin{equation}
	\label{eq: displacement discrete form}
	\int _{\Omega} \mbf{N}_{v}^{\mrm{T}} \left[2 \mbf{N}_{v}\hat{v}\psi^+(\boldsymbol{\varepsilon},v) +\frac{p^2}{2} \frac{\partial {1/M(v)}}{\partial v} \right] \mathd V - \int _{\Omega}\mbf{N}_{v}^{\mrm{T}} \frac{G_c}{4c_n} \frac{n (1-\mbf{N}_{v} \hat{v})^{n-1}}{\ell_s} \mathd V 
	-\int _{\Omega} \frac{G_c}{2c_n} \ell_s \mbf{B}_{v}^{\mrm{T}} \mbf{B}_{v} \hat{v} \mathd V = 0 
\end{equation}
\begin{multline}
	\label{eq: pressure discrete form}
	\int_{\Omega} \mbf{N}_{p}^{\mrm{T}}\mbf{N}_{p} \left\{\frac{1}{M(v)} \frac{\hat{p}^{n,i} - \hat{p}^{n-1}}{\Delta t} + \frac{\alpha^2(v)}{k_\mrm{eff}(v)} \frac{\hat{p}^{i,j}-\hat{p}^{i,j-1}}{\Delta t} \right\} \mathd V
	+\int_{\Omega}  \frac{\mbfs{K}}{\mu} \mbf{B}_{p}^{\mrm{T}} \mbf{B}_{p}^{\mrm{T}} \hat{p}^{i,j} \mathd V \\
	= \int_{\Omega} \mbf{N}_{p}^{\mrm{T}} \left[Q - \alpha(v) \frac{\mbfs{\delta}:(\mbfs{\varepsilon}^{i,j-1}-\mbfs{\varepsilon}^{i-1})}{\Delta t}
	\right]\mathd V - \int_{\partial \Omega_q} \mbf{N}_{p}\mbfs{\bar{q}} \mathd S
\end{multline}
where the permeability is updated using Eq.~\eqref{eq:Kf} before each $p$ solution.

\section{Verification}
\label{sec:verification}
In this section, we verify our implementation against known closed form solutions.
The first two problems test transient responses of the hydro-mechanical part of the implementation.
The third problem tests the fracture aperture computation for a static crack.
Finally, the last problem verifies propagation of a hydraulic fracture in a toughness dominated regime in plane strain.
\subsection{One-dimensional consolidation problem} 
\begin{figure}[htp!]
    \centering
    \includegraphics[scale=0.5]{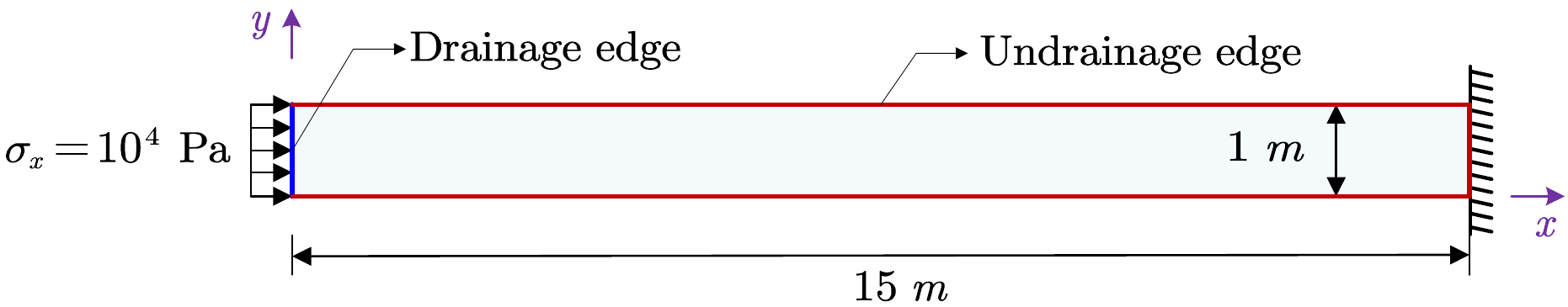}
    \caption{Geometry and boundary conditions of the  one-dimensional consolidation problem}
    \label{fig:geometry_1}
\end{figure}

We start with the 1D consolidation problem to verify the feasibility of the hydro-mechanical modules. The geometry and boundary of this 1D consolidation problem are illustrated in Fig. \ref{fig:geometry_1}. A constant stress $\sigma_x=2\times10^6$ Pa is applied on the left boundary and the right boundary is fixed, $u(L,t) = 0$. 
All the boundaries are assigned with no-flow (undrained) boundary except for the left edge that is drained, i.e., $p=0$ Pa. 
Table \ref{tab:case1 parameter} lists the material parameters used in the simulation.

\begin{table}[htbp]
    \centering
    \begin{tabular}{ccc}
    \hline
      Parameter   &  Value & Unit\\
      \hline
     Young's modulus for porous material ($E$)  &  $3.0\times10^8$ & Pa \\
     Poisson's ratio for porous material ($\nu$)  & 0 & --\\
     Porosity for porous material ($\phi_\mrm{m}$) & 0.3 & --\\
    Fluid compressibility ($c_\mrm{f}$) & $1\times10^{-9}$ & Pa \\
     Permeability ($K_\mrm{m}$) & $2\times 10^{-12}$ & m$^2$ \\
    Fluid viscosity ($\mu$) & $10^{-3}$ & Pa$\cdot$s \\
     Biot coefficient ($\alpha_\mrm{m}$) & 1 & --\\
       \hline
    \end{tabular}
    \caption{Parameters for the 1D consolidation problem}
    \label{tab:case1 parameter}
\end{table}
The analytical solutions of this problem are given as \cite{wang2000theory}
\begin{align*}
  &  p(x,t)=\frac{4d\sigma_x}{\pi}\sum_{m=1}^{\infty}\left\{ \frac{1}{2m+1}\mrm{exp}\left(-\frac{(2m+1)^2\pi^2}{4L^2}ct\right)\mrm{sin}\left(\frac{(2m+1)\pi x}{2L}\right)\right\} \\
 &   u(x,t)=c_m d \sigma_x \left[ L-x-\frac{8L}{\pi^2}\sum_{m+1}^{\infty}\left\{ \frac{1}{(2m+1)^2}\mrm{exp}\left(-\frac{(2m+1)^2\pi^2}{4L^2}ct\right)\mrm{cos}\left(\frac{(2m+1)\pi x}{2L}\right)\right\} \right] + b\sigma_x(L-x)
\end{align*}
where $L=15$ m, $d=\frac{a-b}{a\alpha_\mrm{m}}$, $c = \frac{K_\mrm{m}}{(a \alpha_\mrm{m}^2+S)\mu}$ and $c_m=\frac{a-b}{d}$ with $a=\frac{(1+\nu)(1-2\nu)}{E(1-\nu)}$, $b=\frac{a}{1+a\alpha_\mrm{m}^2/S}$ and $S=\frac{1}{N}+\phi_\mrm{m} c_\mrm{f}$. 

\begin{figure}[htp!]
    \centering
    \subfigure[Pressure distribution]{\includegraphics[scale=0.5]{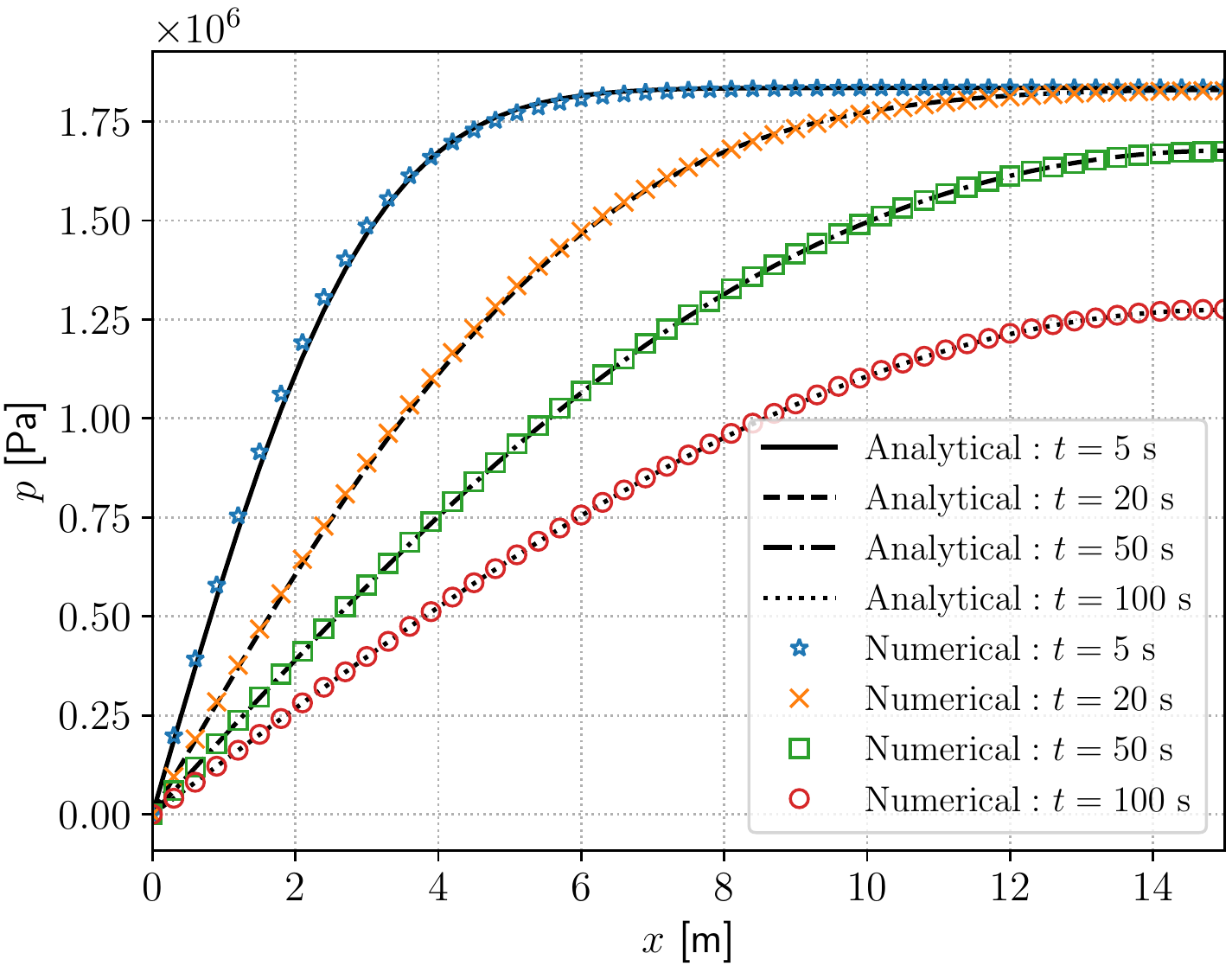}}
    \subfigure[Displacement distribution]{\includegraphics[scale=0.5]{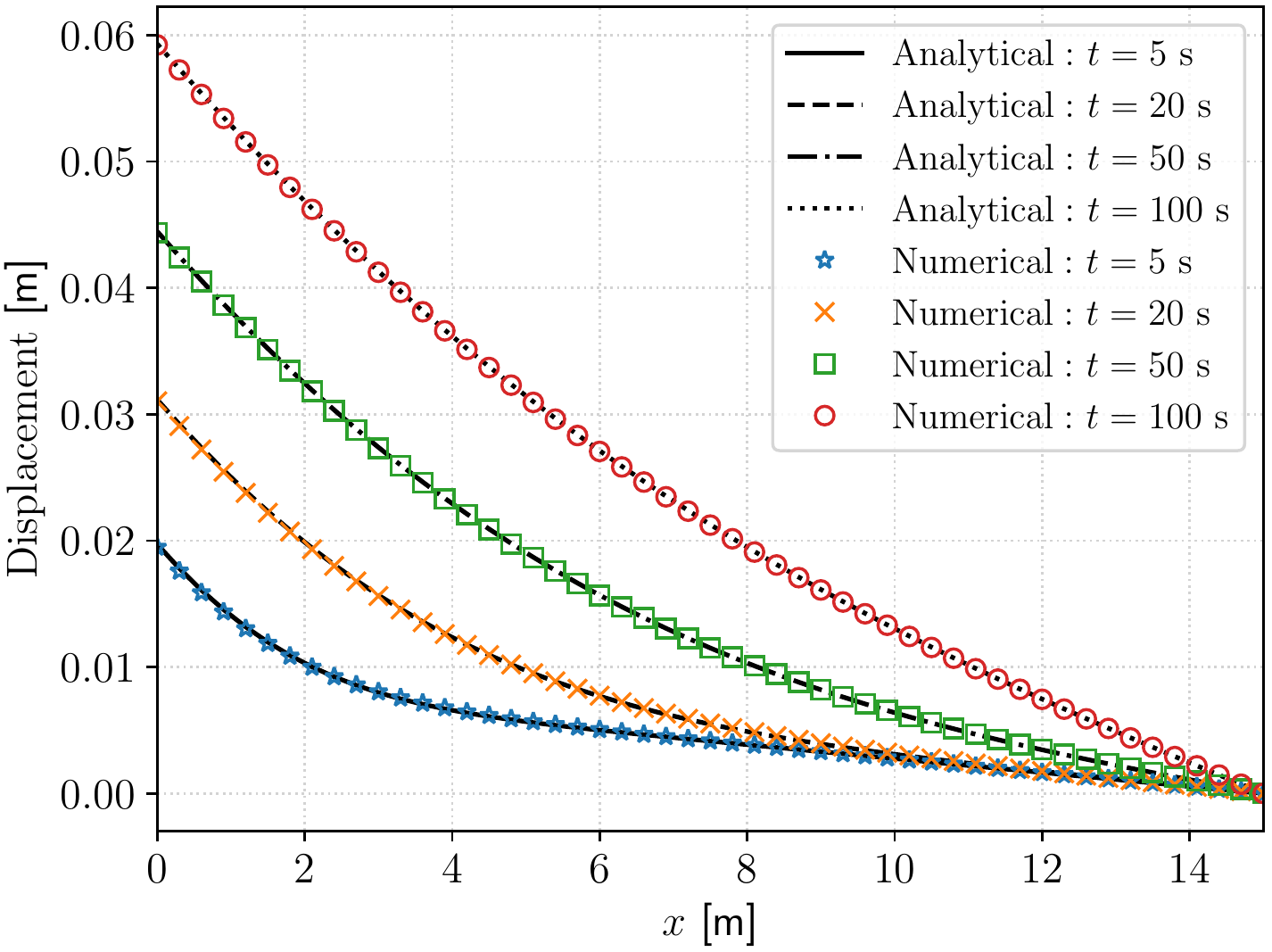}}
    \caption{Comparison of analytical and numerical results of pressure and displacement distributions in the one-dimensional consolidation problem.}
    \label{fig:pressure}
\end{figure}

In this case, damage evolution is not involved, and only displacement and pressure fields are solved. The time increment is taken as $\Delta t=1$ s. 
We can observe a good agreement between the analytical and numerical results for pressure and displacement distributions (Fig. \ref{fig:pressure}).

\subsection{Pressure distribution in a single crack}
\begin{figure}[htp!]
    \centering
    \includegraphics[scale=0.7]{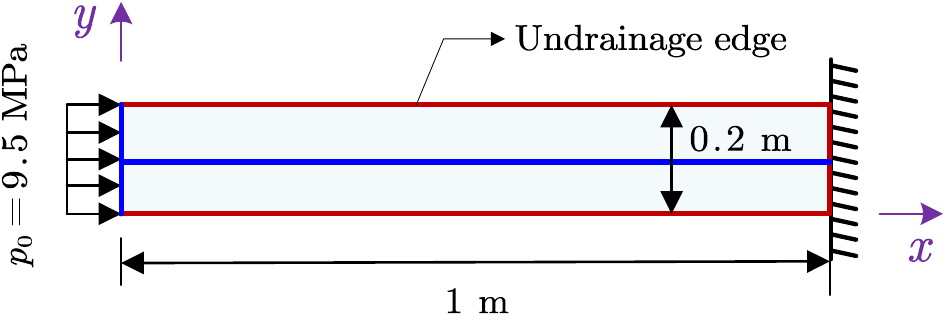}
    \caption{Geometry and boundary conditions for a rock sample with a single crack.}
    \label{fig:geometry_2}
\end{figure}

In this section, the variation of pressure distribution in a single crack with a constant pressure boundary is investigated. The geometry and boundary conditions of this example are illustrated in Fig. \ref{fig:geometry_2}. A Dirichlet condition for pressure, i.e., $p_0=9.5$ MPa is applied on the left edge, and other edges are impermeable. The parameters used in the simulation are listed in Table \ref{tab:case2 parameter} following \cite{zhou2019phase}. The initial crack is represented by setting the initial phase-field value to 0 in the crack domain. 

\begin{table}[htbp]
    \centering
    \begin{tabular}{ccc}
    \hline
      Parameter  & Value & Unit\\
      \hline
       Porosity ($\phi_\mrm{m}$)  & $2.0\times 10^{-5}$ & --\\
       Fluid compressibility ($c_\mrm{f}$) &  $4.55\times10^{-10}$  & Pa \\
       Permeability ($K_\mrm{m}$) &   $10^{-20}$ & m$^2$ \\
       Fluid viscosity ($\mu$) &   $10^{-3}$ & Pa$\cdot$s \\
        Biot coefficient ($\alpha_\mrm{m}$) &  $2\times 10^{-5}$ & --\\
        Regularization parameter ($\ell$) &  $0.01$  & -- \\
       \hline
    \end{tabular}
    \caption{Parameters for the 1D consolidation problem}
    \label{tab:case2 parameter}
\end{table}

The analytical solution is given as \cite{yang2018hydraulic}
\begin{equation*}
    \frac{p(x,t)}{p_0} = 1 + \frac{4}{\pi}\sum_{m=0}^{\infty}\left[ \mrm{exp}\left( -(2m+1)^2(t_\mrm{D}/4)\pi^2\right) \mrm{cos}\left( \frac{(2m+1)\pi}{2}\zeta\right)\frac{(-1)^{m+1}}{2m+1}\right]
\end{equation*}
where $\zeta=\frac{L-x}{L}$ and $L=1$ m in this case, and $t_\mrm{D}$ is a dimensionless time defined as
\begin{equation*}
    t_\mrm{D} = \frac{w^2 t}{12\mu c_\mrm{f} L^2}
\end{equation*}
where $c_\mrm{f}$, $\mu$ and $w$ are the fluid compressibility, the fluid viscosity and the aperture of the central crack, respectively. 

In this case, the fracture width $w$ is only activated in the crack domain ($v=0$), and its value is a constant set to ensure $\frac{w^2}{12 \mu L^2}=0.5$. The time increment is then chosen as $\Delta t=0.01$ s, and the total time is 0.5 s.

\begin{figure}[htp!]
    \centering
    \includegraphics[scale=0.5]{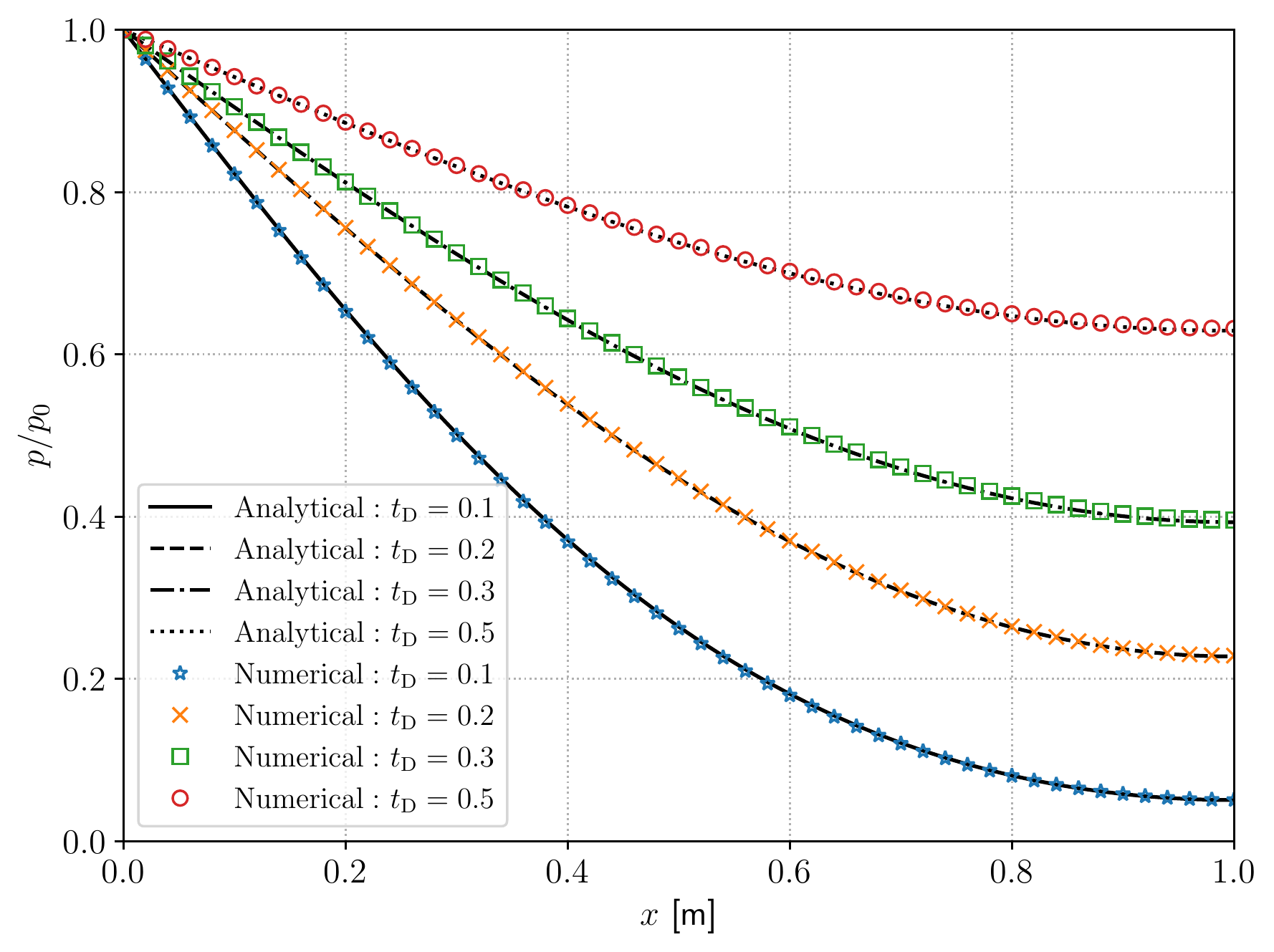}
    \caption{Comparison of analytical and numerical results of pressure distribution along the center of the preexisting crack.}
    \label{fig:p_vs_p0}
\end{figure}

Only the pressure field and phase field are solved, and the mechanics-related part is skipped. The comparison of the numerical and analytical results is shown in Fig. \ref{fig:p_vs_p0}. The numerical results closely match the analytical results.

\subsection{A single-cracked plate subjected to constant internal pressure} \label{sec:sneddon}
\begin{figure}[htpb!]
    \centering
    \includegraphics[scale=0.6]{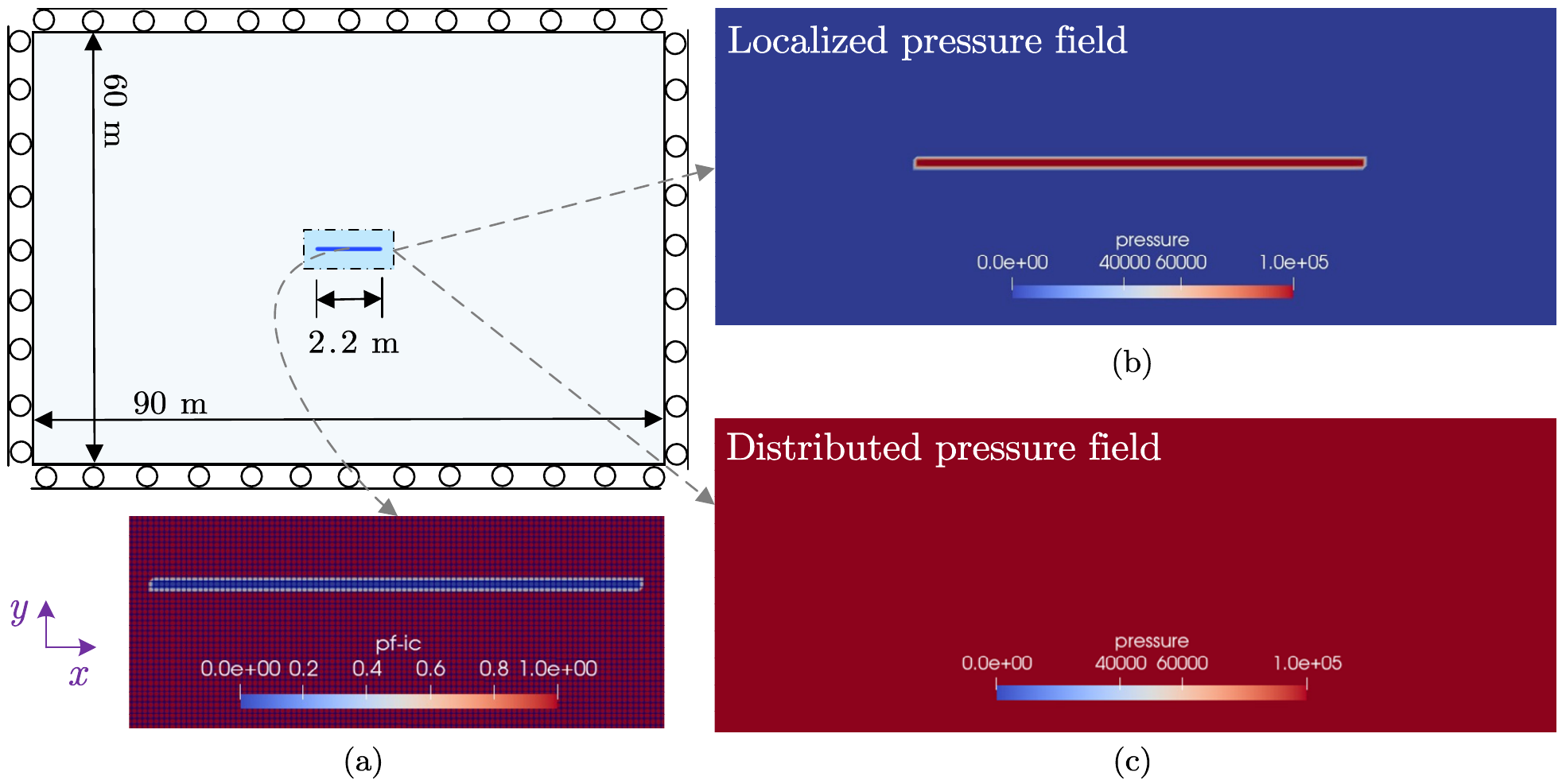}
    \caption{The numerical model of the single-cracked plate subjected to constant internal pressure in which two types of pressure fields are considered, i.e., the localized and distributed.}
    \label{fig:numerical_case3}
\end{figure}

In this section, we test our crack opening computation through Sneddon's internally pressurized line crack $[-a_0, a_0] \times \lbrace 0 \rbrace$ in an infinite domain~\cite{Sneddon-Lowengrub-1969a}.
Though this example has become one of the standard tests for phase--field hydraulic fracture models~\cite{Bourdin2012,Wheeler2014,Heister2015,Lee2017_ls,santillan2017phase, Gerasimov2019,Yoshioka2020,Kienle2021,zhou2018phase}, not all the models recovered the theoretical crack opening accurately.
The computation geometry and boundary conditions are depicted in Fig.~\ref{fig:numerical_case3}a. The crack length is taken as $2a_0=2.2$ m. Here, the initial static crack is represented by assigning $v=0$ to the corresponding elements (i.e., `pf-ic' field in Fig .\ref{fig:numerical_case3}a). 

The analytical solution of crack opening displacement under pressurization is given as~\cite{Sneddon-Lowengrub-1969a}
\begin{equation}
\label{eq: u+}
    u^+(x,0) = \frac{2pa_0}{E^{\prime}} \sqrt{1-\left(\frac{x}{a_0}\right)^2}
    ,
\end{equation}
where $E^{\prime}=E/(1-\nu^2)$, $x$ is the distance to the crack center and $p$ is the internal pressure taking as $p=10^5$ Pa.  
As the solution is based on the assumption of elastic material, we take the material parameters as follows: Young's modulus $E=1.7\times10^{10}$, Poisson's ratio $\nu=0.2$ and  $\alpha_\mrm{m}=0$, $\phi_\mrm{m}=0$.

With the diffused representation of crack by the phase-field model, we need to account for the extra energy near the crack tip as discussed in~\cite{freddi2019fracture, Yoshioka2020}.
Accordingly, a closed-form expression for the effective crack length is given as~\cite{Yoshioka2020}
\begin{equation}
    a_\mrm{eff} = a_0 \left( 1 + \frac{\pi \ell/4}{a_0 \left( h/4c_n \ell +1\right) } \right) 
    .
\end{equation}
In the following computations, the crack lengths are normalized by $a_\mrm{eff}$.
Note that $a_\mrm{eff}$ approaches $ a_0$ as $\ell \to 0$. 
We used \ATone{} in the following examples because \ATone{} was reported to recover crack openings more accurately~\cite{Yoshioka2020}.

For the prescribed pressure field, we consider the localized pressure field (Fig.~\ref{fig:numerical_case3}b) where the pressure is elevated ($p=10^5$ Pa) only in the fully broken elements ($v=0$).
This localized pressure field is supposed to represent the actual hydraulic fracture situation where the fluid pressure within the fracture is distinctively higher than the surrounding.

\subsubsection*{Verification of our proposed model}
\begin{figure}[htpb!]
    \centering
    \subfigure[Given the ratio $\ell/h=4$]{\includegraphics[scale=0.5]{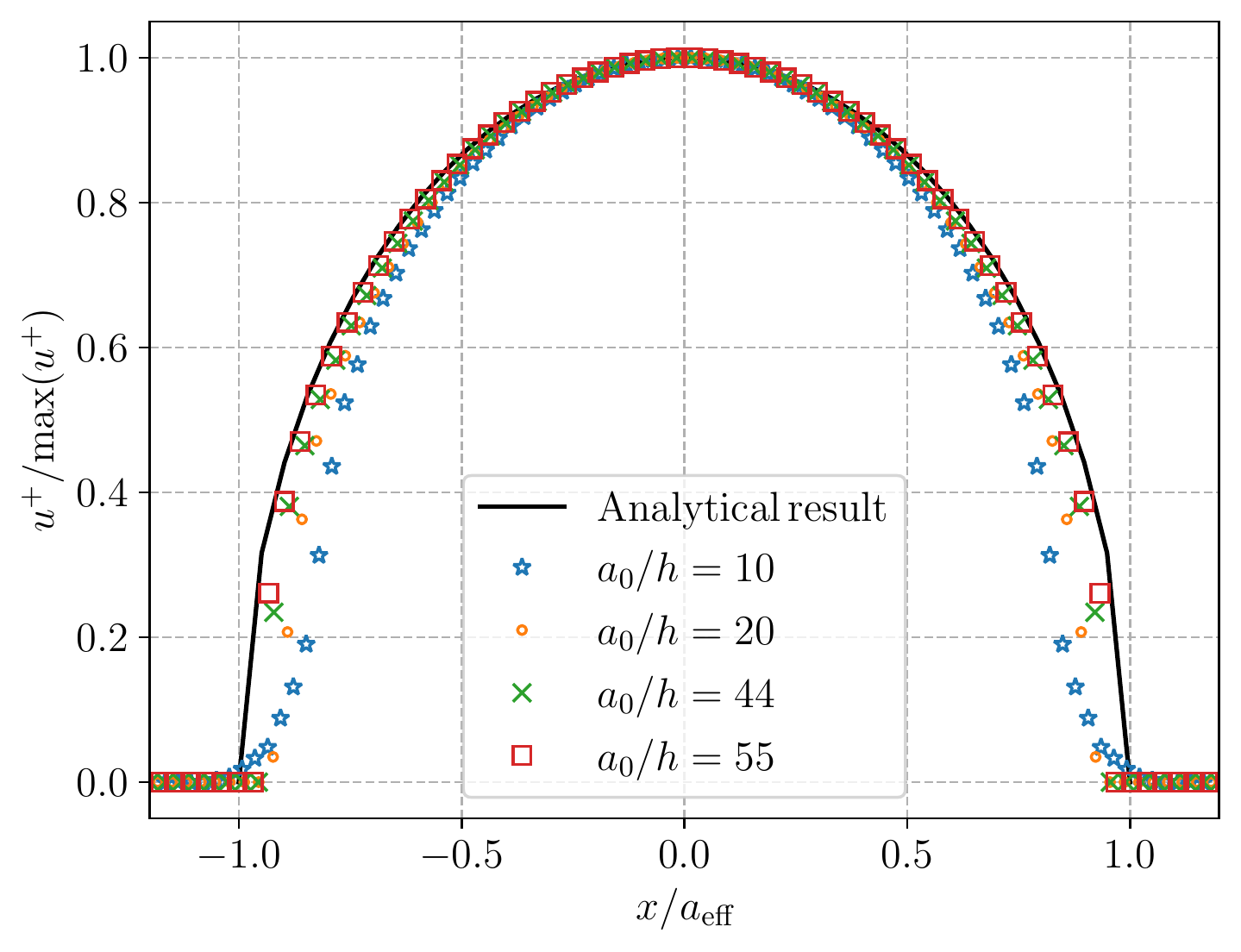} \label{fig:fix ratio}}
    \subfigure[Given the mesh size $h=0.025$ m]{\includegraphics[scale=0.5]{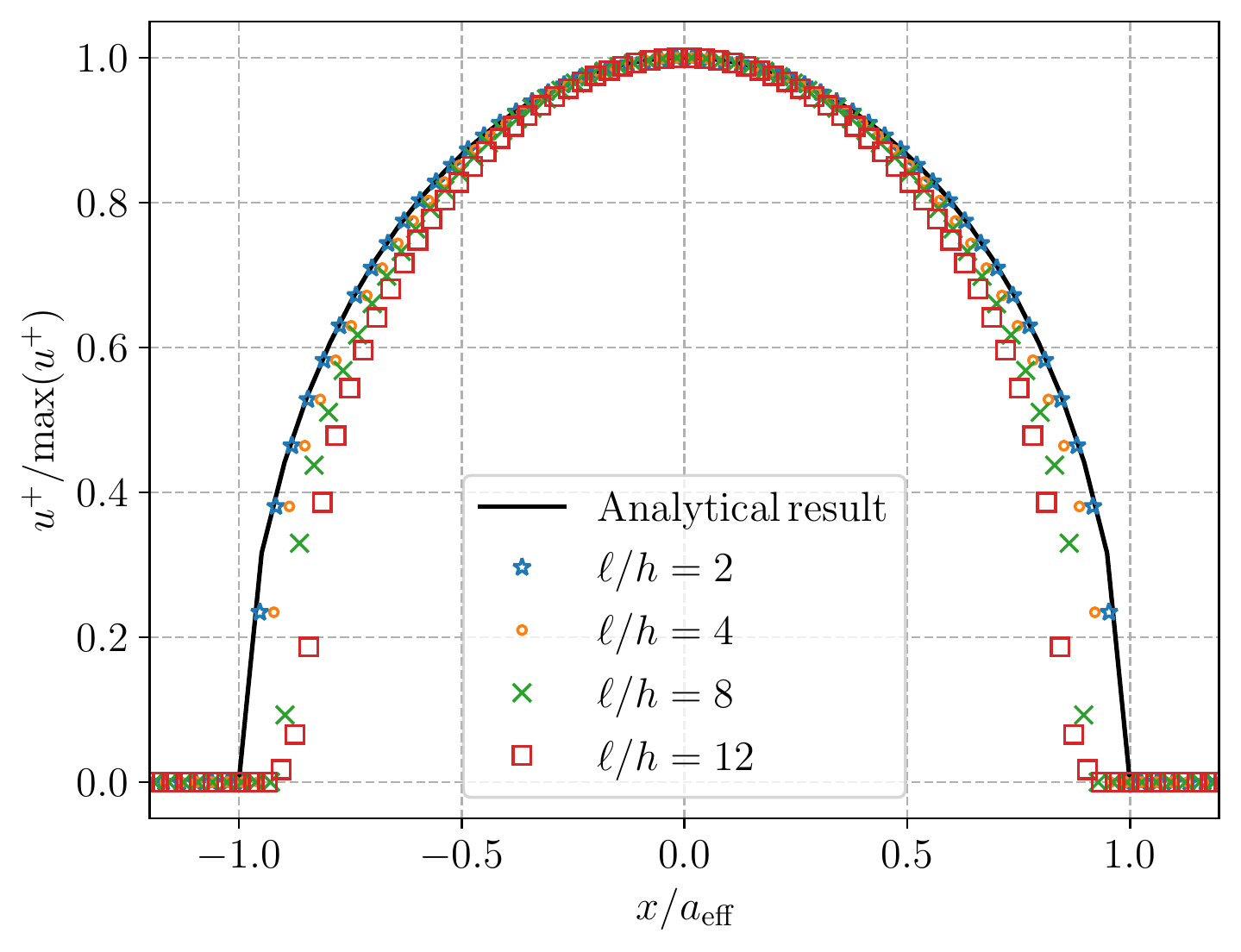} \label{fig:fix mesh size}}
    \caption{Crack opening displacements on various mesh resolutions and length scale parameters for \ATone{} phase-field model with diffused porosity.}
    \label{fig:COD_mesh}
\end{figure}

First, we ascertain the influences of mesh resolution and length scale parameter.
Fig.~\ref{fig:fix ratio} shows the impacts of mesh resolution where we keep the ratio of length scale parameter $\ell$ and  minimum element size (denoted by $h$) as $\ell/h=4$.
The results indicate that roughly 20 elements over the span of the initial crack geometry will give an acceptable result, which is consistent with the conclusion from~\cite{cusini2021simulation}. 
Fig. \ref{fig:fix mesh size} shows the results where the mesh size is kept to $h=0.025$ m while $\ell$ alters. We see that a smaller value of $\ell$ corresponds to a better match with the analytical solution.

\subsubsection*{Comparison against previously proposed models}
Next, we compare the three previously proposed models ($W_0$, $W_1$, and $W_2$) with our proposed model ($W$) using the same mesh resolution ($a_0/h=44$), length scale parameter ($\ell/h=2$) and material parameters. These four models introduce different modifications to the poroelastic strain energy corresponding to different momentum balance equations. 

\begin{figure}[htpb!]
    \centering
    \subfigure[With localized pressure field and $\alpha_\mrm{m}=0$]{\includegraphics[scale=0.5]{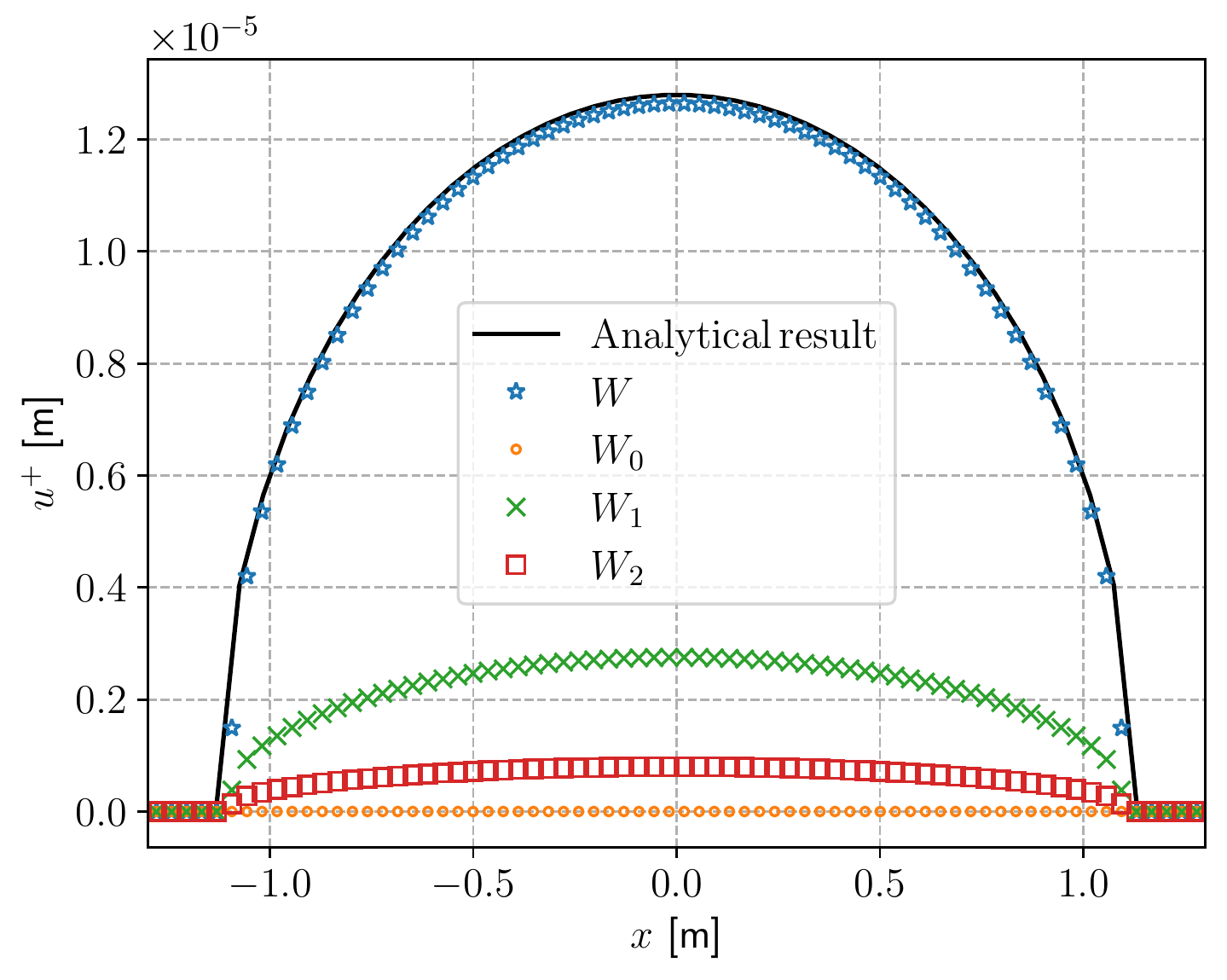}\label{fig:COD_model_localized}}
    \subfigure[With distributed pressure field and $\alpha_\mrm{m}=0$]{\includegraphics[scale=0.5]{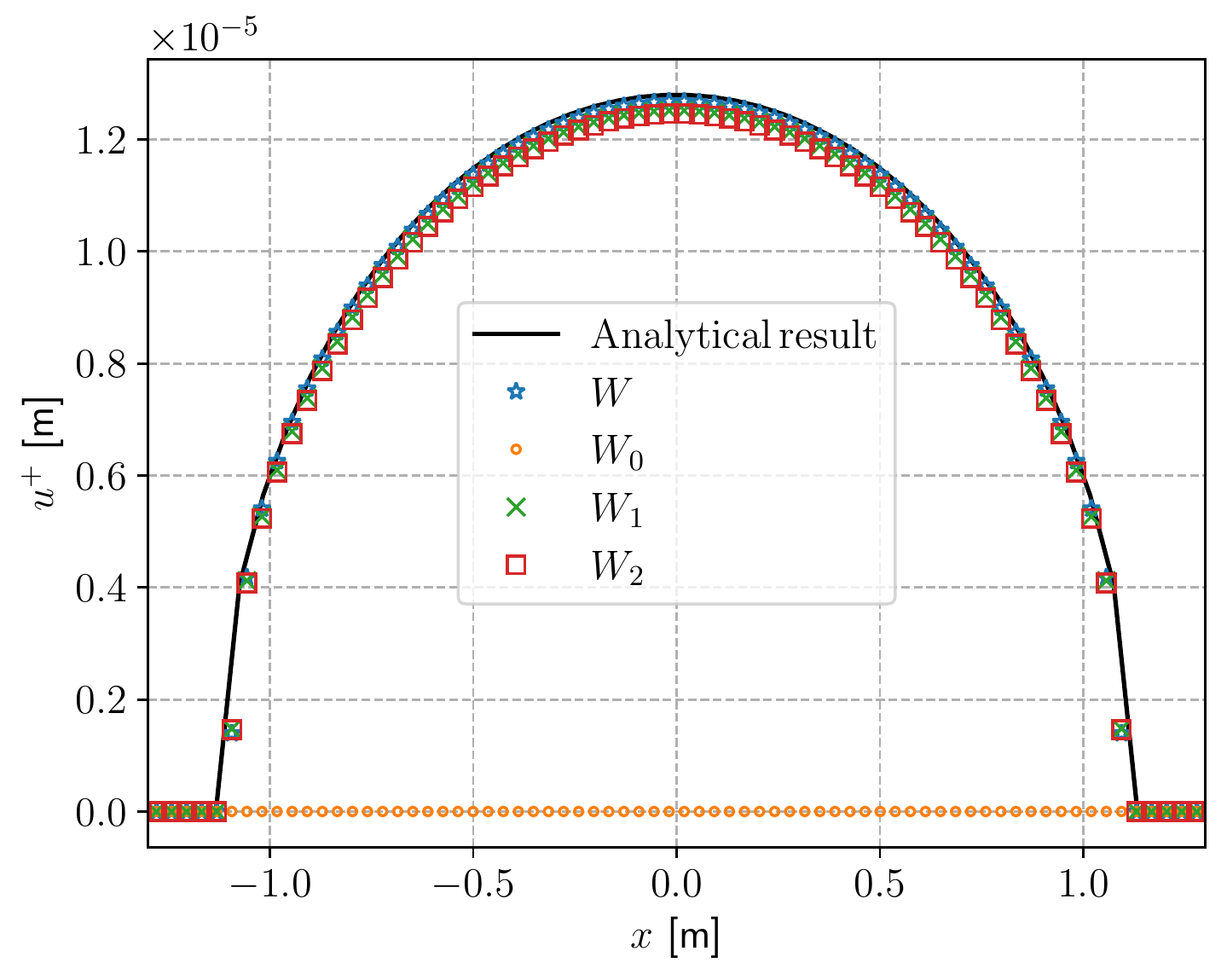}\label{fig:COD_model_distributed}}
    \caption{Crack opening displacements for various phase-field models ($W$, $W_0$, $W_1$, and $W_2$) with (a) localized and (b) distributed pressure fields.}  
    \label{fig:COD_model}
\end{figure}

Fig.~\ref{fig:COD_model_localized} compares crack opening displacement computed using the localized pressure field.
Our proposed model with diffused poroelasticity matches the analytical result, while the other three models fall short of it.
However, we should consider the two underlying assumptions in these models.

The first assumption is about the pressure field.
For $W_1$ and $W_2$, the deformation depends on $\nabla v$ (or $\nabla p$ after integration by parts in Eqs.~\eqref{eq:W1_reg_final}~and~\eqref{eq:W2_reg_final}).
In other words, without diffusing the fluid pressure in the crack, the deformation in the surrounding media will be underestimated. 
Thus, for further comparisons, we consider the distributed pressure field where the fluid pressure ($p=10^5$ Pa) is uniformly distributed in the entire domain (Fig.~\ref{fig:numerical_case3}c).
Fig.~\ref{fig:COD_model_distributed} shows computed cracking openings with the distributed pressure field and we can see that the models with $W$, $W_1$, and $W_2$ match the analytical solution, but not with $W_0$.

\begin{figure}[htpb!]
    \centering
    \subfigure[With localized pressure field]{\includegraphics[scale=0.5]{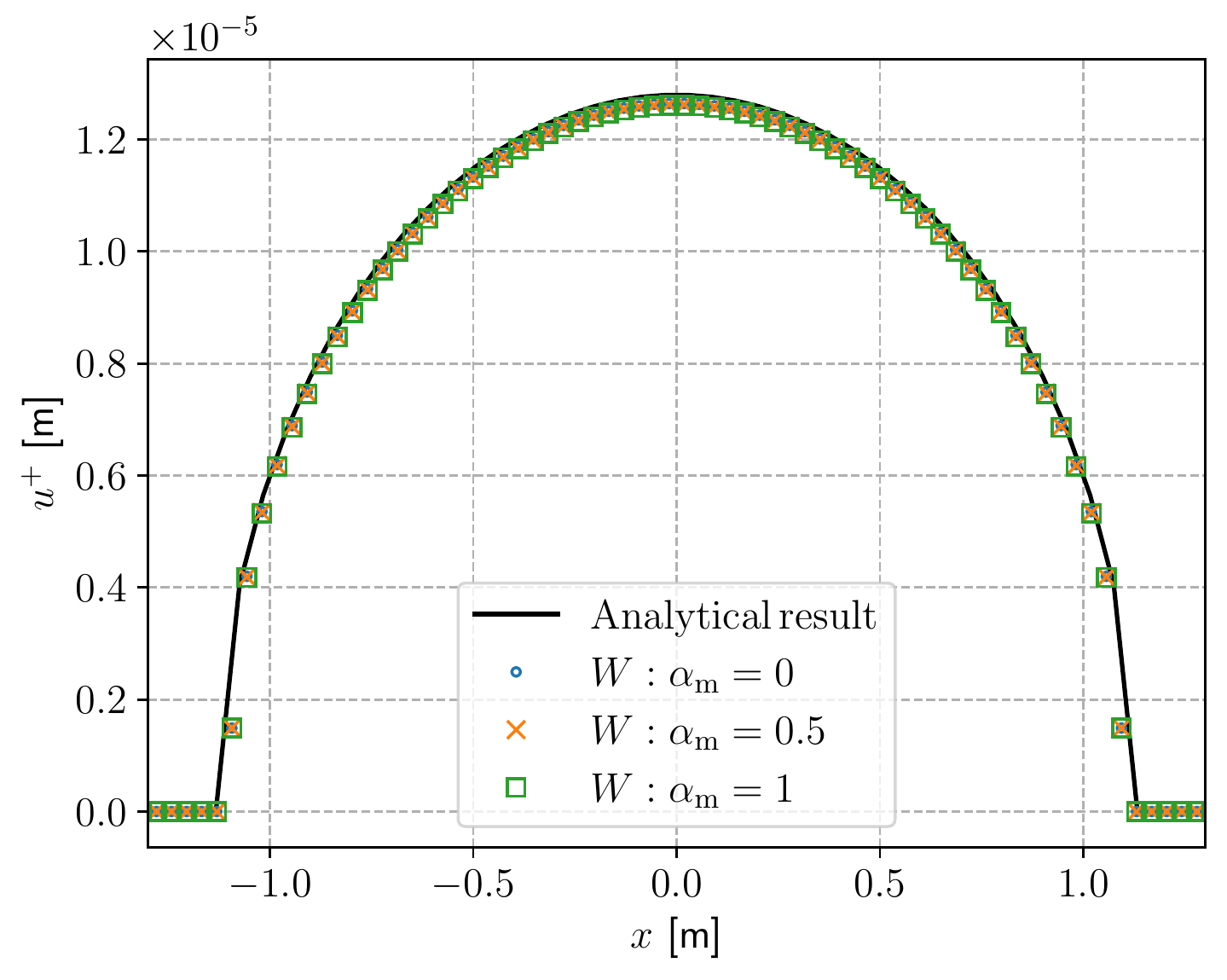}}
    \subfigure[With localized pressure field]{\includegraphics[scale=0.5]{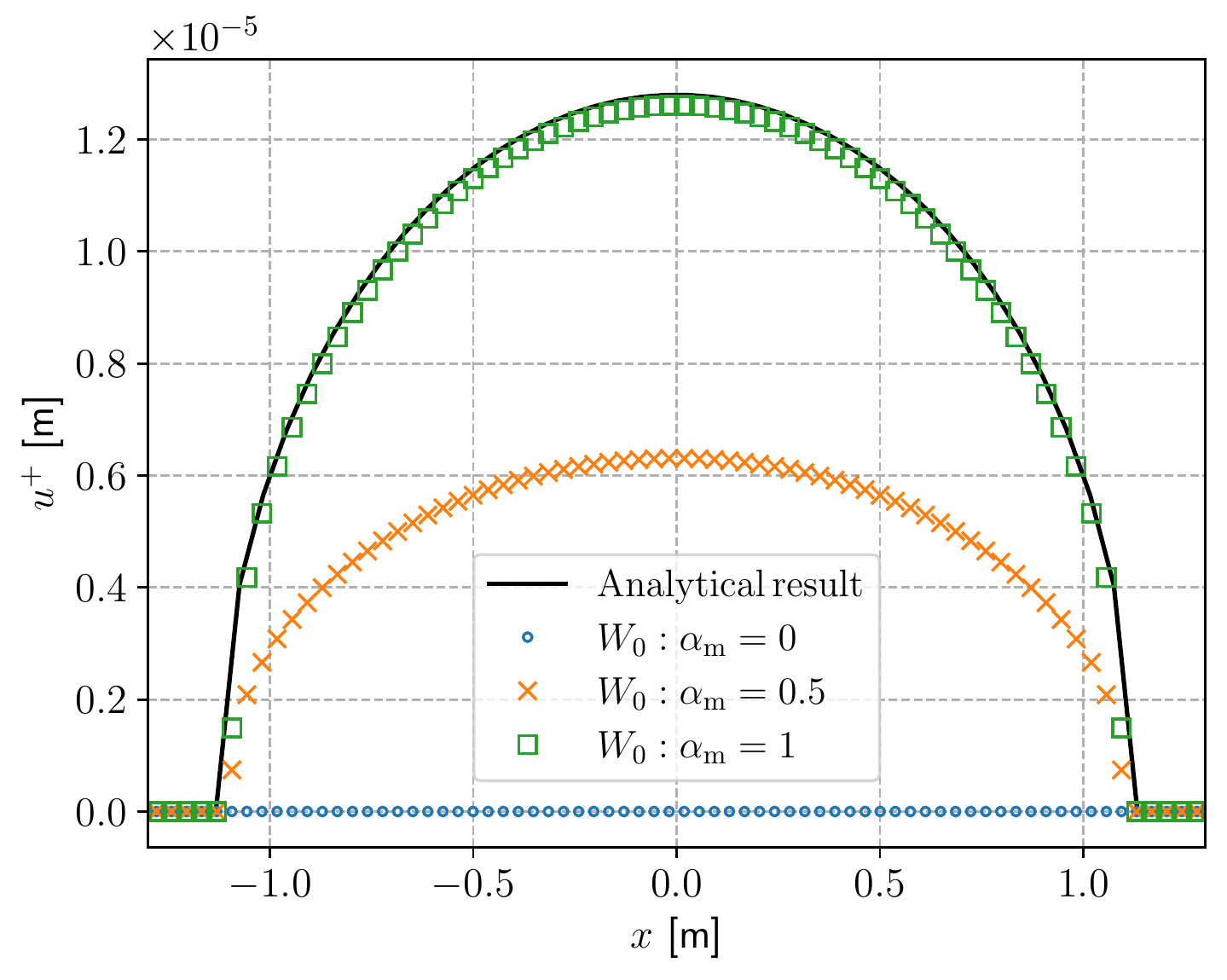}}
    \subfigure[With localized pressure field]{\includegraphics[scale=0.5]{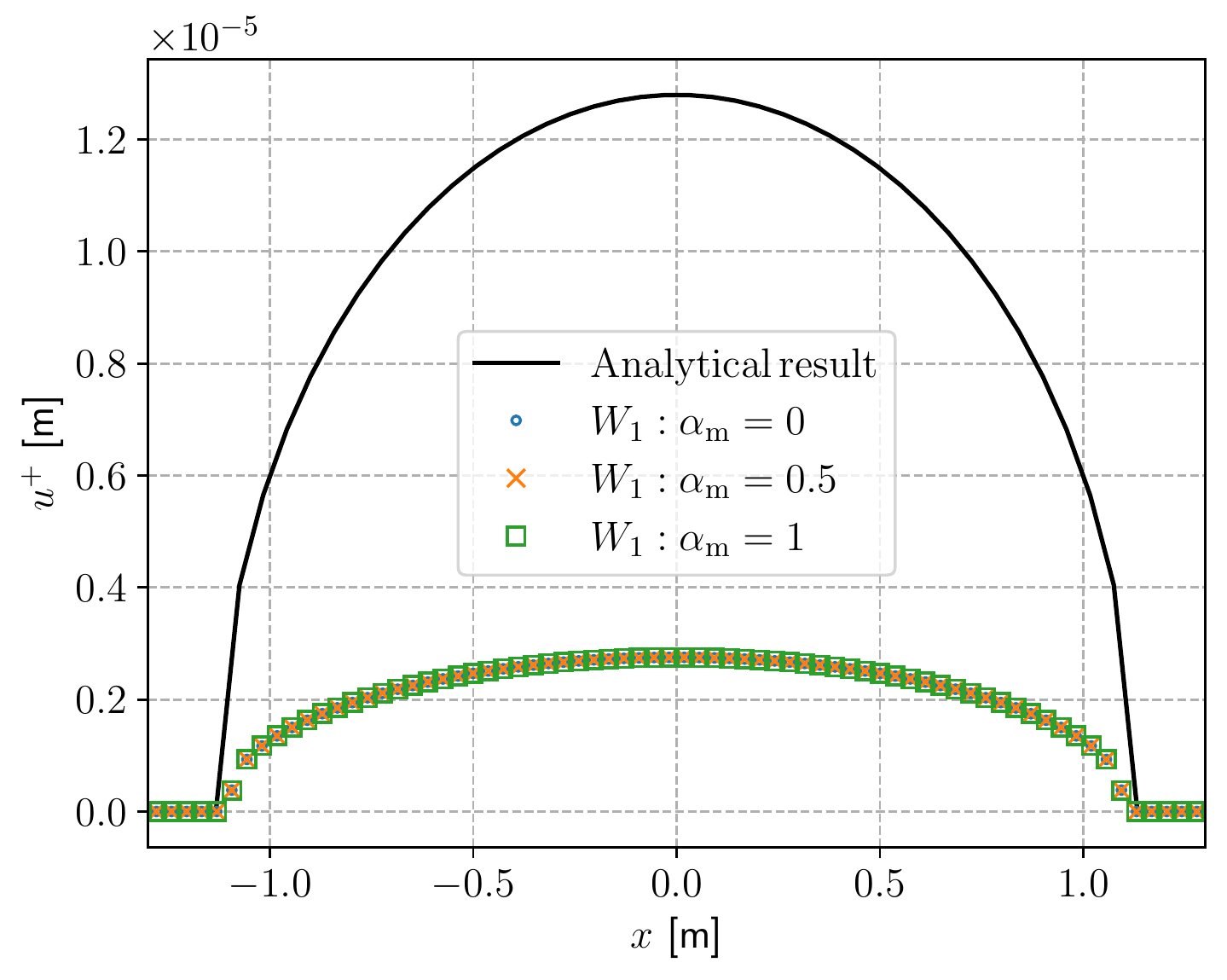}}
    \subfigure[With localized pressure field]{\includegraphics[scale=0.5]{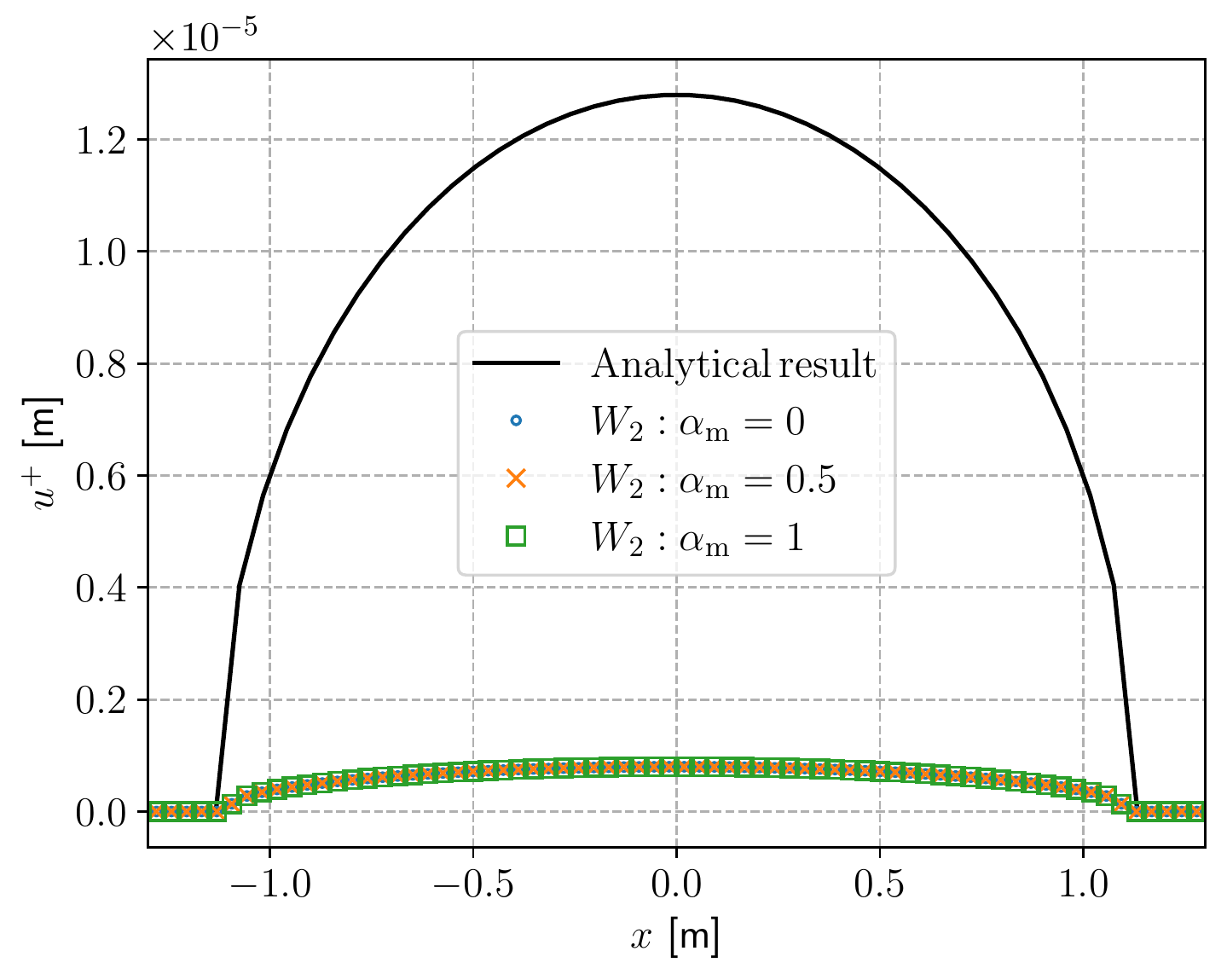}}
    \caption{The effect of Biot's coefficient $\alpha_\mrm{m}$ on the crack opening for various models  ($W$: The proposed model Eq.~\eqref{eq:poroelasticityVariational};  $W_0$: The model Eq.~\eqref{eq:W0_reg_final} used in Zhou et al \citep{zhou2019phase} and Heider and Sun \citep{heider2020phase}; $W_1$: The model Eq.~\eqref{eq:W2_reg_final} proposed by Chukwudozie et al \citep{Chukwudozie2019}; $W_2$: The model Eq.~\eqref{eq:W1_reg_final} proposed by Mikeli\'{c} et al \citep{Mikelic2015_Multi}).}
    \label{fig:COD_alpha}
\end{figure}

The second assumption is about Biot's coefficient.
In particular, the model with $W_0$ requires that $\alpha_\mrm{m} = 1$ to guarantee the stress continuity condition on the crack interface (Eq.~\eqref{eq:momentum}--2).
On the other hand, the stress continuity on the crack interface is explicitly considered in the $W$, $W_1$, and $W_2$ models.
Figs.~\ref{fig:COD_alpha} show computed cracking opening displacements with three different Biot's coefficients ($\alpha_\mrm{m} = 0$, 0.5, and 1.0) with the localized pressure field.
The model with $W_0$ can now predict the analytical crack opening profile with $\alpha_\mrm{m} = 1$ accurately, but the results are sensitive to the value of $\alpha_\mrm{m}$ unlike the other models. 



To be fair, we should note that we have tested these models in the conditions for which they were not intended for.
The models with $W_1$ and $W_2$ diffuse a pressure field beyond the fully fractured domain ($v=0$) using the level--set~\cite{Lee2017_ls} or phase--field calculus~\cite{Chukwudozie2019}.
The setting of this distributed pressure may be unrealistic, but in actual computations, the pressure field only needs to be diffused over the transition zone where $v<1$, and as long as the transition zone is localized (i.e., small $\ell$), the pressure can also be localized to some extent. 

In this simple case, the pressure in the fracture domain is constant, and the pressure field can be easily diffused or reconstructed. 
In practice, however, the pressure profile in the fracture may vary, and that may cause inaccurate deformation for complex fracture networks. 
For this reason, the model should be able to estimate the crack opening correctly with the localized pressure field irrespective of Biot's coefficient such as $W$. 



\subsection{Fluid-driven fracture propagation}

\begin{figure}[htp!]
    \centering
    \includegraphics[scale=0.6]{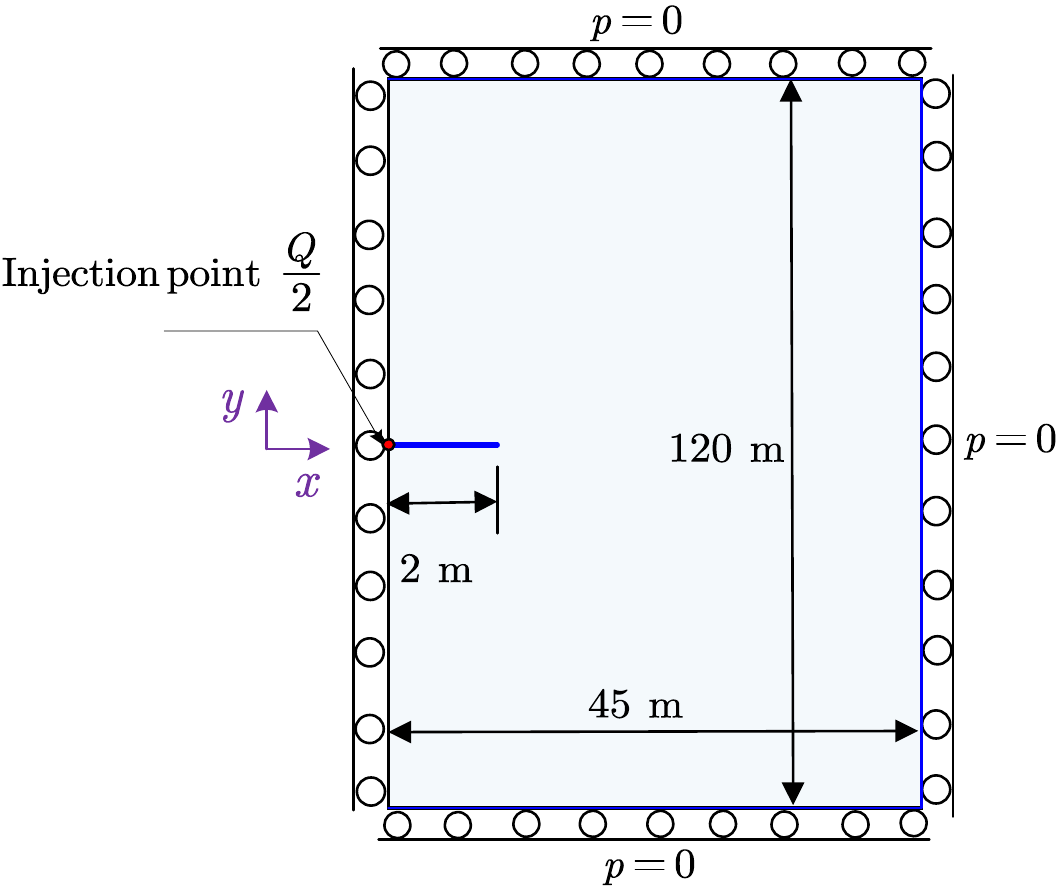}
    \caption{Geometry and boundary conditions for the KGD problem (without scaling).}
    \label{fig:KGD}
\end{figure}

We now move on to the hydraulic fracturing problem under a plane strain condition, known as the KGD (Kristianovich-Geertsma-de Klerk) model \cite{detournay2003near,hu2010plane}. The geometry and boundary conditions are presented in Fig.~\ref{fig:KGD} where only half of the model is considered because of the symmetry. The size of the computation domain is 45 m $\times$ 120 m and the height is taken relatively large to approximate an infinite plane. 
The KGD problem concerns the propagation of a planar fracture by injecting incompressible fluid into the impermeable elastic domain. Therefore, we take $\alpha_\mrm{m}=0$, $\phi_\mrm{m}=0$ and $c_f=0$, mimicking the elastic solid and incompressible fluid respectively except for the permeability which has to be greater than zero for the stiffness matrix in Eq.~\eqref{eq: pressure discrete form} to be positive definite. The  weighting exponent $\xi$ for fracture permeability in Eq.~\eqref{eq:Kf} is taken as $\xi =10$.
The other parameters are listed in Table \ref{tab:case4 parameter}. 

\begin{table}[htp!]
    \centering
    \begin{tabular}{ccc}
    \hline
      Parameter   &  Value & Unit\\
      \hline
       Young's modulus ($E$)  & $17\times10^9$ & Pa \\
       Poisson's ratio ($\nu$)  & 0.2  & -- \\
       Fluid compressibility ($c_f$) &  $0$ & -- \\
       Fluid viscosity ($\mu$) &   $10^{-8}$ & Pa$\cdot$s \\
       Injection rate (Q) & $2\times 10 ^{-3}$ & m$^2$/s \\
        Critical surface energy release rate ($\Gc$)  & 300 & N/m \\
       \hline
    \end{tabular}
    \caption{Parameters for the hydraulic fracturing problem}
    \label{tab:case4 parameter}
\end{table}
Similarly to the analyses in Section \ref{sec:sneddon}, we use an effective fracture energy parameter $G^\mrm{eff}_\mrm{c}$ that accounts for the diffused representation of crack by the phase-field model~\cite{Bourdin2008, Tanne2018, Yoshioka2020}:
\begin{equation}
    G_\mrm{c}^\mrm{eff} = G_\mrm{c} \left( 1 + \frac{h}{4c_n\ell} \right) 
    .
\end{equation}

Also, to recover the crack length from a diffused-represented crack by the phase-field model, we employ this relationship~\cite{yoshioka2021variational}:
\begin{equation}
    L = \frac{\int_\Omega \frac{\Gc}{4c_n} \left[\frac{(1-\nu)^n}{\ell}+ \ell \nabla v\cdot \nabla v \right]\,\mathd V}{G^\mrm{eff}_\mrm{c}}
    ,
    \label{eq:crack_length}
\end{equation}
where the numerator of the right hand side denotes the crack surface energy which can be readily computed using Gaussian integral on the discretized domain.

For the toughness-dominated regime, Garagash \cite{garagash2006plane} presents an analytical solution (see also \cite{santillan2017phase} Appendix A for a detailed derivation). 
The dimensionless viscosity $\mathcal{M}$ is defined as
\begin{equation}
    \mathcal{M} = \frac{\mu^\prime Q}{E^\prime}\left( \frac{E^\prime}{K^\prime}\right)^4
    ,
\end{equation}
where $\mu^\prime = 12 \mu$, $E^\prime = \frac{E}{1-\nu^2}$ and $K^\prime = \sqrt{\frac{32G^\mrm{eff}_\mrm{c} E^\prime}{\pi}}$. If $\mathcal{M}$ is below $\mathcal{M}_c=3.4\times 10^{-3}$, the propagation regime is toughness dominated. 
And we have $\mathcal{M}=3.8\times 10^{-7}$ with the parameters listed in Table~\ref{tab:case4 parameter}.

We now examine the proposed phase-field model $W$ (Eq.~\eqref{eq:poroelasticityVariational}) with various mesh sizes and matrix permeabilities. The ratio of $\ell$ and $h$ remains constant, i.e., $\ell / h =4$ throughout different mesh sizes. The time increment is taken as $\Delta t=0.01$ s in the first 10 steps and $\Delta t=0.1$ s in the remaining steps. The pressure and  crack opening at the injection point and the length of the crack are recorded during hydraulic fracturing. 

\begin{figure}[htp!]
    \centering
    \subfigure[The effect of mesh size ($K_\mrm{m}=1 \times 10^{-18}$ $\mrm{m}^2$)]{\includegraphics[scale=0.5]{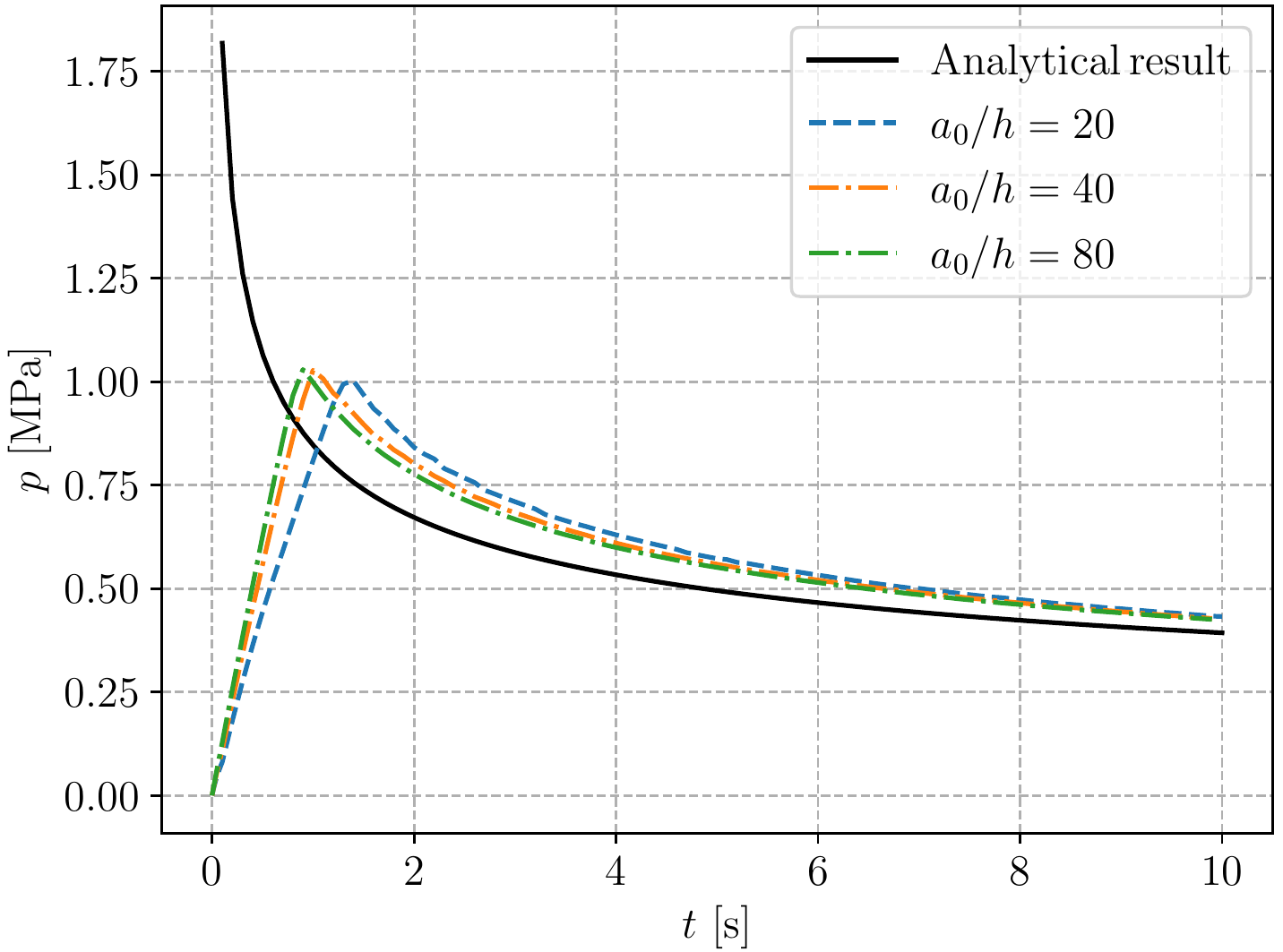} \label{fig:p_mesh_effect}}
    \subfigure[The effect of matrix permeability ($a_0/h=40$)]{\includegraphics[scale=0.5]{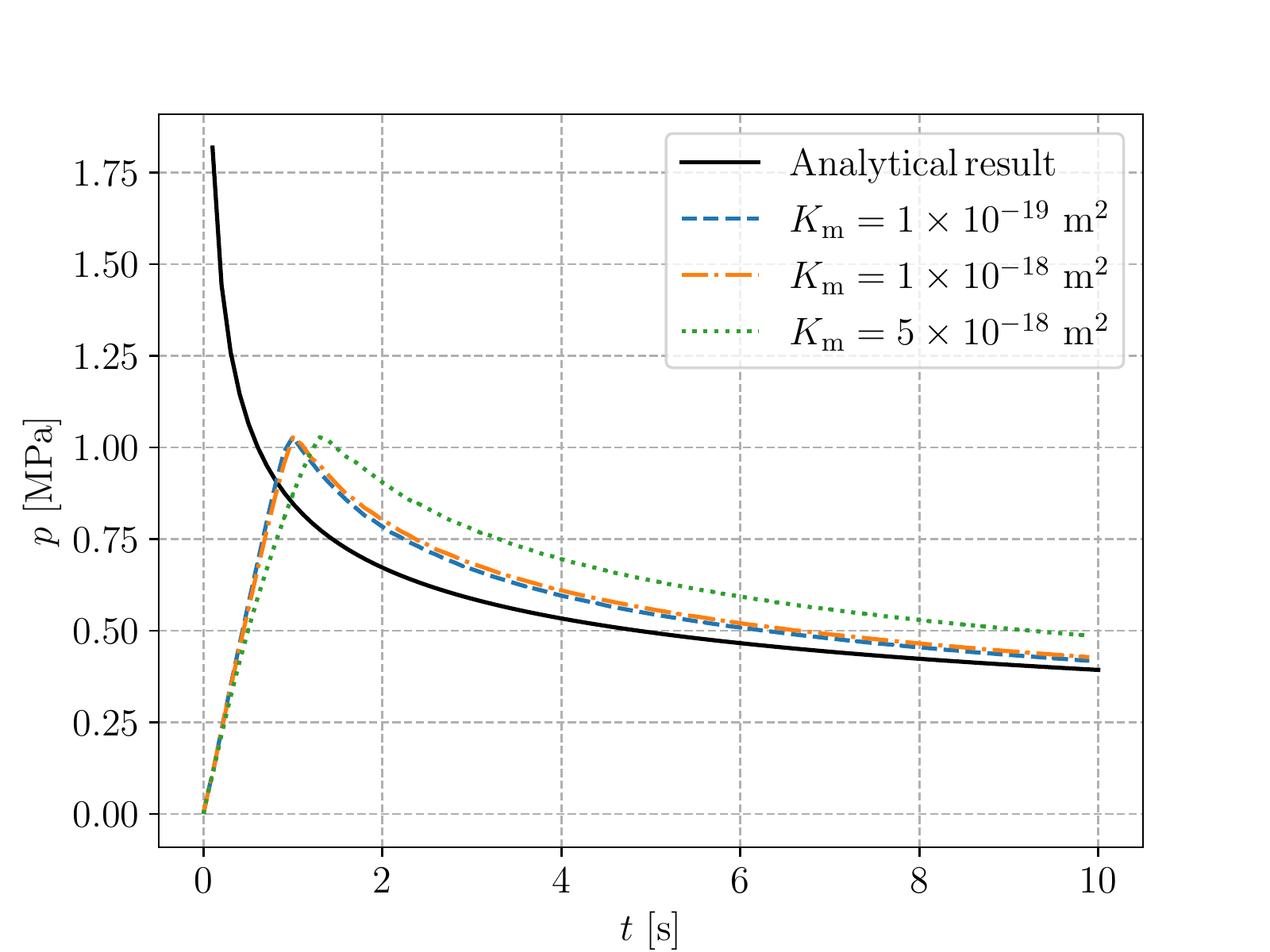} \label{fig:p_K_effct}}
    \caption{Pressures at the injection point for different meshes and permeabilities of the matrix. }
\end{figure}

\begin{figure}[htp!]
    \centering
    \subfigure[The effect of mesh size ($K_\mrm{m}=1 \times 10^{-18}$ $\mrm{m}^2$)]{\includegraphics[scale=0.5]{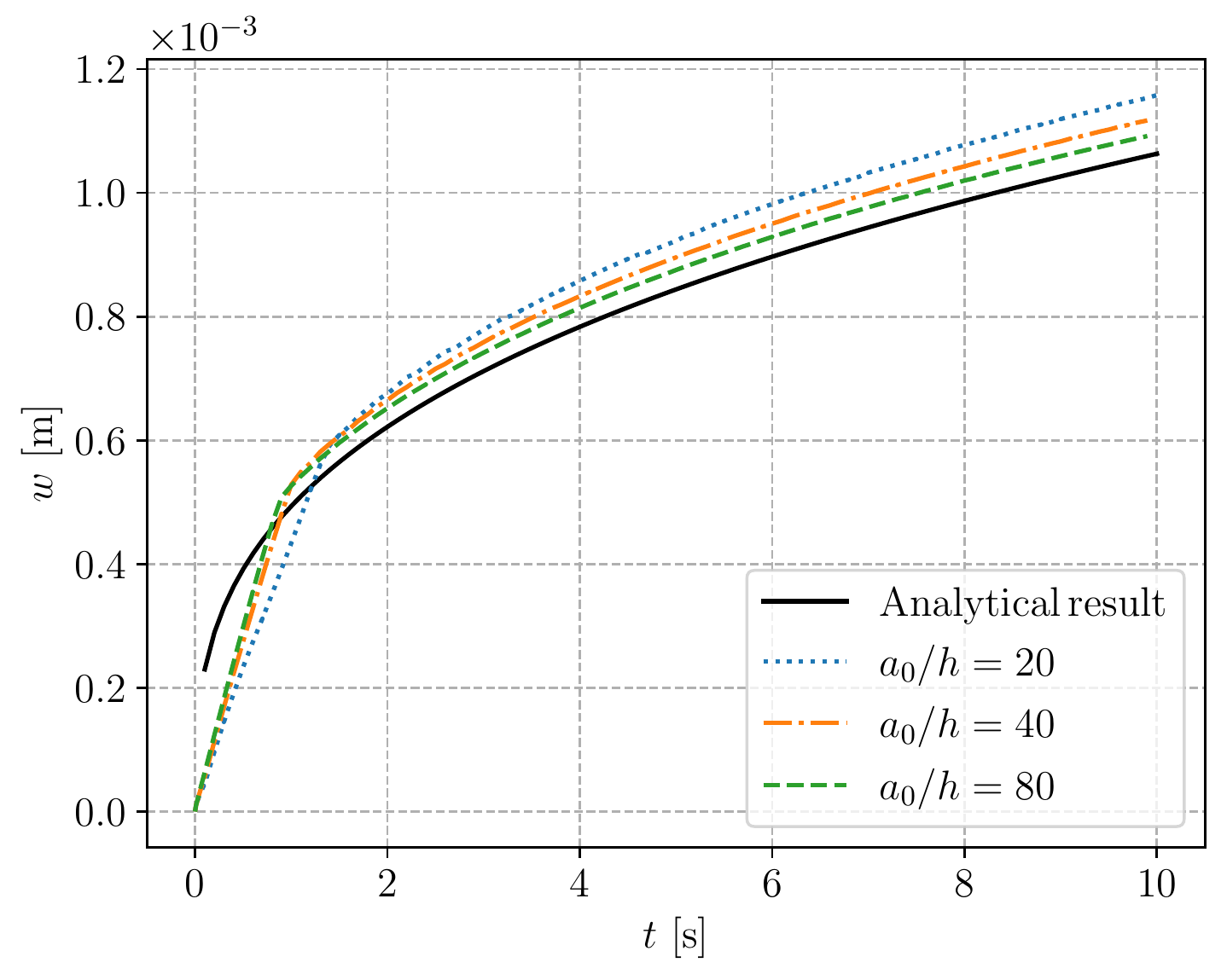} \label{fig:w_mesh_effect}}
    \subfigure[The effect of matrix permeability ($a_0/h=40$)]{\includegraphics[scale=0.5]{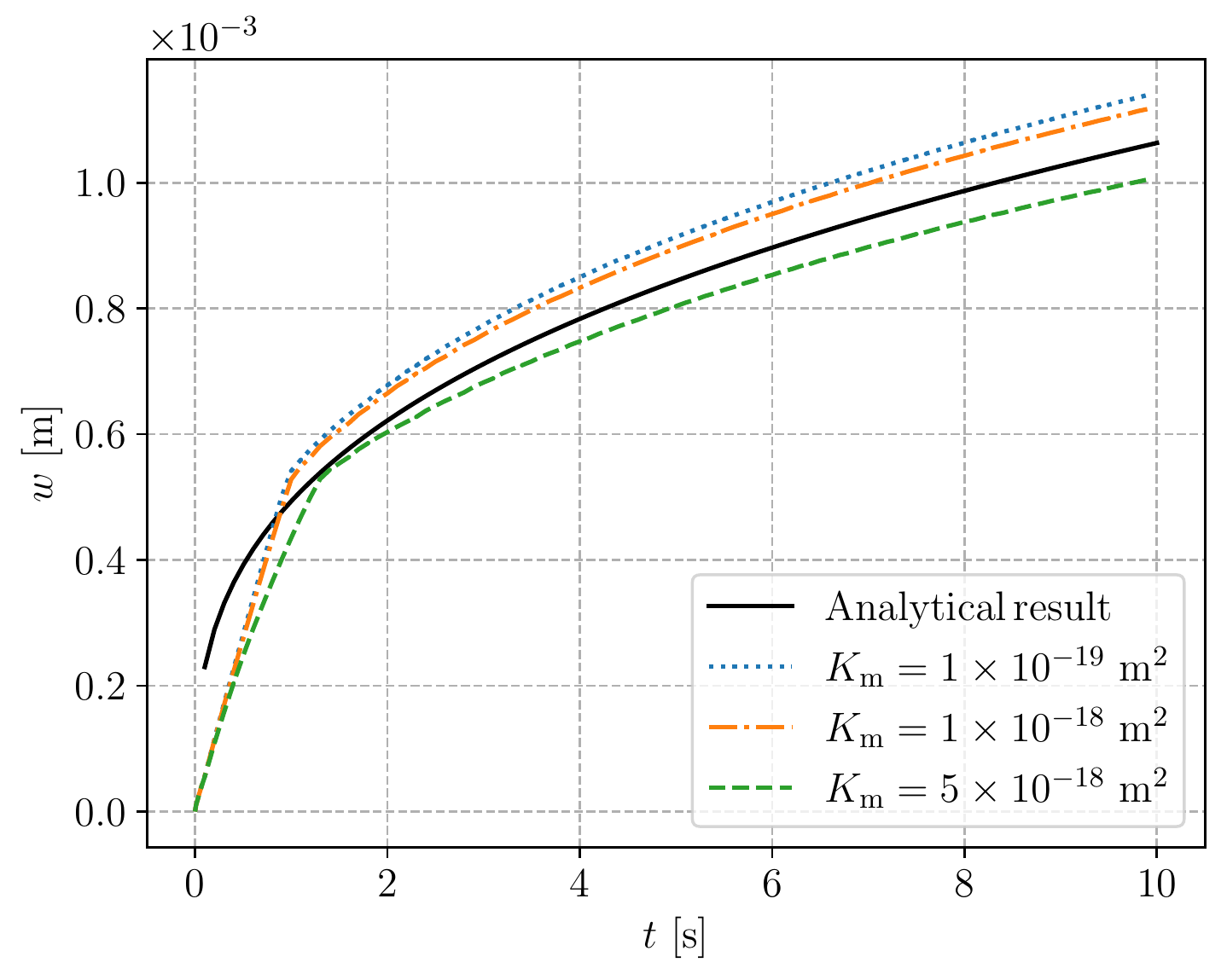} \label{fig:w_K_effect}}
    \caption{Crack opening at the injection point for different meshes and permeabilities of porous material. } 
\end{figure}

\begin{figure}[htp!]
    \centering
    \subfigure[The effect of mesh size ($K_\mrm{m}=1 \times 10^{-18}$ $\mrm{m}^2$)]{\includegraphics[scale=0.5]{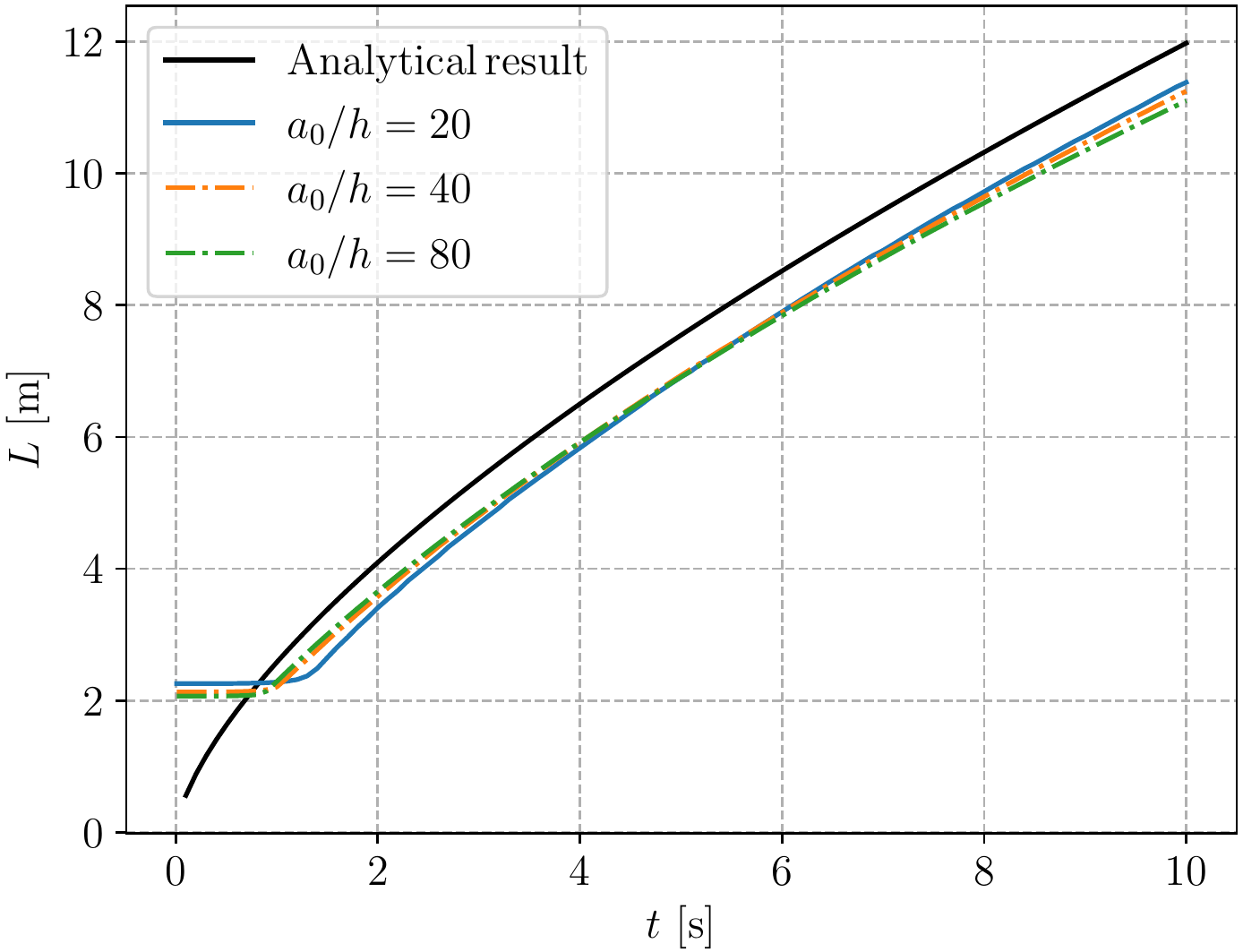} \label{fig:Length_mesh_effect}}
    \subfigure[The effect of matrix permeability ($a_0/h=40$)]{\includegraphics[scale=0.5]{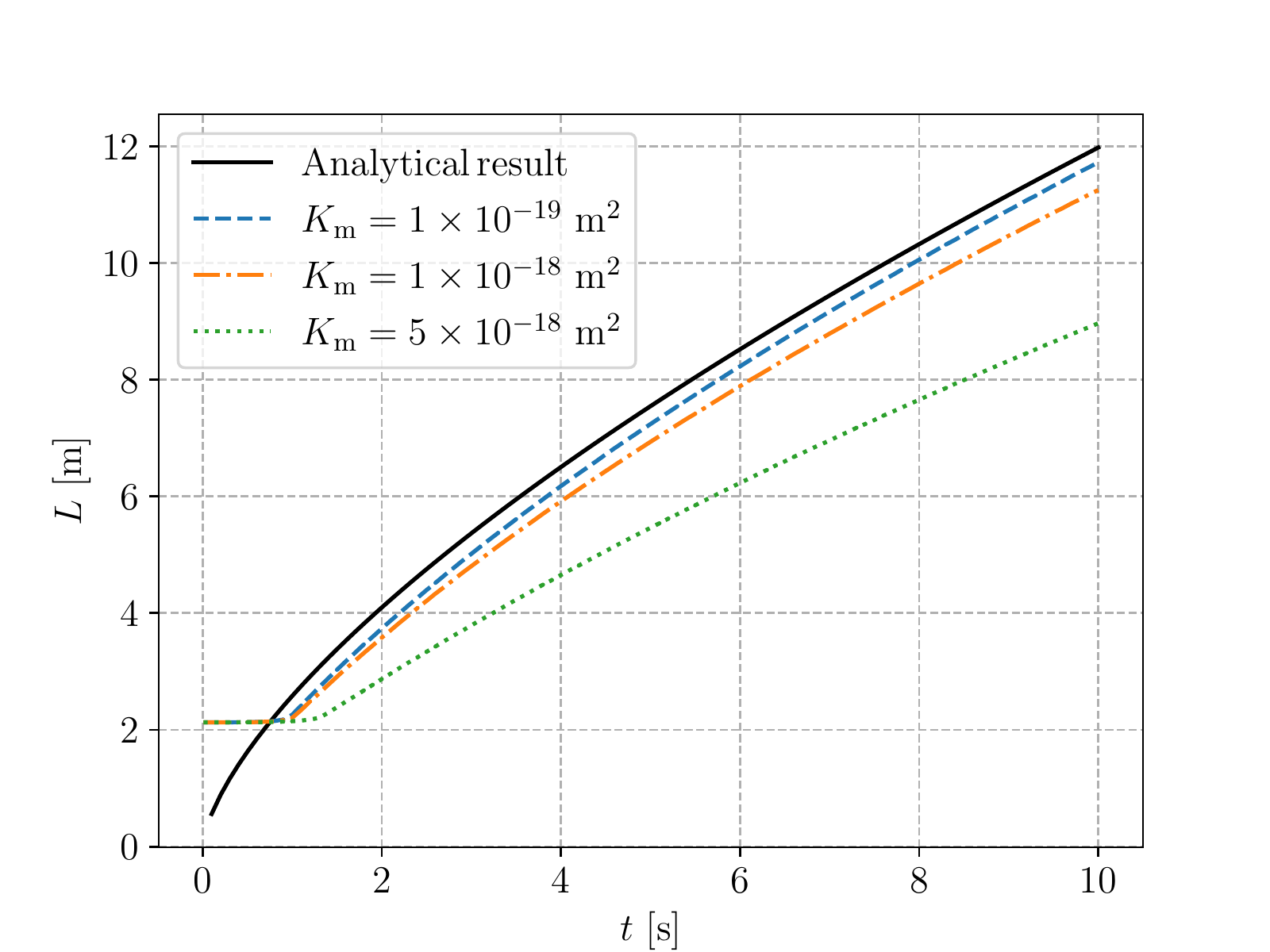} \label{fig:Length_K_effect}}
    \caption{Crack length during the injection for different meshes and permeabilities of porous material. }
\end{figure}

Using $K_\mrm{m}=1\times 10^{-18}$ m$^2$, we compare the results for three different mesh sizes, i.e., $a_0/h=20, 40 \mrm{\, and \,} 80$ in Figs. \ref{fig:p_mesh_effect}, \ref{fig:w_mesh_effect} and \ref{fig:Length_mesh_effect}. 
As we can see, a finer mesh corresponds to a closer approximation to the analytical solutions of the pressure and crack opening at the injection point, but the crack length is less sensitive to the mesh refinement because the expression in Eq.~\eqref{eq:crack_length} already accounts for the effects of phase-field regularization.

Using $a_0/h=40$, we compare the results for three different values of the permeability, i.e., $K_\mrm{m} = 1\times 10^{-19}, 1 \times 10^{-18} \mrm{\, and \,} 5\times 10^{-18}$. Fig.~\ref{fig:p_K_effct} shows that the critical pressure for the onset of fracture propagation is not affected by the permeability of the porous material. However, a higher permeability will delay the onset of fracture propagation because more fluid leaks off to the surrounding porous media. Furthermore, the crack opening and propagation are smaller for the higher permeabilities (Fig. \ref{fig:w_K_effect} and Fig.~\ref{fig:Length_K_effect}). 
To capture the KGD fracture propagation in an impermeable media, the permeability needs to be low (e.g., $K_\mrm{m} = 1 \times 10^{-19}$ m$^2$), but extremely low permeabilities will cause numerical instabilities in solving the hydraulic equation. 

\section{Hydraulic fracture interacted with a natural fracture}
\label{sec:hydraulic_interface}
\begin{figure}[htpb!]
    \centering
    \includegraphics[scale=0.8]{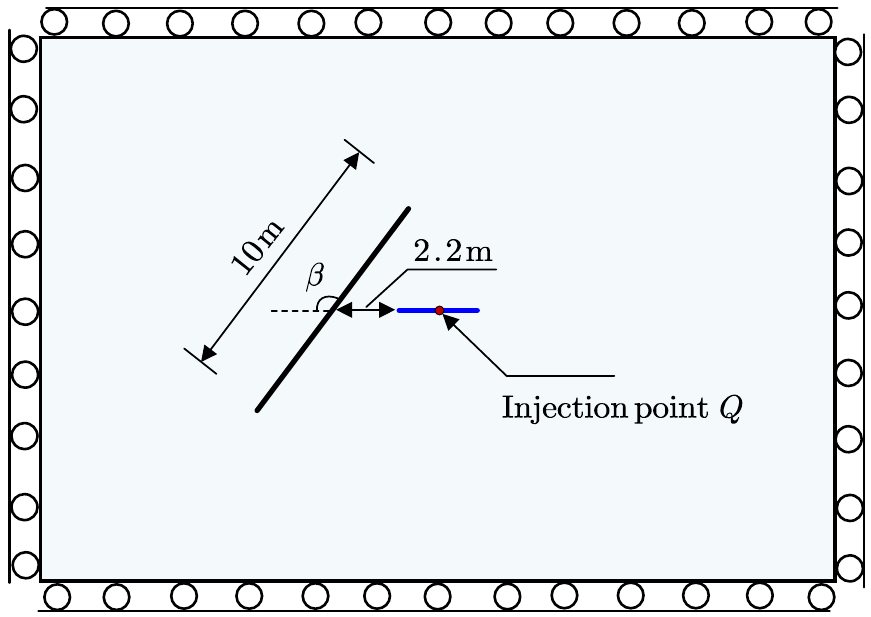} 
    \caption{Geometry and boundary conditions for the hydraulic fracture interactions with a natural fracture example (without scaling).}
    \label{fig:interaction_hydraulic}
\end{figure}

In this section, we demonstrate the proposed model's ability to simulate hydraulic fracture propagation under the presence of a natural fracture. 
We represent a natural fracture using the interface modeling method proposed by~\cite{Hansen-Dorr2019,yoshioka2021variational} by assigning the effective interface fracture toughness.
This method provides a simple but robust description of the deflection and penetration of a crack impinging into an interface. The geometry of this numerical model is shown in Fig.~\ref{fig:interaction_hydraulic} where a natural fracture with weaker fracture toughness is placed 2.2 m away from the initial hydraulic fracture. The sizes of the domain and initial hydraulic fracture are the same as the setting in Fig. \ref{fig:numerical_case3}, and the length of this natural fracture is 10 m. We consider two different inclination angles, $\beta=90^\circ \mrm{\, and \,} 165^\circ$, and two different interface toughness, $\Gc^{\mrm{int}}/\Gc$=0.2 and 0.5. The parameters are the same as those listed in Table \ref{tab:case4 parameter} except that we take Biot's coefficient $\alpha_\mrm{m} =0.6$, porosity $\phi_\mrm{m}=0.01$ and permeability $K_m=5\times10^{-18}$ to account for poroelastic effects. The  weighting exponent $\xi$ for fracture permeability  is taken as $\xi =50$ in this case.

An incompressible fluid is injected at the center of the initial hydraulic fracture with the injection rate of $2\times 10 ^{-3}$ m$^2$/s. The crack will propagate after the pressure reaches a critical value and then impinge on the natural fracture. 

\begin{figure}[htpb!]
    \centering
    \subfigure[Phase field]{\includegraphics[width=7cm]{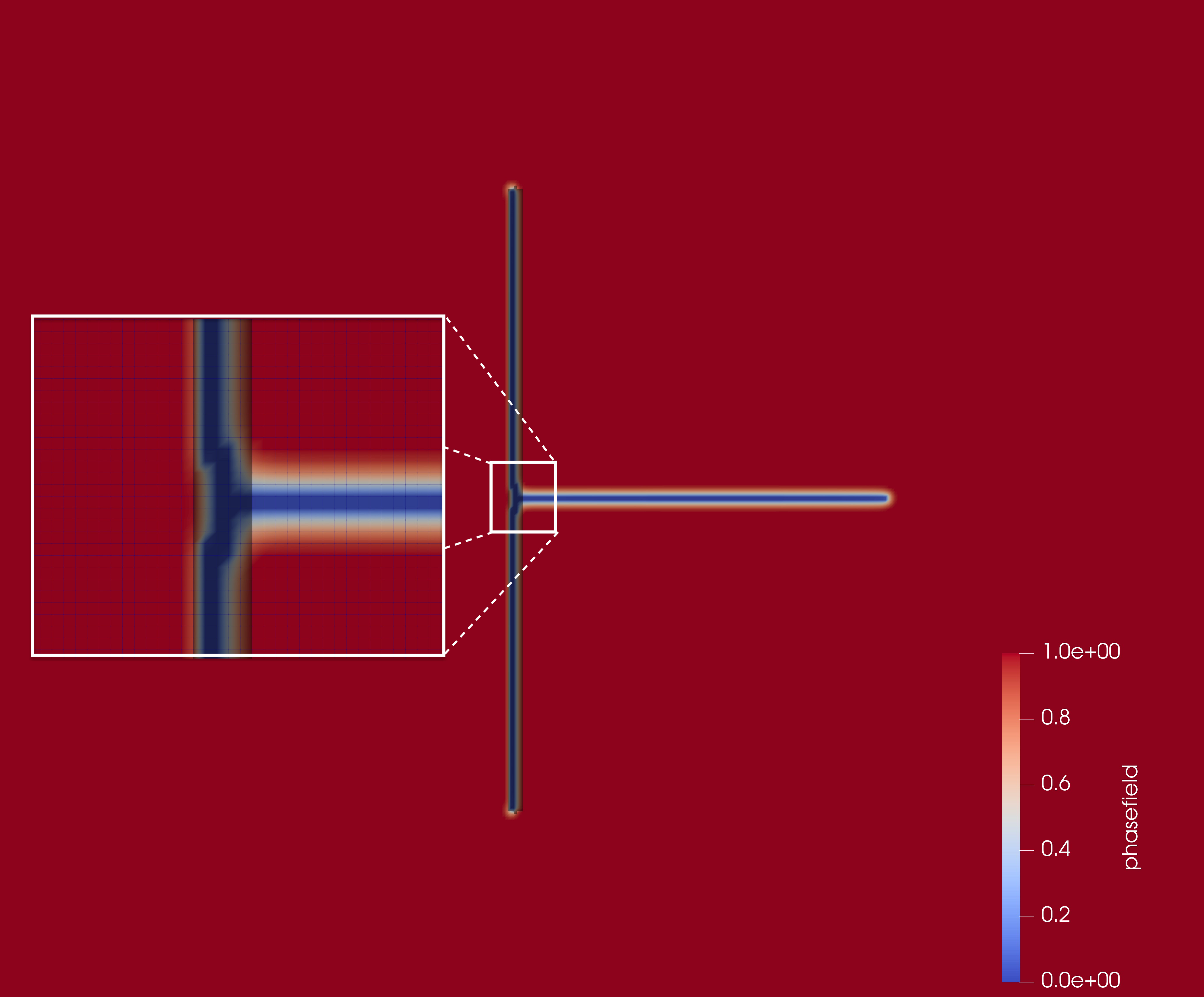} \label{fig:pf_combined_int_1}}
    \subfigure[Pressure]{\includegraphics[width=7cm]{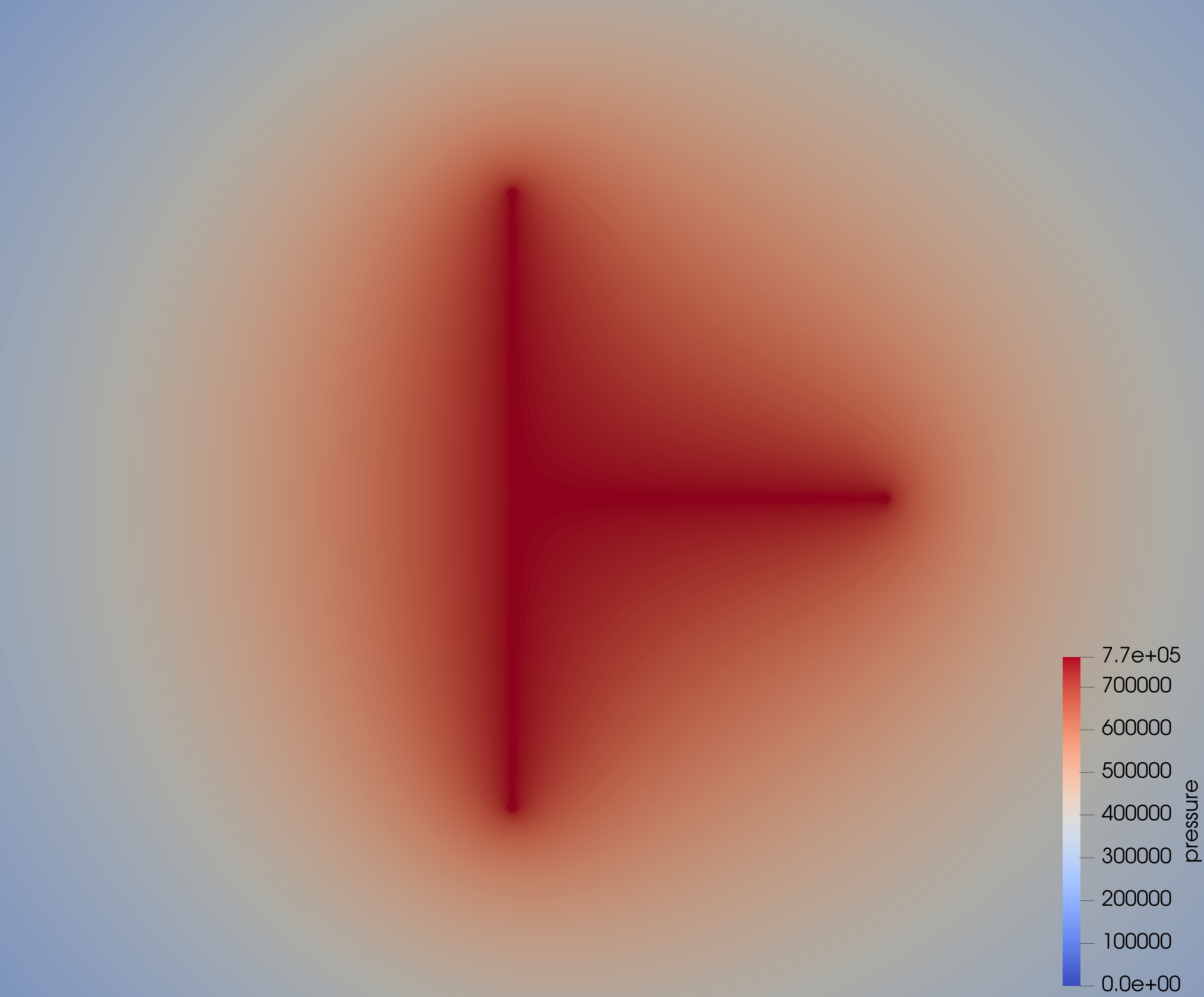} \label{fig:pressure_int_1}}
    \subfigure[Width (Interpolated onto node)]{\includegraphics[width=7cm]{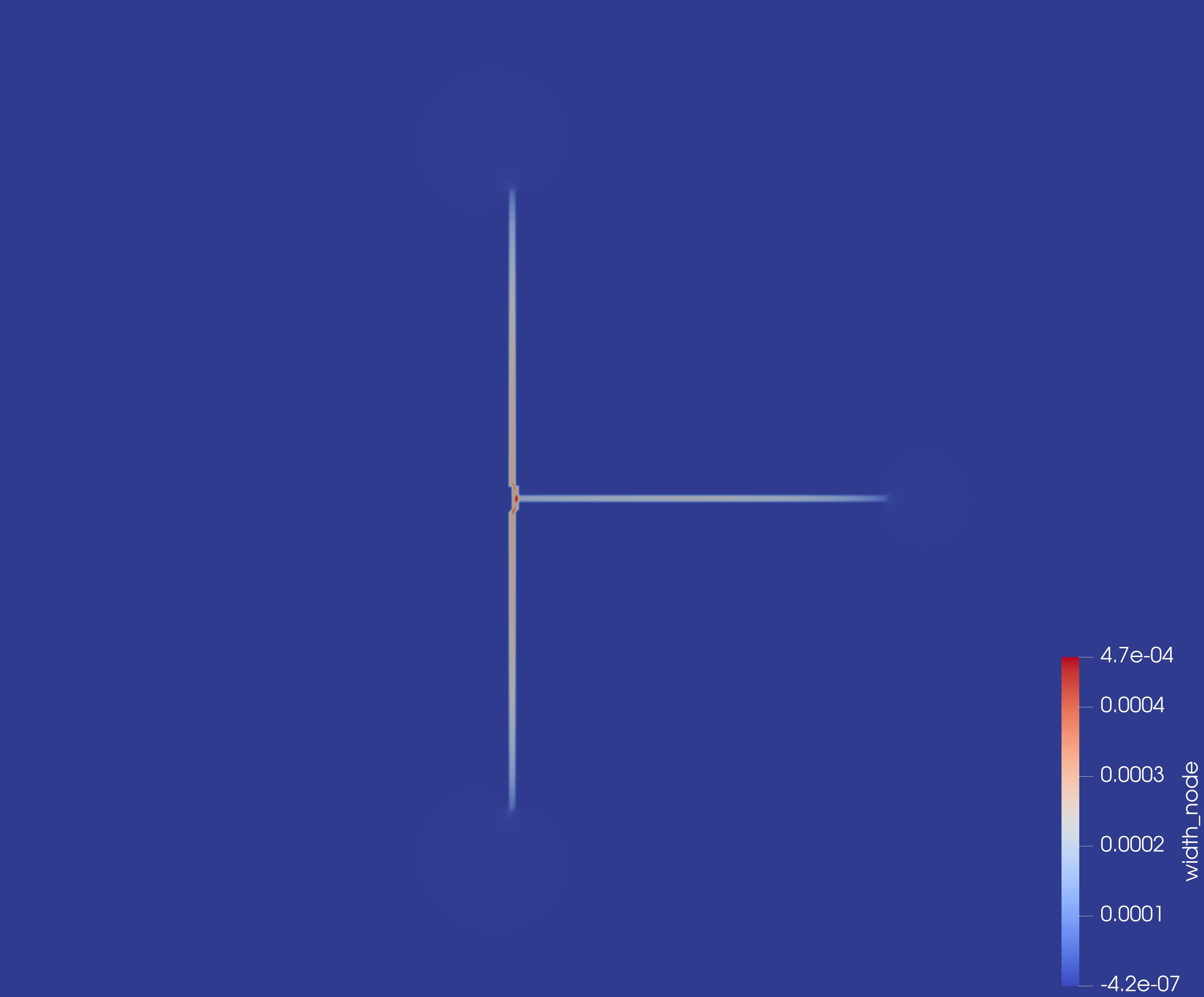} \label{fig:width_int_1}}
    \subfigure[Boit's coefficient (Calculated by Eq. \eqref{eq:alpha} )]{\includegraphics[width=7cm]{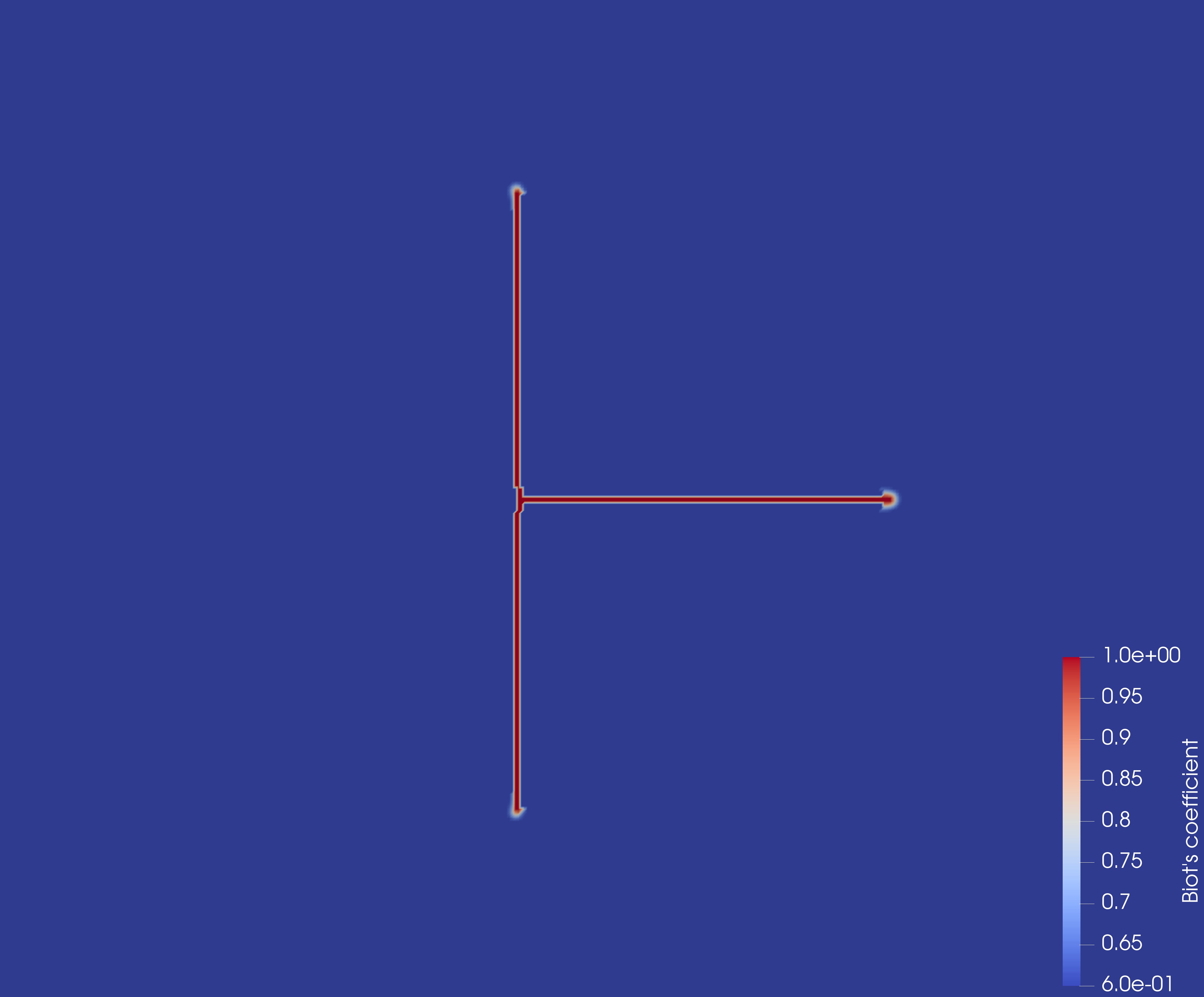} \label{fig:alpha_int_1}}
    \caption{Phase-field (a), pressure (b), width (b) and Biot's coefficient (d) profiles for the hydraulic fracturing interacting with the natural fracture at $t=10$ s. ($\beta=90^\circ$, $\Gc^\mrm{int}/\Gc=0.2$)  }
    \label{fig:results_90_02}
\end{figure}

\begin{figure}[htpb!]
    \centering
    \subfigure[Phase field]{\includegraphics[width=7cm]{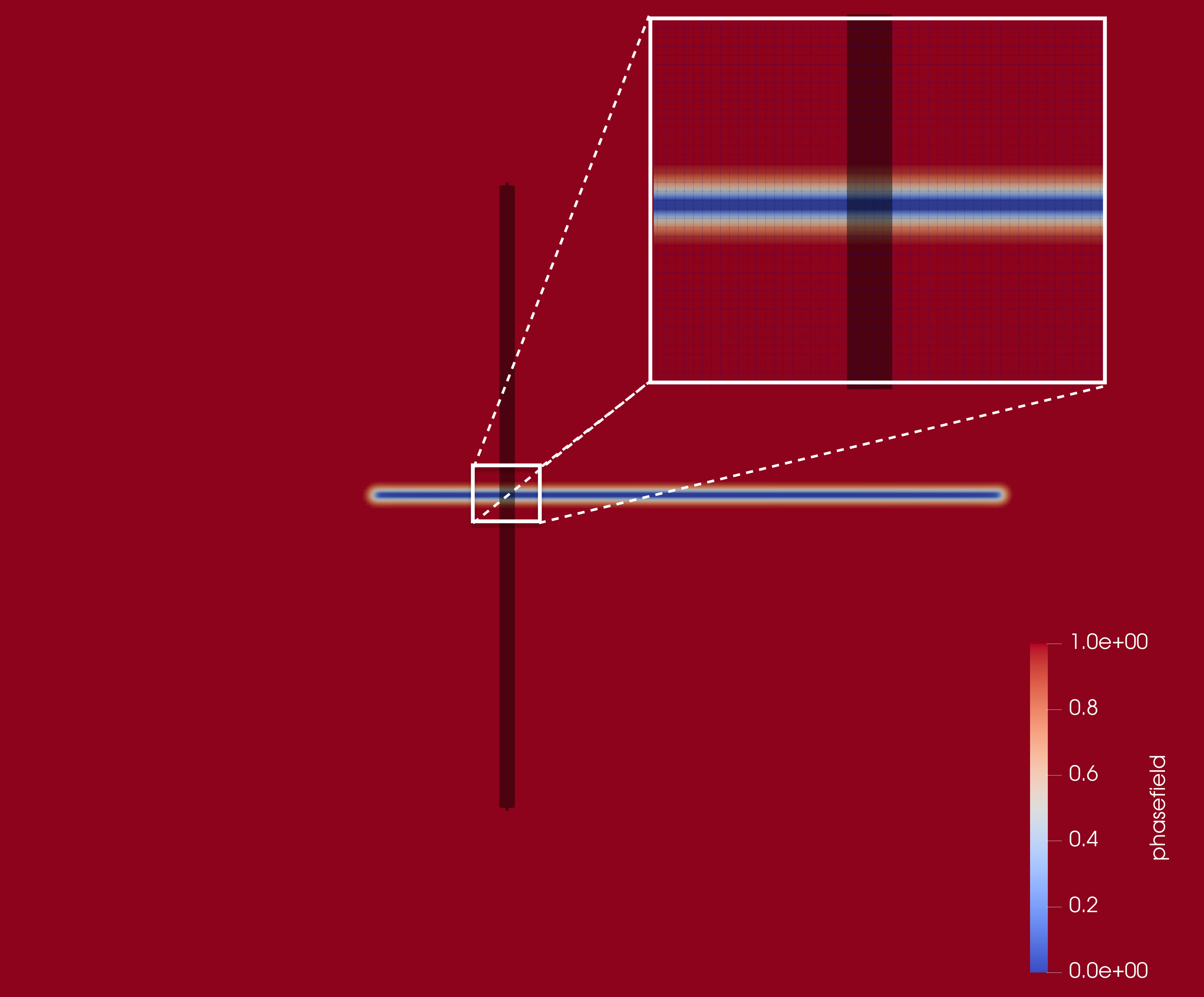} \label{fig:pf_combined_int_2}}
    \subfigure[Pressure]{\includegraphics[width=7cm]{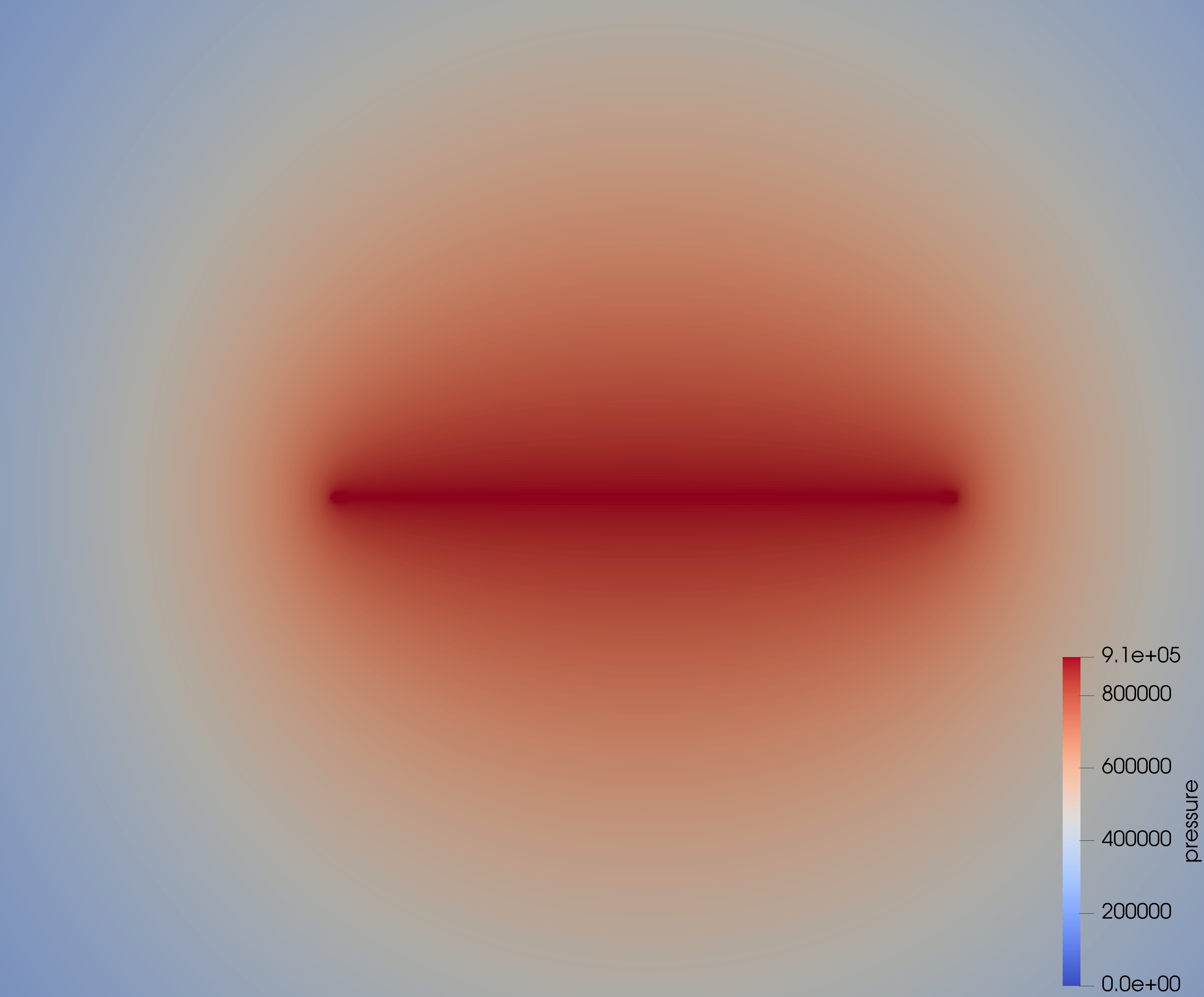} \label{fig:pressure_int_2}}
    \subfigure[Width (Interpolated onto node)]{\includegraphics[width=7cm]{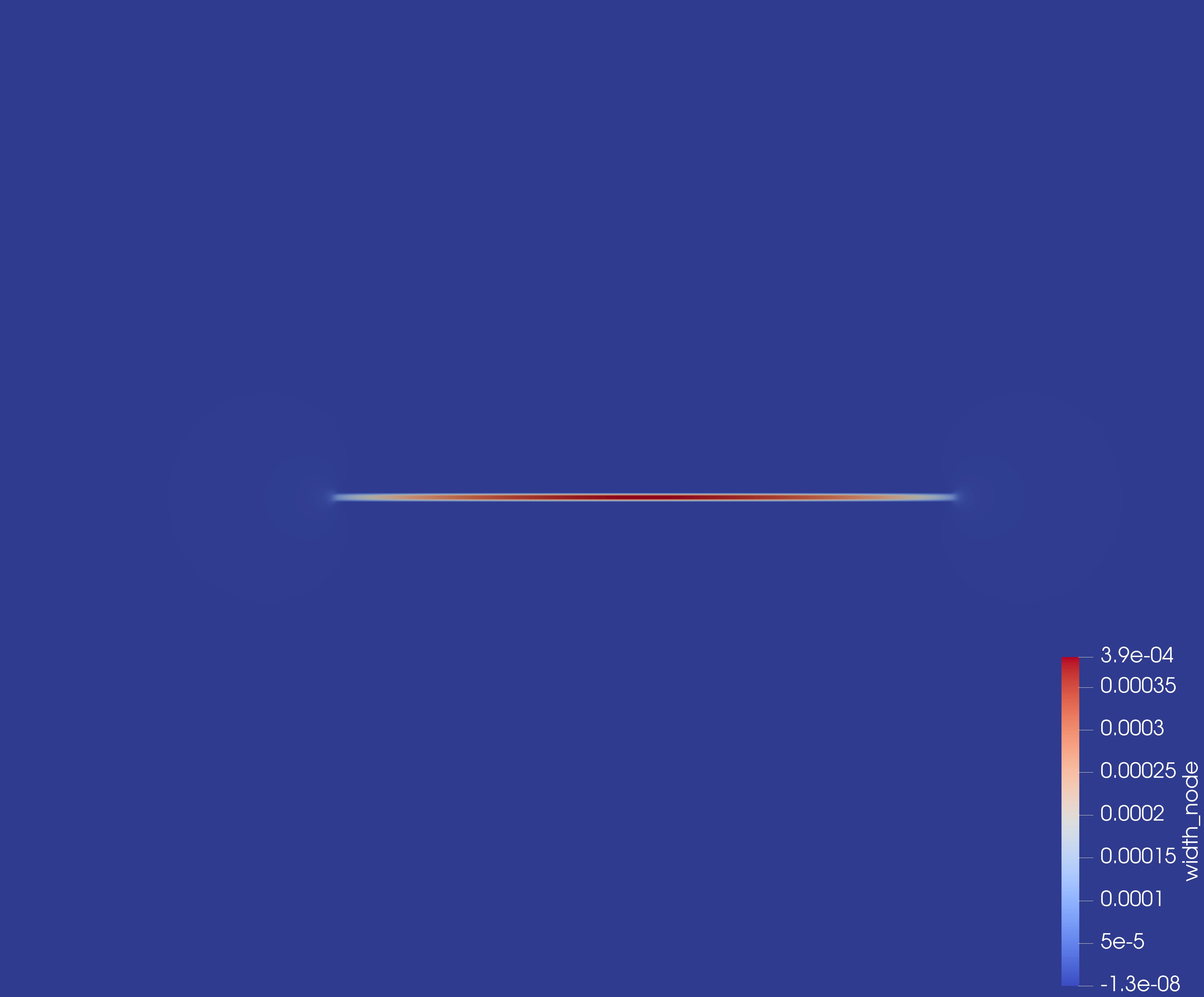} \label{fig:width_int_2}}
    \subfigure[Boit's coefficient (Calculated by Eq. \eqref{eq:alpha} )]{\includegraphics[width=7cm]{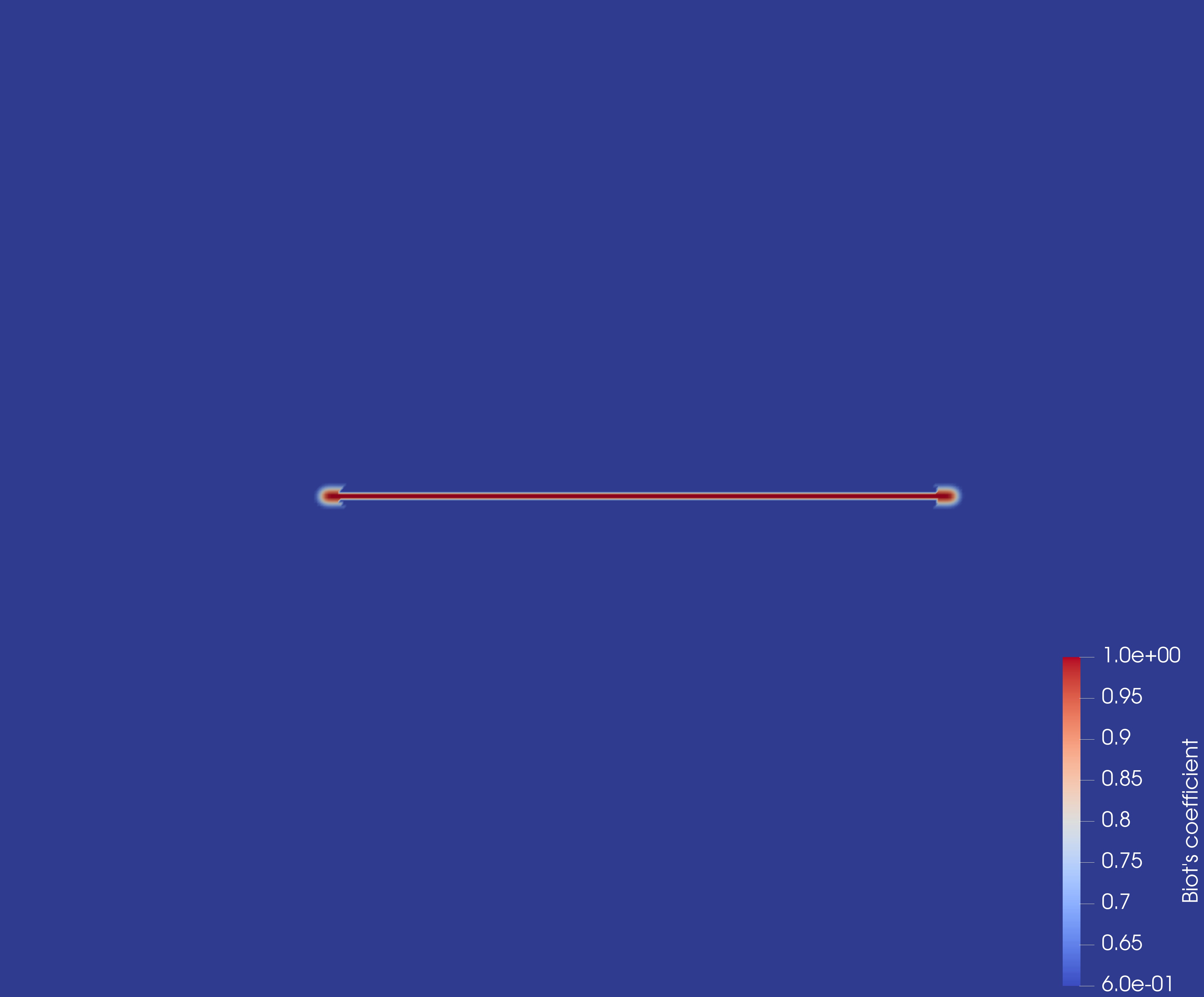} \label{fig:alpha_int_2}}
    \caption{Phase-field (a), pressure (b), width (b) and Biot's coefficient (d) profiles for the hydraulic fracturing interacting with the natural fracture at $t=10$ s. ($\beta=90^\circ$, $\Gc^\mrm{int}/\Gc=0.5$)  }
    \label{fig:results_90_05}
\end{figure}

We first investigate the effect of natural fracture with an inclination angle of $\beta=90^{\circ}$. 
Figs.~\ref{fig:results_90_02} and \ref{fig:results_90_05} show that the hydraulic fracture branches along the weaker interface ($\Gc^{\mrm{int}}=0.2\Gc$), but bypasses the relatively stronger interface ($\Gc^{\mrm{int}}=0.5\Gc$). 
As can be seen, the cementation state of the interface (represented by different interface toughness) impacts the final hydraulic channel network and the resulting pressure field. 

Furthermore, Figs. \ref{fig:alpha_int_1} and \ref{fig:alpha_int_2} show that Biot's coefficient transitions smoothly from 0.6 to 1 around the crack tip while it changes sharply in the vicinity of induced hydraulic fracture.
As Biot's coefficient (Eq.~\eqref{eq:alpha}) depends both on the phase-field variable and on the strain decomposition strategy, pressurizing a hydraulic fracture will induce compressive strains in the nearby reservoir.
Consequently, we have $H (\Tr{\bm{\eps}} ) = 0$ and obtain $\alpha = \alpha_\mrm{m}$. 
And the unilateral contact effect of the fracture on Biot's coefficient can be well captured by the proposed model. 

\begin{figure}[htpb!]
    \centering
    \subfigure[$\beta=90^\circ$]{\includegraphics[scale=0.5]{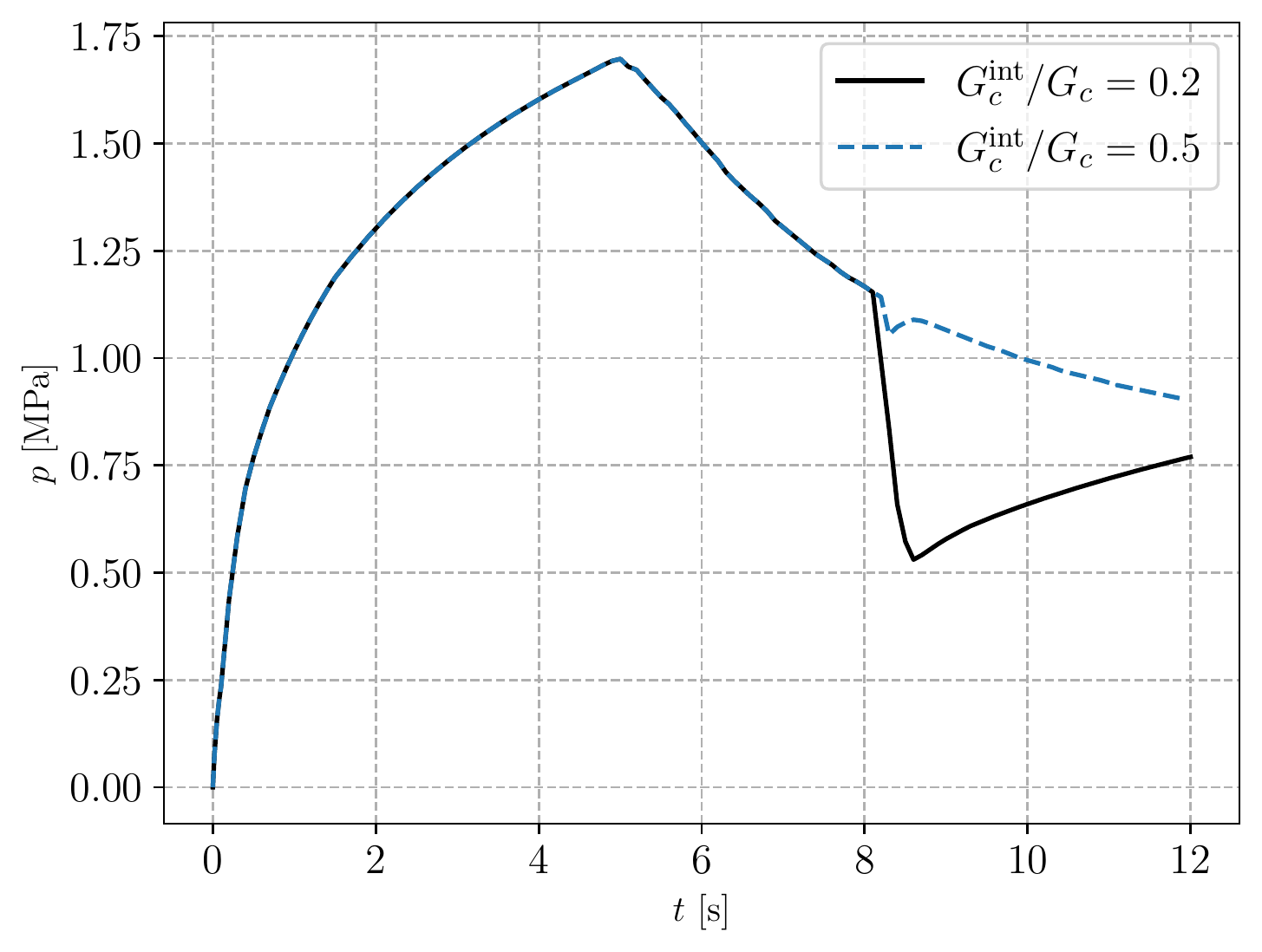} \label{fig:pressure_90}}
    \subfigure[$\beta=165^\circ$]{\includegraphics[scale=0.5]{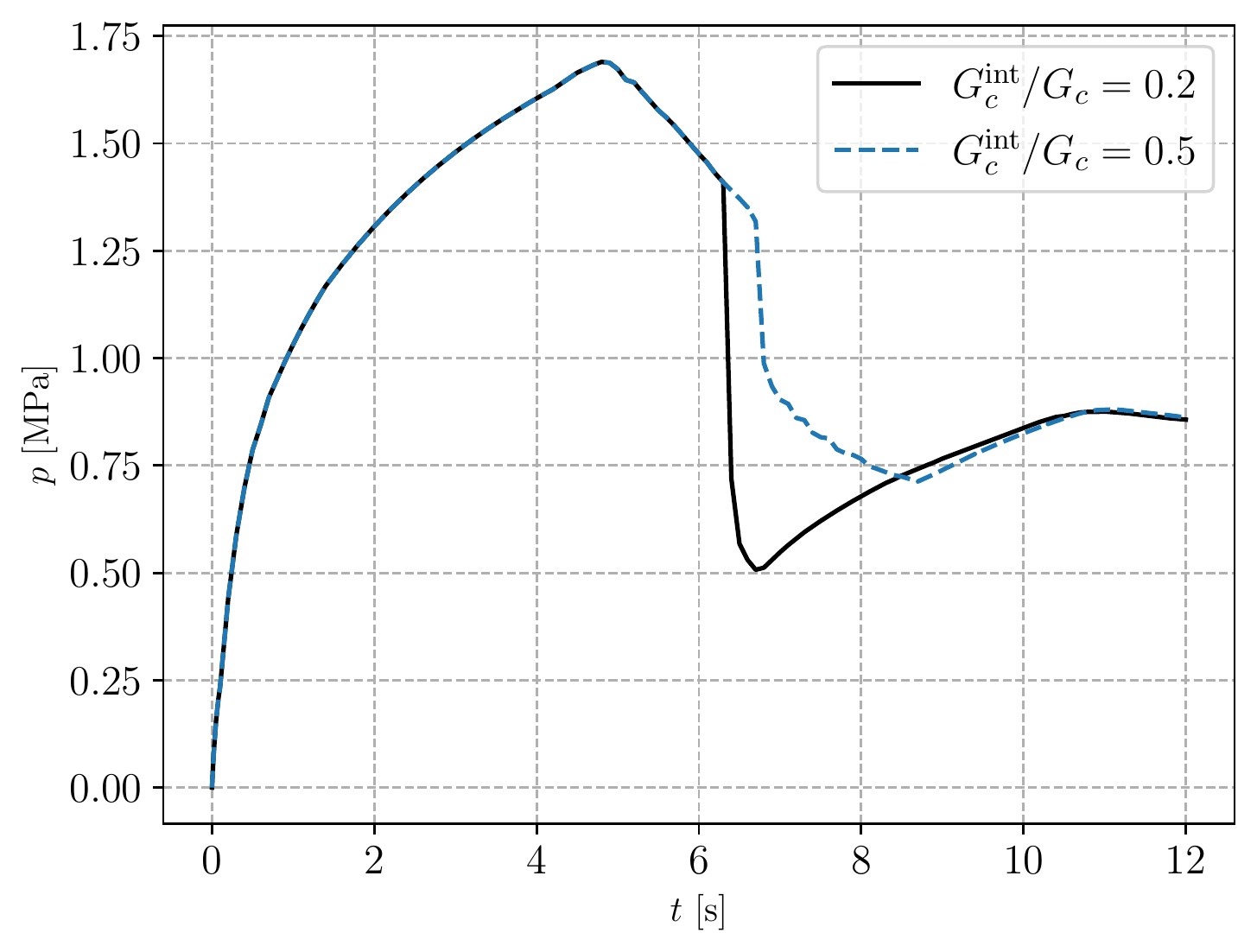} \label{fig:pressure_165}}
    \caption{Pressure responses at the injection point for natural fractures with  inclination angles being $\beta=90^\circ \mrm{\, and \,} 165^\circ$.}
\end{figure}

Fig. \ref{fig:pressure_90} shows the pressure responses for two interface toughness with $\beta=90^\circ$.
These two responses are identical in the early time until the hydraulic fracture starts interacting with the natural fracture.
For $\Gc^{\mrm{int}}/\Gc=0.2$, we see a sudden drop when the hydraulic fracture branches along the interface. 
Once the branched fractures reach the reservoir media, the pressure starts to increase. 
In contrast, for $\Gc^{\mrm{int}}/\Gc=0.5$, we see only a slight pressure drop when the hydraulic fracture bypasses the weak interface but then it follows back the original pressure decline trend.

\begin{figure}[htpb!]
    \centering
    \subfigure[Phase field]{\includegraphics[width=7cm]{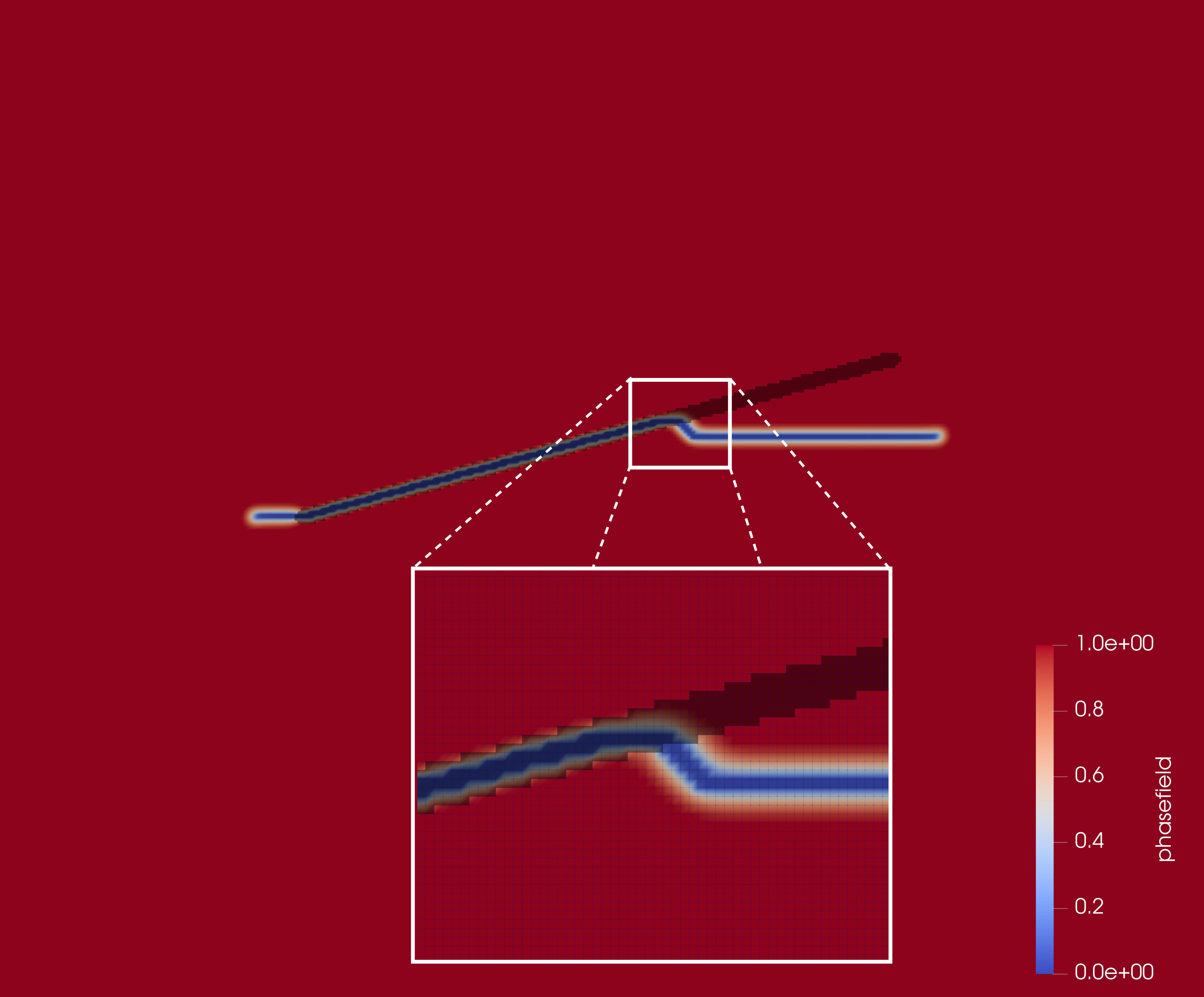} \label{fig:pf_combined_int_3}}
    \subfigure[Pressure]{\includegraphics[width=7cm]{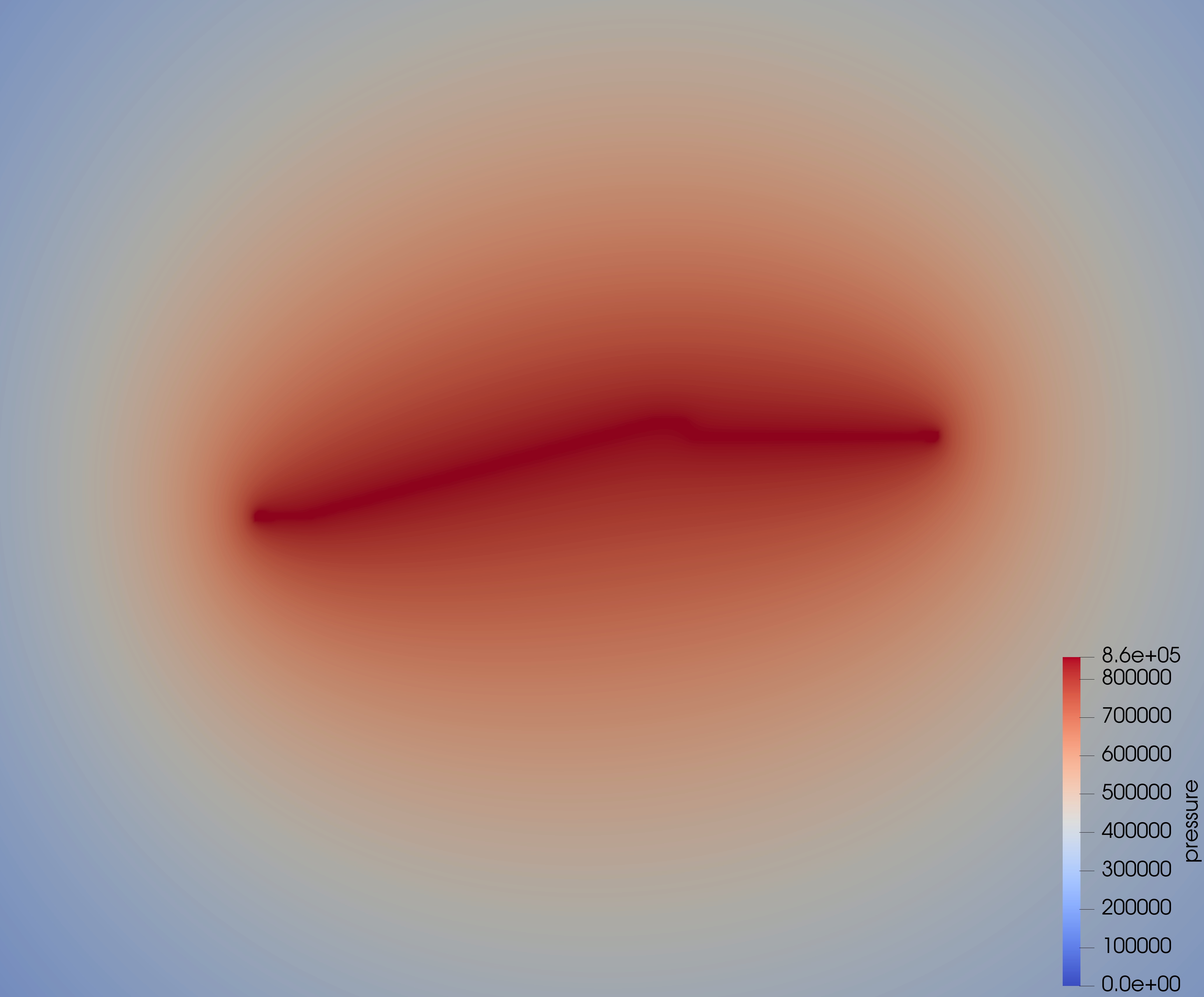} \label{fig:pressure_int_3}}
    \subfigure[Width (Interpolated onto node)]{\includegraphics[width=7cm]{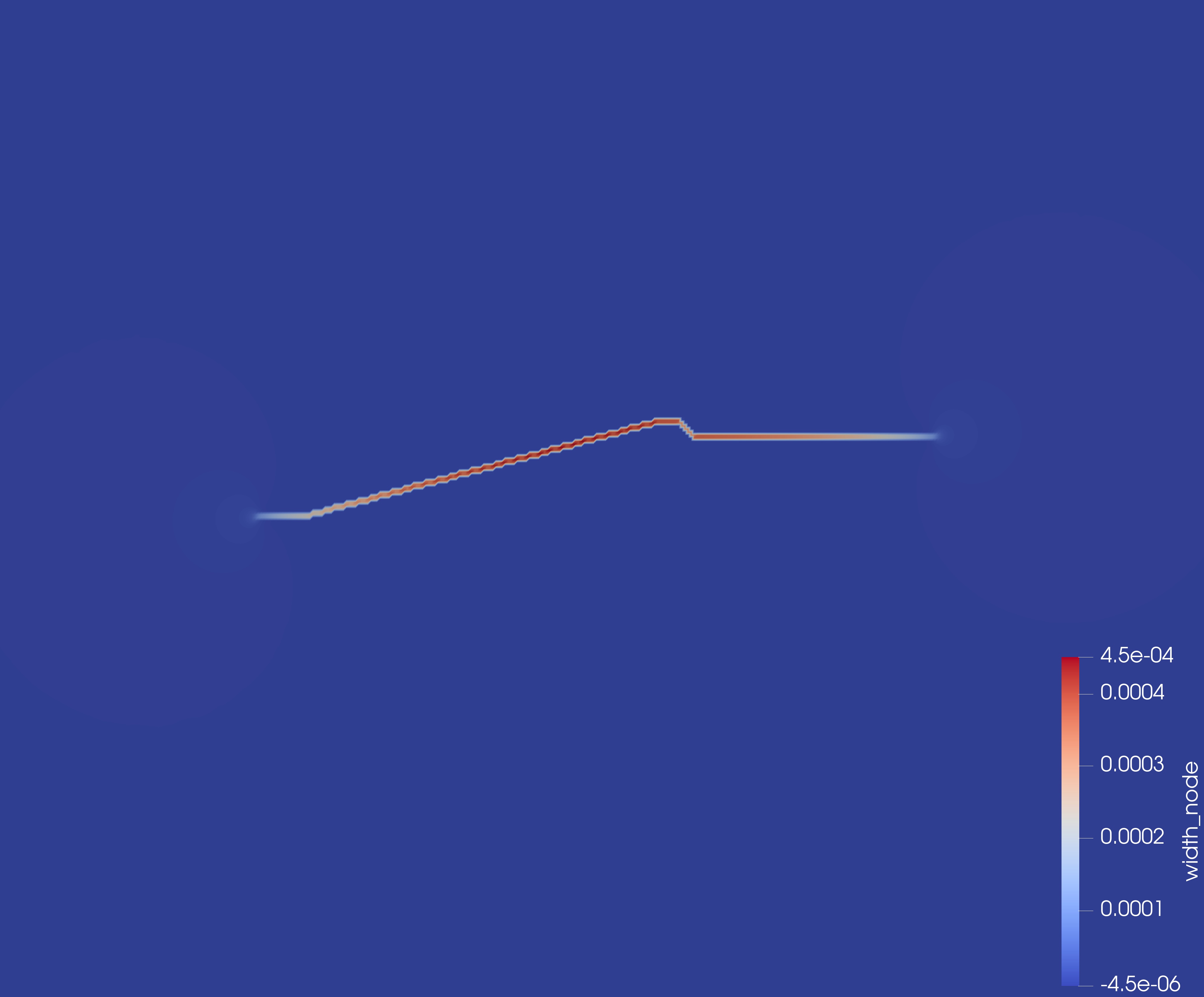} \label{fig:width_int_3}}
    \subfigure[Boit's coefficient (Calculated by Eq. \eqref{eq:alpha} )]{\includegraphics[width=7cm]{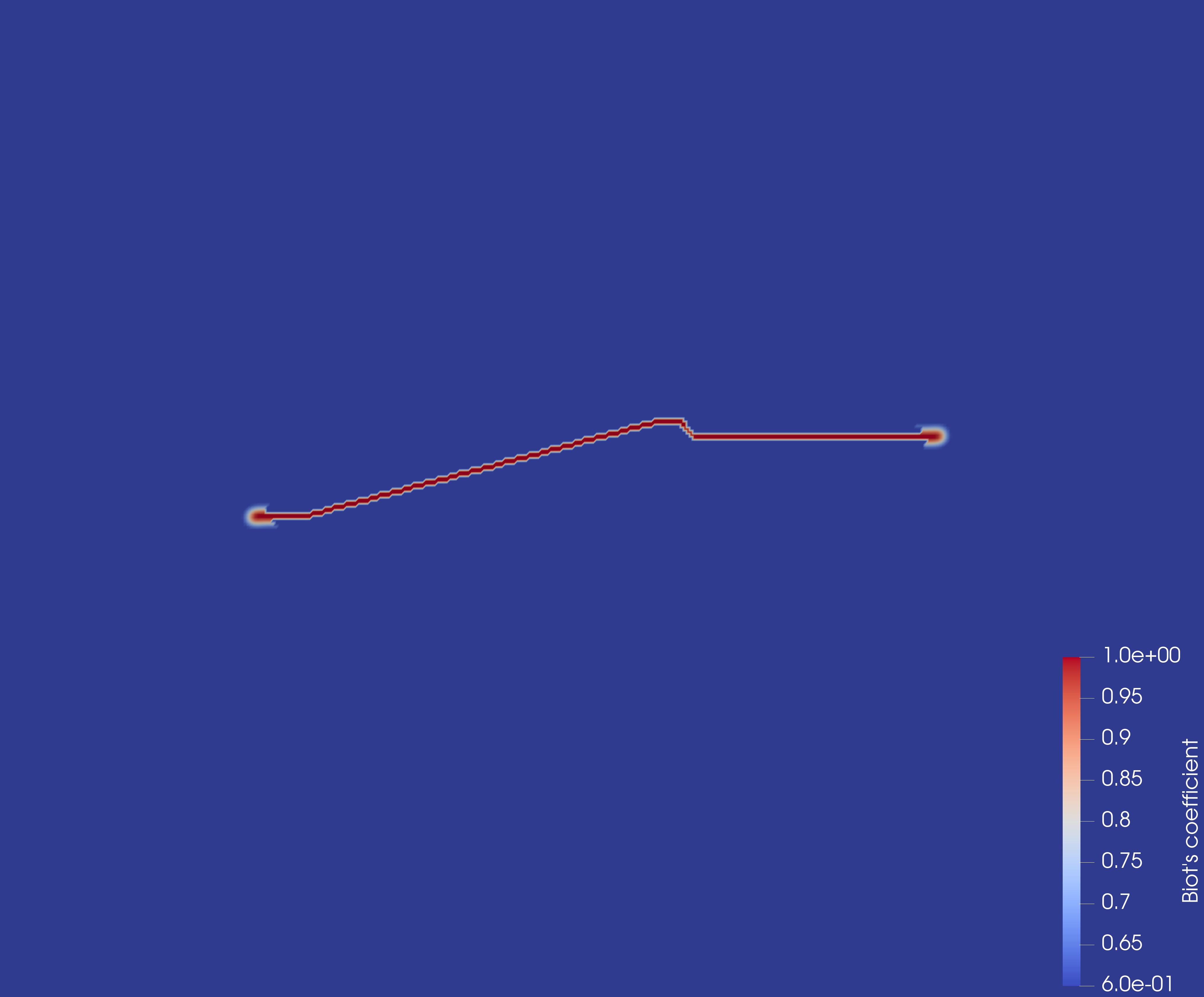} \label{fig:alpha_int_3}}
    \caption{Phase-field (a), pressure (b), width (b) and Biot's coefficient (d) profiles for the hydraulic fracturing interacting with the natural fracture at $t=10$ s. ($\beta=165^\circ$, $\Gc^\mrm{int}/\Gc=0.2$)  }
\end{figure}

\begin{figure}[htpb!]
    \centering
    \subfigure[Phase field]{\includegraphics[width=7cm]{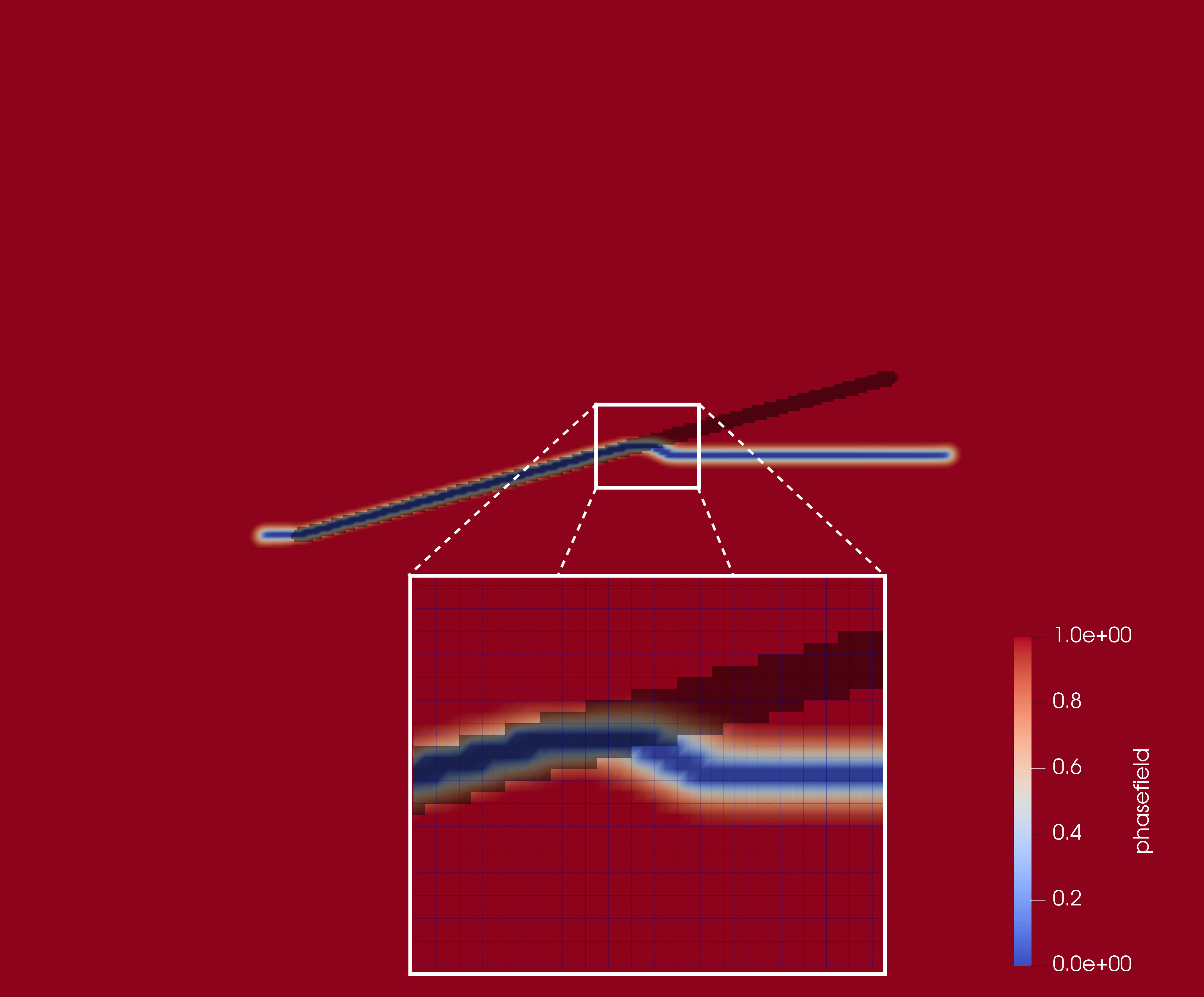} \label{fig:pf_combined_int_4}}
    \subfigure[Pressure]{\includegraphics[width=7cm]{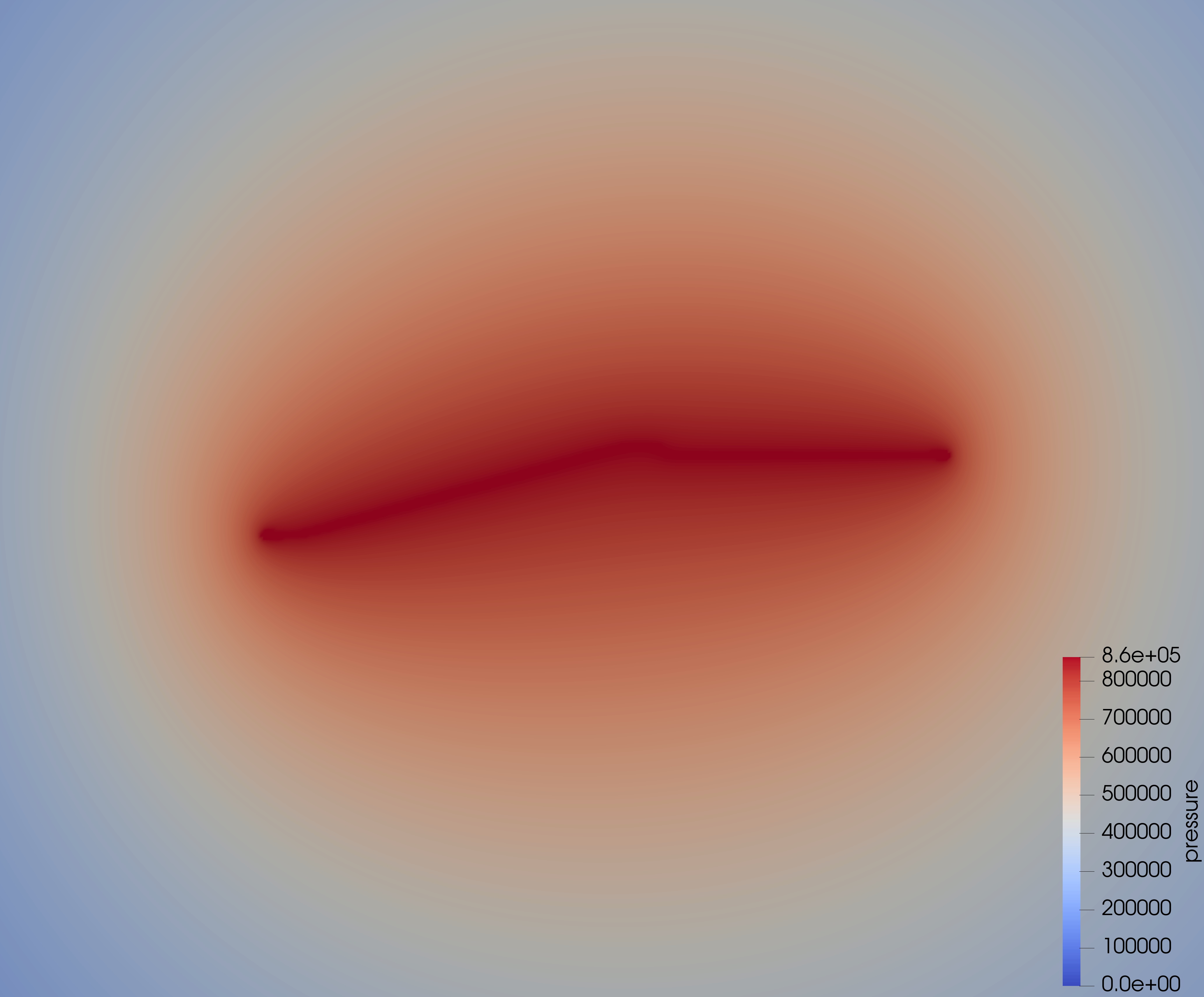} \label{fig:pressure_int_4}}
    \subfigure[Width (Interpolated onto node)]{\includegraphics[width=7cm]{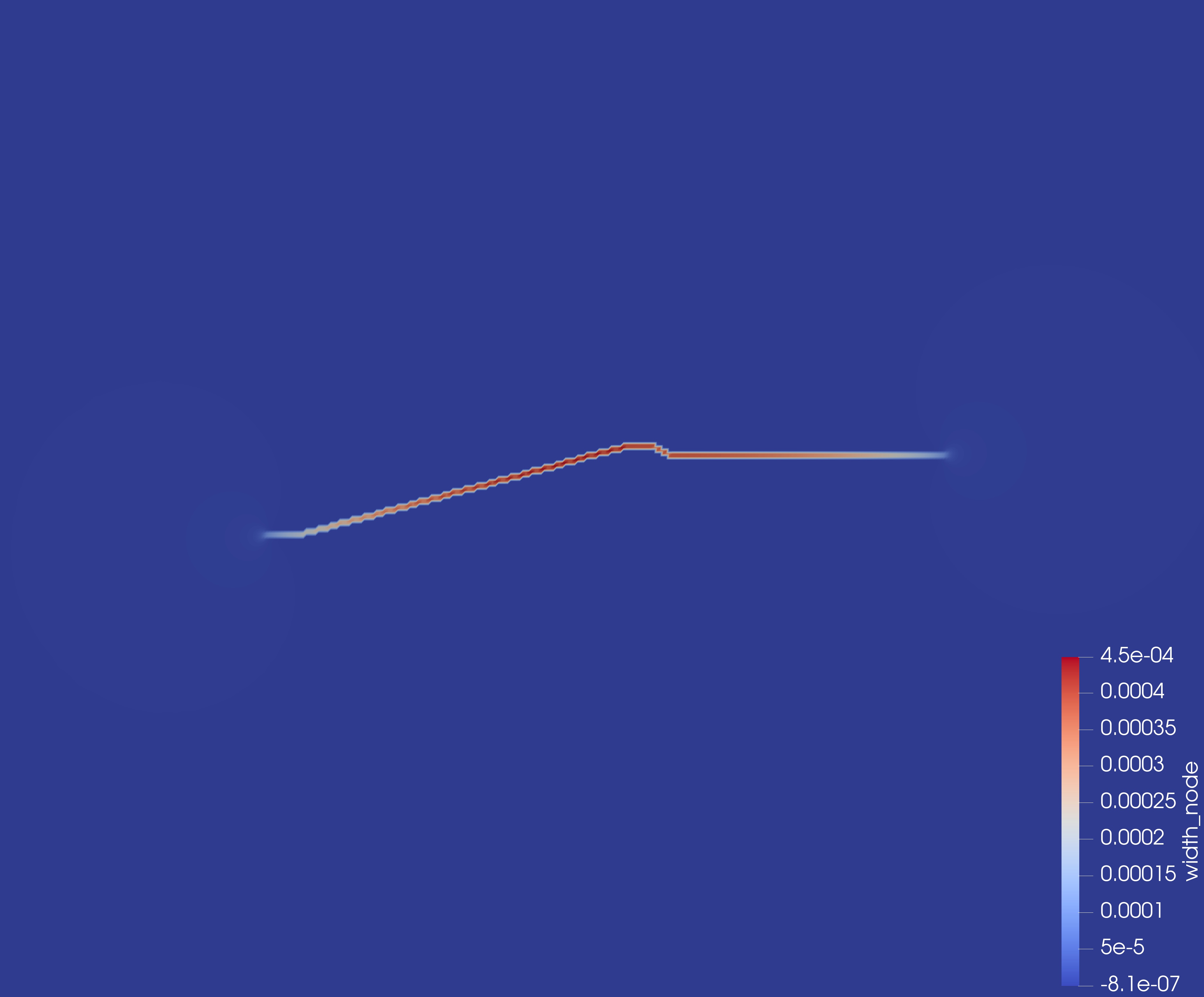} \label{fig:width_int_4}}
    \subfigure[Boit's coefficient (Calculated by Eq. \eqref{eq:alpha} )]{\includegraphics[width=7cm]{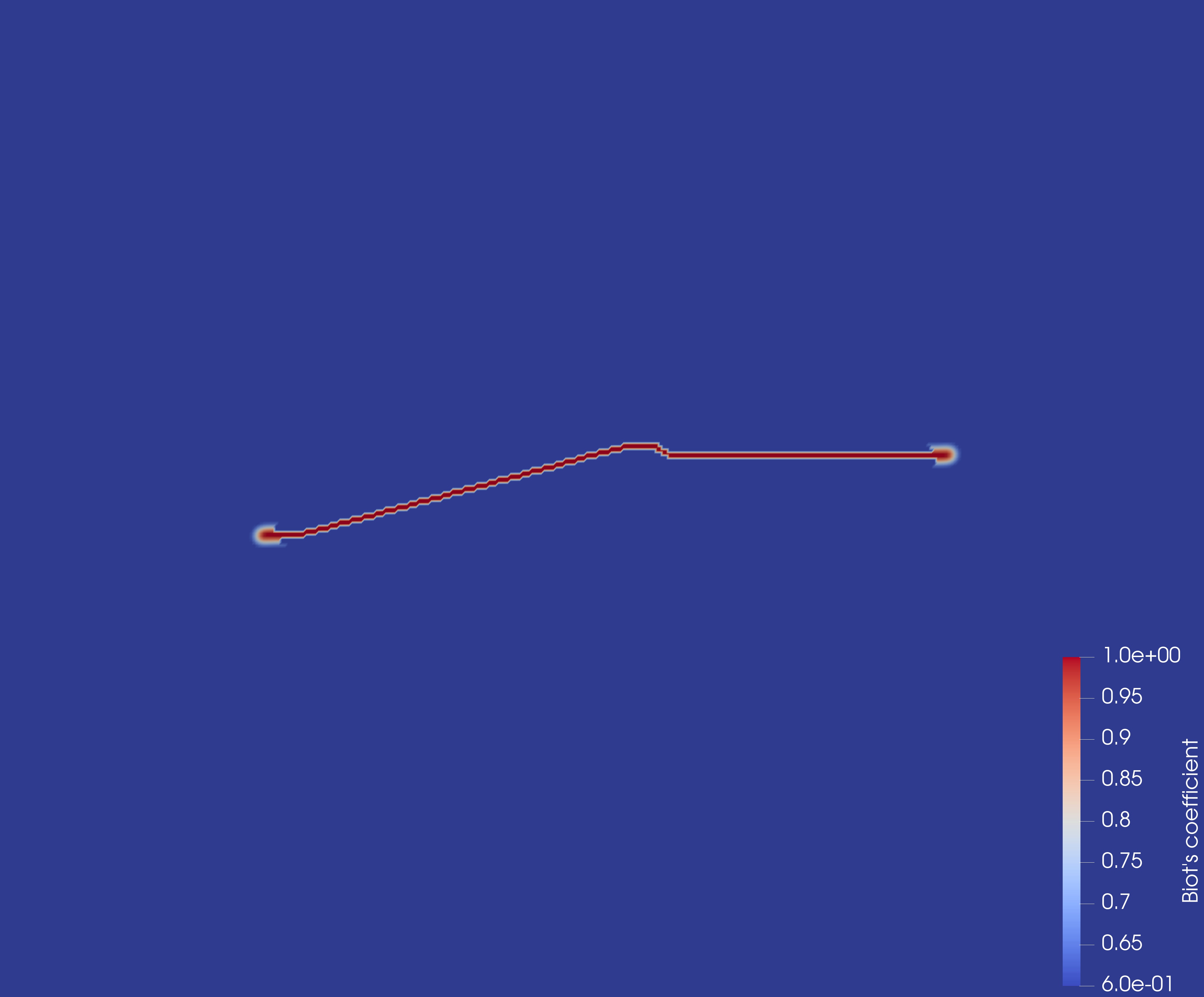} \label{fig:alpha_int_4}}
    \caption{Phase-field (a), pressure (b), width (b) and Biot's coefficient (d) profiles for the hydraulic fracturing interacting with the natural fracture at $t=10$ s. ($\beta=165^\circ$, $\Gc^\mrm{int}/\Gc=0.5$)  }
\end{figure}

With $\beta=165^\circ$, the natural fracture attracts the hydraulic fracture as shown in Figs. \ref{fig:pf_combined_int_3} and \ref{fig:pf_combined_int_4}. The smaller $\Gc^{\mrm{int}}/\Gc$, the easier and earlier the hydraulic fracture is attracted into the natural fracture (\ref{fig:pressure_165}). 
After the hydraulic fracture propagates through the natural fracture, the pressure increases until the propagation resumes in the reservoir media (see Figs.~\ref{fig:pf_combined_int_3}~and~\ref{fig:pf_combined_int_4}).

The distributions of Biot's coefficient also show a sharp change in the vicinity of developed hydraulic fracture and a smooth transition in the crack tip. 
The examples demonstrate that our approaches for width Eq.~\eqref{eq:width} and crack normal vector Eq.~\eqref{eq:eigen_strain} computations are well applicable to complex hydraulic fractures. 

\section{Conclusions}
In this paper, we proposed a hydro-mechanical phase-field model with diffused poroelasticity derived from micromechanical analysis. 
Our proposed Biot's parameters depend not only on the phase-field variable but also on the type of energy decomposition, and the poroelastic strain energy is degraded through Biot's parameters rather than degrading the fluid pressure terms. 
The theoretical analysis shows that our proposed model automatically ensures the stress continuity on the crack interface (i.e., $\mbfs{\sigma}\cdot \mbfs{n}_\Gamma = p\mbfs{n}_\Gamma$ with $v=0$). 
Our proposed model was compared against the three commonly-used phase-field models for poroelastic media, and it can more accurately recover crack opening displacements irrespective of the pressure profile or Biot's coefficient.
The proposed model was also verified against the KGD fracture solution in the toughness-dominated regime. 
Furthermore, we simulated hydraulic fracture interactions with a natural fracture to demonstrate its capability for dealing with complex hydraulic fracture behaviours.

The proposed model is implemented in an open source code, OpenGeoSys~\cite{ogs:6.4.3}.
As future study, the model is to be extended for thermo-hydro-mechanical and phase-field coupling problems.
On the other hand, to deal with the loss of local mass conservation of fluid for the traditional continuous Galerkin finite element method, the mixed finite element  method~\cite{feng2021phase} or enriched Galerkin finite element method \cite{lee2016locally} may be needed in the future.

\section*{CRediT author statement}
All authors have contributed equally.

\section*{Acknowledgments}
The authors gratefully acknowledge the funding provided by the Deutsche Forschungsgemeinschaft (DFG, German Research Foundation), Germany through the HIGHER project (Grant No. PA 3451/1-1). TY's contribution is funded by the China and Germany Postdoctoral Exchange Program (Grant No. ZD202137) and the National Natural Science Foundation of China (Grant No. 12202137). The first author (TY) would like to express his gratitude to Dr. Mostafa Mollaali for his help in visualization and Dr. Yijun Chen for the informative discussions. 

\section*{Conflict of interest}
The authors declare that they have no known competing financial interests or personal relationships that could have appeared to influence the work reported in this paper.

\bibliographystyle{apalike}
\bibliography{reference} 
\end{document}